\documentclass[12pt,a4paper]{article}

\usepackage[margin=25.4mm]{geometry}

\usepackage[utf8]{inputenc}
\usepackage[T1]{fontenc}
\usepackage[]{microtype}
\usepackage{relsize,comment}
\usepackage{todonotes}
\usepackage[colorlinks=false,citecolor=black,hidelinks]{hyperref}

\usetikzlibrary{arrows}

\usepackage{stmaryrd}
\SetSymbolFont{stmry}{bold}{U}{stmry}{m}{n} 

\usepackage{xcolor,scrtime,graphicx}
\usepackage{amsmath,amsfonts,amsmath,eucal,mathrsfs,amssymb}

\usepackage{enumerate,todonotes}
\usepackage{bef_alex}

\def\Bigset#1#2{\Big\{\; #1 \; \Big| \; #2 \; \Big\} }

\numberwithin{equation}{section}
\numberwithin{figure}{section}

\usepackage{bm} 
\renewcommand{\FG}{\bm}
\renewcommand{\mafo}{\mathrm}
\renewcommand{\ti}{{\times}}

\usepackage[final,notref,notcite]{showkeys} 

\makeatletter
\renewcommand*\env@cases[1][1.2]{%
  \let\@ifnextchar\new@ifnextchar
  \left\lbrace
  \def\arraystretch{#1}%
  \array{@{\,}c@{\ }l@{}}%
}
\makeatother

\usepackage[normalem]{ulem}
\usepackage{tocloft}
\newcommand{\MMM}{\mathcal M}
\newcommand{\ArgM}{\mathsf M}

\newcommand{\sfH}{\mathsf H}
\newcommand{\sfW}{\mathsf W}
\newcommand{\HK}{\mathsf H\!\!\mathsf K}
\def\div{\mathop{\mafo{div}}}
\newcommand\trace{\mathop{\mafo{tr}}}
\newcommand{\oti}{{\otimes}}

\newcommand{\mfh}{\mathfrak h}

\newcommand{\msd}{\mathsf{d}}
\newcommand{\calET}{\calE\kern-3pt \calT}

\newcommand{\DDD}{\mathfrak{D}}

\newcommand{\AAA}{}
\newcommand{\EEE}{\color{black}}
\newcommand{\GGG}{}
\newcommand{\MAT}{}
\let\psxi\xi
\let\one\tau

\newcommand{\bfarctan}{\bm{\mathrm{arctan}}}

\newcommand{\TRAGRO}{{transport-growth}}
\newcommand{\TTRAGRO}{{Transport-growth}}

\newtheorem{problem}[theorem]{Problem}

\newcommand{\sfd}{\mathsf d}
\newcommand{\sfx}{\mathsf x}
\newcommand{\sfr}{\mathsf r}

\newcommand{\dist}{\mathrm{dist}}
\newcommand{\supp}{\operatorname{supp}}
\newcommand{\cP}{\calP}

\newcommand{\OOmega}{{\R^d}}
\newcommand{\Leb}[1]{\calL^{#1}}
\newcommand{\ggamma}{\boldsymbol \gamma}

\newcommand{\STAR}{{\bm*}}
\newcommand{\mfo}{\mathfrak o}
\newcommand{\mfC}{\mathfrak C}
\newcommand{\mfz}{\mathfrak z}

\newcommand{\bftan}{\mbox{\boldmath$\tan$}}
\newcommand{\bfsin}{\mbox{\boldmath$\sin$}}

\newcommand{\Mcost}[3]{\mathcal C(#1,#2;#3)}

\newcommand{\res}{\mathop{\hbox{\vrule height 7pt width .5pt depth 0pt
\vrule height .5pt width 6pt depth 0pt}}\nolimits}

\renewcommand{\calE}{\mathscr E}
\newcommand{\Ell}{\mathrm L}
\newcommand{\GRAD}{\bm g}
\newcommand{\LLL}{\mathbb L}
\newcommand{\RRR}{\mathsf R}
\newcommand{\geo}[3]{\mathrm{geo}_{#1}\big(#2,#3\big)}

\newcommand{\tcos}[2]{\cos_{#1}\!{(#2)}}
\newcommand{\tcosq}[2]{\cos^2_{#1}\!{\left(#2\right)}}

\newcommand{\Gtrafo}{\mathsf G}
\newcommand{\chGtrafo}{\check{\mathsf G}}
\newcommand{\IGtrafo}{{\mathsf G}^{-1}}
\newcommand{\cIGtrafo}{\check{\mathsf G}^{-1}}
\newcommand{\HopfLax}[2]{\mathscr P_{\kern-2pt #1}\kern1pt #2}
\newcommand{\RopfLax}[2]{\mathscr R_{\kern-2pt #1}\kern1pt #2}
\newcommand{\Int}{\operatorname{int}}
\newcommand{\Cl}{\operatorname{cl}}
\newcommand{\Ext}[2]{\operatorname{ext}_{#1}(#2)}
\newcommand{\Bdry}[2]{\operatorname{bdry}_{#1}(#2)}
\newcommand{\pit}{{\pi/2}}
\newcommand{\sfA}{{\mathsf A}}
\newcommand{\sfB}{{\mathsf B}}
\newcommand{\bigup}{\vphantom{\big|}}

\newcommand{\MR}[1]{\relax}

\begin{document}
\title{Fine properties of geodesics and geodesic $\lambda$-convexity \\
for the Hellinger--Kantorovich distance} 

\author{Matthias Liero\thanks{Weierstraß-Institut f\"ur Angewandte Analysis 
               und Stochastik, Berlin, Germany.
               ORCID: \url{https://orcid.org/ 0000-0002-0963-2915}}
, 
Alexander Mielke$^{*,}$\thanks{Humboldt Universität zu Berlin,
           Germany. ORCID:~\url{https://orcid.org/0000-0002-4583-3888}}, 
and 
Giuseppe Savar\'e\thanks{Department of Decision Sciences and BIDSA, 
            Bocconi University, Via Roentgen 1, 20136 Milan, Italy. 
            Email: \texttt{giuseppe.savare@unibocconi.it},
            ORCID: \url{https://orcid.org/0000-0002-0104-4158}}
   } 

\date{August 29, 2022; revised September 22, 2023} 

\maketitle
\vspace{-2em}

\begin{abstract}\footnotesize
  We study the fine regularity properties of optimal potentials for the dual
  formulation of the Hellinger--Kantorovich problem ($\HK$), providing sufficient
  conditions for the solvability of the primal Monge formulation.  We also
  establish new regularity properties for the solution of the Hamilton--Jacobi
  equation arising in the dual dynamic formulation of $\HK$, which are
  sufficiently strong to construct a characteristic \TRAGRO\ flow
  driving the geodesic interpolation between two arbitrary positive measures.
  
  These results are applied to study relevant geometric properties of $\HK$
  geodesics and to derive the convex behaviour of their Lebesgue density along
  the transport flow.
  
  Finally, exact conditions for functionals defined on the space of measures
  are derived that guarantee the geodesic $\lambda$-convexity with respect to
  the Hellinger--Kantorovich distance.  Examples of geodesically convex
  functionals are provided.
\end{abstract}

\vspace{-2.5em}
{\setlength{\cftbeforesecskip}{3pt}\footnotesize \renewcommand{\contentsname}{}\tableofcontents}

\section{Introduction}
\label{s:Intro}

In \cite{LiMiSa16OTCR, LiMiSa18OETP} the Hellinger--Kantorovich distance (in
\cite{KoMoVo16NOTD, CPSV18IDOT, CPSV18UOTG} it is also called
Wasserstein--Fisher--Rao distance or Kantorovich--Fisher--Rao distance in
\cite{GalMon17JKOS}) was introduced to describe the interaction between optimal
transport and optimal creation and destruction of mass in a convex domain of
$\R^d$.  Here we further investigate the structure of (minimal) geodesics, and
we fully analyze the question of geodesic $\lambda$-convexity of integral
functionals with respect to this distance.

The Hellinger--Kantorovich distance can be considered as a combination, more
precisely the inf-convolution, of the Hellinger--Kakutani distance on the set
of all measures (cf.\ e.g.\,\cite{Schi97GAFV}) and the $\rmL^2$
Kantorovich--Wasserstein distance, which is well-known from the theory of
optimal transport, see e.g.\ \cite{AmGiSa08GFMS,Vill09OTON}. Throughout this
text, we denote by $\MMM(\OOmega)$ all nonnegative and finite
Borel measures endowed with the weak topology induced by the canonical duality
with the continuous functions $\rmC_0(\OOmega)$ decaying at infinity. 
While the $\rmL^2$ Kantorovich--Wasserstein distance $\sfW(\mu_0,\mu_1)$ of
measures $\mu_0$, $\mu_1 \in \MMM(\OOmega)$ requires $\mu_0$ and $\mu_1$ to have
the same mass to be finite, the Hellinger--Kakutani distance, which is defined
via
\[
\sfH(\mu_0,\mu_1)^2 = \int_\OOmega  \big(\sqrt{\theta_0} - \sqrt{\theta_1}\big)^2 \dd
(\mu_0{+}\mu_1), \text{ where } \theta_j = \frac{\rmd
\mu_j}{\rmd(\mu_0{+}\mu_1)},
\]
has the upper bound $\sfH(\mu_0,\mu_1)\leq \mu_0(\OOmega)+ \mu_1(\OOmega)$, with
equality if $\mu_0$ and $\mu_1$ are mutually singular.

As a generalization of the dynamical formulation of the Kantorovich--Wasserstein distance
(see \cite{BenBre00CFMS}), the Hellinger--Kantorovich distance $\HK_{\alpha,\beta}$ can 
be defined in a dynamic way via
\begin{align}
 \label{eq:HKdynamical}
 \HK_{\alpha,\beta}(\mu_0,\mu_1)^2 = \inf \bigg\{ &\int_{t=0}^1
  \int_\OOmega  \big(\alpha |\Upsilon(t,x)|^2 {+}\beta\xi(t,x)^2\big)
                  \dd\mu_t(x) \dd t  \\
 \nonumber
              &  \bigg| \ \mu\in \rmC\big([0,1];\MMM(\OOmega)\big),  \
                \mu_{t=0}=\mu_0,\ \mu_{t=1}=\mu_1,
                \text{ (gCE) holds} \bigg\},
\end{align}
where   $\Upsilon:(0,1)\times \OOmega\to \R^d$ and $\xi:(0,1)\times \OOmega\to
\R$ are Borel maps characterizing  
the generalized continuity equation 
\[
\text{(gCE) } \qquad   \frac{\partial}{\partial t} \mu + \alpha\,\div \big( \mu
\Upsilon \big) =\beta\, \xi \mu,
\]
  formulated in a distributional sense.
The parameters $\alpha>0$ and $\beta>0$ allow us to control the relative
strength of the Kantorovich--Wasserstein part and the
Hellinger--Kakutani part, i.e.\ $\HK_{\alpha,\beta}$ is the
inf-con\-volution of $ \HK_{\alpha,0}= \frac1{\sqrt \alpha} \sfW$ and
$ \HK_{0,\beta}= \frac1{\sqrt \beta}\sfH$,  see
\cite[Rem.\,8.19]{LiMiSa18OETP}. 
Subsequently, we will restrict to the standard case $\alpha=1$ and
$\beta=4$, since the general case can easily be obtained by scaling
the underlying space $\OOmega$. We will shortly write $\HK$ instead of
$\HK_{1,4}$.

It is a remarkable fact, deeply investigated in \cite{LiMiSa18OETP}, that the
$\HK$ distance has many interesting equivalent characterizations, which
highlight its geometric and variational character.  A first one arises from the
dual dynamic counterpart of \eqref{eq:HKdynamical} in terms of subsolutions of
a suitable Hamilton--Jacobi equation:
\begin{equation}
  \label{eq:16intro}
  \begin{aligned}
    \frac 12\HK^2(\mu_0,\mu_1)= \sup\bigg\{\int_\OOmega
    \xi(\tau,\cdot)\dd\mu_1 &-\int_\OOmega\xi(0,\cdot)\dd\mu_0\,\Big|\, \xi\in
    \rmC^\infty_\rmc([0,1]\times \R^d),
    \\
    & \frac\partial{\partial t}\xi+
    \frac12 |\nabla \xi|^2+2 \xi^2\le0 \quad\text{in }[0,1]\times\R^d\bigg\}.
  \end{aligned}
\end{equation}
By expressing solutions of \eqref{eq:16intro}
in terms of a new formula of Hopf--Lax type,
one can write a static duality representation
\let\one1%
\begin{align}
  \nonumber
  \HK^2(\mu_0,\mu_\one)&=
                         \sup\Big\{\int_\OOmega (1-\ee^{-2\varphi_\one})
                         \dd\mu_\one
                         -\int_\OOmega (\ee^{2\varphi_0}-1)
                    \dd\mu_0\,\Big|
\,\\&\hspace{5em}\varphi_0,\varphi_\one\in\rmC_b(\OOmega),~
  \varphi_\one(x_\one){-}\varphi_0(x_0)\leq \Ell_1(x_\one{-}x_0)\Big\}
\label{eq:HK.dual.intro}
\end{align}
associated with the convex cost function
$\Ell_1(z):=\frac12\log(1+\tan^2(|z|))$
which forces $|z|<\pi/2$. 
\GGG Notice that it is possible to write \eqref{eq:HK.dual.intro}
in a symmetric form with respect to $\varphi_0,\varphi_1$ just by
changing the sign of $\varphi_1$.

It is remarkable that \eqref{eq:HK.dual.intro} 
can be interpreted as the dual problem
of the static
Logarithmic Entropy Transport (LET) variational formulation
of $\HK$.
By introducing the logarithmic entropy density
$F:\left[0,\infty\right[\to[0,\infty[$ via
\begin{equation}
\label{eq:LED-intro}
F(s) := s\log s - s + 1\quad\text{for $s>0$}\quad \text{and} \quad  F(0):=1,
\end{equation}
we get 
\begin{equation}
  \label{eq:LET-intro}
  \HK^2(\mu_0,\mu_\one)=
  \min\Big\{
  \int_\OOmega F(\sigma_0)\dd\mu_0+ \int_\OOmega F(\sigma_1)\dd\mu_1+
   \iint_{\OOmega\times\OOmega}2\Ell_1(x_0{-}x_1)\dd\bfeta 
   \Big\}
\end{equation}
where the minimum is \AAA taken over all \GGG positive finite Borel measures
$\bfeta$ in $\R^d\times \R^d$ whose marginals
$(\pi_i)_\sharp \bfeta=\sigma_i\mu_i$ are absolutely continuous with respect to
$\mu_i$.  \EEE

The subdifferential 
\[ 
\rmD \Ell_1(z)= \pl\Ell_1(z)=
\bftan(z) := \tan\big(|z|)\,\frac{z}{|z|}
\]
and its inverse $w \mapsto \bfarctan(w)$ will play an important role.  We
continue to use bold function names for vector-valued functions constructed
from real-valued ones as follows:
\begin{equation}
  \text{for a map $f:\R\to \R$ with
$f(0)=0$ we set 
$\bm f : \R^d \to \R^d$ via $\bm f(x):=f(|x|)
\frac{x}{|x|}$}.\label{eq:5}
\end{equation}
A \GGG fourth \EEE crucial formula, which we will extensively study in the
present paper, is related to the primal Monge formulation of Optimal Transport,
and clarifies the two main components of $\HK$ arising from transport and \AAA
growth or decay \EEE effects.  Its main ingredient is the notion of \TRAGRO\
pair $(\bfT,q):\R^d\to \R^d\times [0,\infty)$ acting on measures
$\mu\in \MMM(\R^d)$ as
\begin{equation}
  \label{eq:171}
  (\bfT,q)_\star \mu:=
  \bfT_\sharp(q^2\cdot \mu),\quad
  \big((\bfT,q)_\star \mu\big)(A):=
  \int_{\bfT^{-1}(A)}q^2\dd\mu\quad
  \text{for every Borel set }A\subset \R^d.
\end{equation}
The Monge formulation of $\HK$ then looks for
the optimal pair $(\bfT,q)$
among the ones transforming $\mu_0$ into $\mu_1$
by $(\bfT,q)_\star\mu_0=\mu_1$ which minimizes the
conical cost
\begin{equation}
  \label{eq:172intro}
  \Mcost \bfT q{\mu_0}:=\int_{\OOmega} 
    \Big(1+q^2(x)-2q(x)
    \cos_\pit\big(|\bfT(x){-}x|\big)\Big)\dd\mu_0(x),
\end{equation}
where $\cos_\pit(r):=\cos\big( \min\{r,\pit\}\big)$.  As for the usual
Monge formulation of 
optimal transport, the existence of an optimal \TRAGRO\ pair
$(\bfT,q)$ minimizing \eqref{eq:172intro} requires more restrictive properties
on $\mu_0,\mu_1$ which we will carefully study.  It is worth noticing that the
integrand in \eqref{eq:172intro} has a relevant geometric interpretation as the
\GGG square distance $\sfd^2_{\pi,\mfC}$, 
where $\sfd_{\pi,\mfC}$ is the distance \EEE
on the cone space $\mfC$ over $\R^d$ (cf.\ \eqref{eq:7}) between the points $[x,1]$ and
$[\bfT(x),q(x)]$ and suggests that $\HK$ induces a distance in $\MMM(\R^d)$
which plays a similar role than the $\rmL^2$ Kantorovich--Rubinstein--Wasserstein
distance in $\mathcal P_2(\R^d)$. The dynamic formulation
\eqref{eq:HKdynamical}, moreover, suggests that its minimizers
$(\mu_t)_{t\in [0,1]}$ should provide minimal geodesics in $(\MMM(\R^d),\HK)$
which behave like \TRAGRO\ interpolations between $\mu_0$ and
$\mu_1$.

Inspired by the celebrated paper \cite{Mcca97CPIG}, we want to study the
structure of such minimizers and to characterize integral functionals which are
convex along such kind of interpolations.

\subsection{Improved regularity of potentials and geodesics}
\label{su:IntroGeod}

In the first part of the paper we will exploit the equivalent formulations of
$\HK$ in order to obtain new information on the regularity and on the fine
structure of the solutions to \eqref{eq:HK.dual.intro}, \eqref{eq:16intro}, and
\eqref{eq:172intro}.

More precisely, we will initially prove in Section \ref{se:2nd.Optim} that the
optimal $\HK$ potential $\varphi_0$ is locally semi-convex outside a closed
$(d{-}1)$-rectifiable set, so that when $\mu_0\ll\Leb d$ and $\mu_1$ is
concentrated in a neighborhood of $\supp(\mu_0)$ of radius $\pit$ the Monge
formulation \eqref{eq:172intro} has a unique solution.

After the transformation $\xi_0:=\frac 12(\ee^{\varphi_0}-1)$
(which linearizes the second integrand in
the duality formula \eqref{eq:HK.dual.intro}), 
we also obtain a family of maps, for $t\in [0,1]$,
\begin{equation}
  \label{eq:173}
  \bfT_{0\to t}(x)=x+\bfarctan\Big(\frac{t\nabla\xi_0}{1{+}2t\xi_0(x)}\Big),\quad
  q^2_{0\to t}(x):=(1{+}2t\xi_0(x))^2+t^2|\nabla\xi_0(x)|^2,
\end{equation}
with the following properties:
\begin{enumerate}
\item
$(\bfT_{0\to1},q_{0\to 1})$ is the unique solution of
\eqref{eq:172intro} and provides the beautiful formula
\begin{equation}
  \label{eq:172}
  \HK^2(\mu_0,\mu_1)=
  \int_{\R^d} \Big(4\xi_0^2+|\nabla\xi_0|^2\Big)\dd\mu_0,
\end{equation}
showing that the (closure of the) space of $\mathrm C^1_c(\R^d)$ functions
with respect to the Hilbertian norm
\begin{equation}
  \|\xi\|_{H^{1,2}(\R^d,\mu)}^2=
  \int_{\R^d}\Big(4\xi^2+|\nabla\xi|^2\Big)\dd\mu\label{eq:174}
\end{equation}
provides the natural notion of tangent space $\mathrm{Tan}_\mu \MMM(\R^d)$ and
a nonsmooth Riemannian formalism in $(\MMM(\R^d),\HK)$ as for the Otto calculus
in $(\mathcal P_2(\R^d),W_2)$.
\item
The curve $\mu_t=(\bfT_{0\to t},q_{0\to t})_\star\mu_0$ is an
explicit characterization of the geodesic interpolation solving
\eqref{eq:HKdynamical}.  A crucial fact is that for $\mu_0$-a.e.\,$x$ the curve
$[\bfT_{0\to t}(x),q_{0\to t}(x)]$ is a geodesic in the cone space $\mfC$
interpolating the points $[x,1]$ and $[\bfT_{0\to 1}(x),q_{0\to 1}(x)]$.
\end{enumerate}
It is then natural to investigate if the potential $\xi_0$ can be used to build
an optimal solution $\xi_t$ of \eqref{eq:16intro}, 
which should at least formally solve the Hamilton-Jacobi equation
\begin{equation}
\label{eq:added0}
    \partial_t\xi_t+\frac 12|\nabla\xi_t|^2+2\xi_t^2=0\quad
    \text{on the support of $\mu$ in 
        $(0,1)\times \R^d$.}
\end{equation}
This problem will be investigated
in Section \ref{sec:HJ}, by a detailed analysis of the regularity of the
 forward solutions to \eqref{eq:16intro} provided by the 
generalized Hopf--Lax
formula (see \eqref{eq:23})
\begin{equation}
\xi_t(x)=\xi(t,x) = \big(\HopfLax t{ \xi_0}\big)(x)=
  \frac 1{t} \HopfLax 1{\big(t\xi_0(\cdot)\big)}(x)=\inf_{y\in\OOmega}
\frac{1} 
{2t}
\Big(1-\frac{\tcosq{\pit}{|x{-}y|}}{1+
2t\xi_0(y)}\Big).\label{eq:23intro}
\end{equation}
It is well known that one cannot expect smoothness of such a solution;
however, the particular structure of transport duality 
suggests that the final value $\xi_1$ given by \eqref{eq:23intro} 
corresponds to the optimal potential $\varphi_1$ of the dual formulation \eqref{eq:HK.dual.intro} via the transformation $\xi_1=\frac 12(1-\mathrm e^{-2\varphi_1})$, so that 
the initial and final optimal potentials 
$\xi_0$ and $\xi_1$ are simultaneously linked by the forward-backward relation
\begin{equation}
\label{eq:added1}
    \xi_1=\mathscr P_1\xi_0,\quad
    \xi_0=\mathscr R_1(\xi_1)\quad\text{where }
    \mathscr R_t(\eta):=-\mathscr P_t(-\eta)\text{ is the backward flow.}
\end{equation}
Following the approach of \cite[Cha.\,7]{Vill09OTON}
(see also \cite[Sec.\,8]{LiMiSa18OETP}) and using the reversibility in time of geodesics, we 
can add to 
the family of forward potentials $\xi_t$ 
given by \eqref{eq:23intro}
the crucial information 
provided by the backward solutions $\bar\xi_t$
starting from $\xi_1$:
\begin{equation}
    \bar\xi_t:=\mathscr R_{1-t}\xi_1=
    -\mathscr P_{1-t}\big({-}\xi_1\big)\quad 
    \text{for } t\in [0,1].
\end{equation}
In general, $\xi_t$ and $\bar \xi_t$ do not coincide for $t\in (0,1)$ but still satisfy 
\begin{equation}
    \xi_t(x)\ge \bar \xi_t(x)\quad\text{in }(0,1)\times \R^d,\quad
    \xi_0=\bar \xi_0,\quad\xi_1=\bar\xi_1.
\end{equation}
The crucial fact arising from the
optimality condition \eqref{eq:added1},
and the geometric property of the geodesic $(\mu_t)_{t\in [0,1]}$
is that for every $t\in [0,1]$
\begin{equation*}
    \text{the support of $\mu_t$
is contained in the \emph{contact set}}
\quad
\Xi_t:=\bigset{x\in \R^d} {\xi_t(x)=\bar \xi_t(x)}.
\end{equation*}
On the contact set $(\Xi_t)_{t\in [0,1]}$, we can combine the 
(delicate) first- and second-order 
super-differentiability properties of $\xi_t$ arising from the  inf-convolution structure of \eqref{eq:23intro}
with the corresponding sub-differentiability properties 
exhibited by $\bar \xi_t$.

Using tools from nonsmooth analysis, we
are then able to give a rigorous meaning to the characteristic flow associated with
\eqref{eq:added0}, 
i.e.\ to the maps $t\mapsto \bfT(t,\cdot)=\bfT_{s\to t}(\cdot)$, $t\mapsto q(t,\cdot)=q_{s\to t}(\cdot)$ solving
(we omit to write 
the explicit dependence on $x$ when not needed)
\begin{equation}
  \label{eq:157intro}
  \left\{
    \begin{aligned}
      & \dot \bfT(t)=\nabla\xi_t(\bfT(t)),\\
      & \dot q(t)=2\xi_t(\bfT(t))q(t),
    \end{aligned}
  \right.\quad \text{in }(0,1),\quad
  \bfT(s,x)=x,\ q(s,x)=1.
\end{equation}
Moreover, we will prove that $\bfT_{s\to t}$ is a family of bi-Lipschitz maps
on the contact sets obeying a natural concatenation property. As can be
expected, the maps $\bfT_{s\to t},q_{s\to t}$ provide a precise representation
of the geodesics via $\mu_t=(\bfT_{s\to t},q_{s\to t})_\star \mu_s$ for all
$s,t\in (0,1)$. In particular $(\bfT_{s\to t},q_{s\to t})$ is an optimal
\TRAGRO\ pair between $\mu_s$ and $\mu_t$ minimizing the cost of
\eqref{eq:172intro}.  \EEE

Using
this valuable information, in Section~\ref{se:GeodCurves} we obtain various relevant 
structural properties of geodesics in $(\MMM(\R^d),\HK)$ such as
{\em non-branching, localization, and regularization effects}.  In particular,
independently of the regularity of $\mu_0$ and $\mu_1$, we will show that 
for $s\in (0,1)$ the Monge problem between $\mu_s$ and $\mu_0$ or between
$\mu_s$ and $\mu_1$ always admit a unique solution, a property which is well
known in the Kantorovich--Wasserstein framework.
  
Surprisingly enough, despite the lack of global regularity, we will also establish precise formulae for the
first and second derivative of the differential of $\bfT_{s\to t}$ (and thus the second order differential of
$\xi_t$) along the flow, which coincides with the equations that one obtains by formally
differentiation using the joint information of the Hamilton--Jacobi equation \eqref{eq:added0} and 
\eqref{eq:157intro} assuming sufficient regularity.
For instance, differentiating in time the first equation of \eqref{eq:157intro} 
and differentiating in space \eqref{eq:added0} one finds
\begin{equation*}
    \ddot \bfT(t)=\partial_t\nabla\xi_t(\bfT(t))+\rmD^2\xi_t\nabla\xi_t(\bfT(t)),\quad
    \partial_t\nabla\xi_t=
    -\rmD^2\xi_t\nabla\xi_t+4\xi_t\nabla\xi_t,
\end{equation*}
which yield
\begin{subequations}
\begin{equation}
    \ddot \bfT(t)=4\xi_t(\bfT(t))\nabla\xi_t(\bfT(t)).
\end{equation}
For $q(t)$, $\sfB(t):=\rmD\bfT_{s\to t}$,
and its determinant $\delta(t):=\det \sfB(t)$ similar, just more involved, calculations yield 
the crucial second order equations
\begin{align}
\label{eq:ddot1}
\ddot q(t)&=|\nabla\xi_t(\bfT(t))|^2q(t),\\
\label{eq:ddot2}
\ddot\sfB(t)&=-4\Big(\nabla\xi_t\otimes\nabla\xi_t+\xi_t\rmD^2\xi_t\Big)\circ \bfT(t)\cdot\sfB(t),\\
\label{eq:ddot3}
    \ddot\delta(t)&=\Big((\Delta\xi_t)^2-|\rmD^2\xi_t|^2-
    4|\nabla\xi_t|^2-4\xi_t\Delta\xi_t\Big)\circ \bfT(t)\cdot\delta(t).
\end{align}
\end{subequations}
In our case, even though we do not have enough regularity to justify the above formal
computations, 
we can still derive them rigorously by a deeper analysis using the
variational properties of the contact set.
Even if our discussion is restricted to 
the Hellinger--Kantorovich case and uses the particular form of the Hopf--Lax semigroup \eqref{eq:23intro}
and its characteristics \eqref{eq:173}, we think
that our argument applies to more general cases
and may provide new interesting estimates also in the typical balanced case of Optimal Transport.

Such regularity and the related second order estimates
are sufficient to express the Lebesgue density $c_t$ of the
measures $\mu_t$ and thus to obtain crucial information on its behavior along
the flow.  In particular, Corollary \ref{cor:start} shows that $c(t,\cdot)$ is given by
\begin{subequations}
  \label{eq:alpha.delta}
\begin{align}
  \label{eq:alpha.delta.A}
&c(t, y) \big|_{y = \bfT_{s\to t}(x)} = c(s,x) \frac{\alpha_s(t,x)}{\delta_s(t,x)}
\quad \text{with }\\
  \label{eq:alpha.delta.B}
  &\alpha_s(t,x)=(1{+}2(t{-}s)\psxi_s(x))^2 +
    (t{-}s)^2 |\nabla\psxi_s(x)|^2 =q_{s\to t}(x)\\
  &\delta_s(t,x):= \det( \rmD \bfT_{s\to t}(x)), 
\end{align}
and the time-dependent \TRAGRO\ mapping $\bfT_{s\to t},q_{s\to t}$ are given in
terms of $\psxi$ via \eqref{eq:157intro} and the analog of \eqref{eq:173}.
\end{subequations}
In particular, we will show that if $\mu_s\ll \Leb d$ for some $s\in (0,1)$
then $\mu_t\ll\Leb d$ for every $t\in (0,1)$
and combining \eqref{eq:ddot1}, \eqref{eq:ddot2}, and
\eqref{eq:alpha.delta.A}
we will also prove that $c_t$ is a convex function
along the flow maps $\bfT_{s\to t}$.

\subsection{Geodesic \texorpdfstring{$\lambda$}{lambda}-convexity of functionals} 
\label{su:IntroGLCvx}

The second part of the paper is devoted to establish necessary and sufficient
conditions for geodesic $\lambda$-convexity of energy functionals $\calE$
defined for a closed and convex domain $\Omega \subset \R^d$ with
non-empty interior in the form
 \begin{equation}
\label{eq:calE} 
\calE(\mu)=
\text{\footnotesize$\ds\int$}_{\!\!\Omega} \! E(c(x)) \dd x +
E'_\infty \mu^\perp(\Omega)\quad \text{for} \quad \mu=c\Leb d{+}\mu^\perp \text{ with } 
  \mu^\perp \perp \Leb d,
\end{equation}
where $E'_\infty:= \lim_{c\to\infty}E(c)/c\in \R\cup\{+\infty\}$ is the recession 
constant and $E(0)=0$ holds.

In \cite[Prop.\,19]{LiMiSa16OTCR} it was shown that the total-mass
functional $\mathscr M:\mu
\mapsto \mu(\OOmega)$ has the surprising property that it is
\AAA exactly quadratic \EEE along $\HK$ geodesics $\gamma:[0,1]\to
\MMM(\OOmega)$, namely 
\begin{equation}
  \label{eq:calMcvx}
  \mathscr M(\gamma(t)) = (1{-}t)\mathscr M(\gamma(0)) + t \mathscr M(\gamma(1)) -
t(1{-}t) \HK(\gamma(0),\gamma(1))^2 \ \text{ for } t \in [0,1].
\end{equation}
Thus, as a first observation we see that a density function $E$
generates a geodesically $\lambda$-convex functional $\calE$ if and
only if $E_0: c\mapsto E(c)- \lambda c$ generates a geodesically convex
functional (i.e.\ geodesically $0$-convex). Hence, subsequently we
can restrict to $\lambda=0$.

To explain the necessary and sufficient conditions on $E$ for $\calE$ 
to be geodesically convex,  
we first look at the 
differentiable case, and  we define the
shorthand notation 
\[
\eps_0(c)=E(c), \quad \eps_1(c)=c E'(c), \quad \eps_2(c) = c^2 E''(c).
\]
For the Kantorovich--Wasserstein distance $\sfW$ the necessary and
sufficient conditions are the so-called McCann conditions \cite{Mcca97CPIG}:
\begin{equation}
  \label{eq:McCann.cond}
  \begin{aligned}
\eps_2(c) &\geq \frac{d{-}1}d \big(
\eps_1(c)-\eps_0(c)\big) \geq 0 \text{ for all }c>0\\[0.6em]
&\Longleftrightarrow \ \left\{ \ba{l} r \mapsto r^{d}E(r^{-d}) 
\text{ is lower semi-continuous and convex and } \\[0.3em]
   r \mapsto (d{-}1)r^{d}E(r^{-d}) 
  \text{ is non-increasing on }  {]0,\infty[},
  \ea\right.
\end{aligned}
\end{equation}
see also \cite[Prop.\,9.3.9]{AmGiSa08GFMS}. 
For the Hellinger--Kakutani distance we simply need the condition
\begin{equation}
  \label{eq:Hell.cond}
2 \eps_2(c) + \eps_1(c) \geq 0 \quad
    \Longleftrightarrow \quad 
\Big( r \mapsto E(r^2) \text{ is convex }\Big). 
\end{equation}
In the case of differentiable $E$, our main result yields the
following necessary and sufficient conditions for geodesic convexity
of $\calE$ on $(\MMM(\OOmega),\HK)$, see  Proposition
\ref{pr:N.E.conds},  
\begin{equation}
  \label{eq:HK.cond.diff}
(d{-}1)\big(\eps_1(c)-\eps_0(c)\big)\geq 0 \quad \text{and}\quad \bbB(c) \geq 0 
\quad \text{for all }c>0,
\end{equation}
where the matrix $\bbB(c)\in \R^{2\ti 2}_\text{sym}$ is given by 
\[
  \bbB(c) := 
\bma{cc} \eps_2(c)-\frac{d-1}d\big(\eps_1(c){-}\eps_0(c)\big) & \eps_2(c)
-\frac12\big(\eps_1(c){-}\eps_0(c) \big) \\[0.3em]  \eps_2(c)
-\frac12\big(\eps_1(c){-}\eps_0(c) \big) & \eps_2(c) + \frac12 \eps_1(c) 
\ema.
\]
We immediately see that the non-negativity of the diagonal element
$\bbB_{11}(c)$ gives the first McCann condition in \eqref{eq:McCann.cond}, and
$\bbB_{22}(c) \geq 0 $ gives \eqref{eq:Hell.cond}. However, the
condition $\bbB(c)\geq 0$ is strictly stronger, since e.g.\ it implies that the
additional condition $(d{+}2)\eps_1(c)-2\eps_0(c)\geq 0$ holds, see
\eqref{eq:HK.cond3/2}. This condition means that $c \mapsto c^{-2/(d+2)}E(c)$ has to be
non-decreasing, which will be an important building block for the main geodesic
convexity result.

Indeed, our main result in Theorem~\ref{th:Geod.Cvx} is formulated
for general lower semi-continuous and convex functions $E:{[0,\infty[}\to
\R\cup\{\infty\}$ without differentiability assumptions. The
conditions on $E$ can be formulated most conveniently 
in terms of the auxiliary function $N_E:{]0,\infty[}^2\to
\R\cup\{\infty\}$ defined via  
\begin{subequations}
  \label{eq:HK.cond.N}
\begin{equation}
  \label{eq:N.def}
N_E(\rho,\gamma )= \big(\frac\rho\gamma\big)^d
E\Big(\frac{\gamma^{2+d}}{\rho^d}\Big).  
\end{equation}
Then, $\calE$ defined in \eqref{eq:calE} is geodesically convex if and
only if $N_E$ satisfies 
\begin{align}
  \label{eq:HK.cond.N.b}
&N_E:{]0,\infty[}^2\to \R\cup\{\infty\}  \text{ is convex, and} \\ 
  \label{eq:HK.cond.N.c} 
&  \rho \mapsto  (d{-}1)N_E(\rho,\gamma) \text{ is non-increasing}. 
\end{align} 
\end{subequations} 
The McCann conditions \eqref{eq:McCann.cond} are
obtained by looking at $N_E(\cdot,\gamma)$ for fixed $\gamma$, while the
Hellinger--Kakutani condition \eqref{eq:Hell.cond} follows by looking
at $s \mapsto N_E(s\rho,s\gamma)$ for fixed $(\rho,\gamma)$. 

The proof of the sufficiency and necessity of condition 
\eqref{eq:HK.cond.N} for geodesic convexity of $\calE$ is based on the explicit
representation \eqref{eq:alpha.delta} of the geodesic curves giving 
\[
\calE(\mu(t))=\int_\Omega E(c(t,y))\dd y = \int_\Omega e(t,x)
\dd x \ \text{ where }
 e(t,x):=\delta_s(t,x)\,E\Big(c_s(x)\frac{\alpha_s(t,x)}{\delta_s(t,x)}\Big).
\] 
By definition, we have $\alpha_s(t,x)\geq 0$, and Corollary \ref{cor:start}
guarantees $\delta_s(t,x)>0$. Hence, we can introduce the two functions
\[
\gamma(t,x)=\big(c_s(x)\alpha_s(t,x)\big)^{1/2} \quad \text{ and }
\quad \rho(t,x) = \big(c_s(x)\alpha_s(t,x) 
  \big)^{1/2} \delta_s(t,x)^{1/d}, 
\]
which connect the densities $e(t,x)$ with the function $N_E$ defined in 
\eqref{eq:N.def} in the form
\[
e = \delta \,E\big(\,c\,\frac{\alpha}{\delta}\,\big) = N_E(\rho,\gamma).
\]

For smooth $E$ we have smooth $N_E$ and may show convexity of $t\mapsto
e(t,x)$ via 
\begin{align*}
\pl_t^2 e(t,x)& =:\ddot e = \Big\langle \rmD^2
N_E(\rho,\gamma)  \binom{\dot\rho}{\dot\gamma},
\binom{\dot\rho}{\dot\gamma} \Big\rangle + \Big\langle  \rmD
N_E(\rho,\gamma), \binom{\ddot\rho}{\ddot\gamma} \Big\rangle\geq 0.
\end{align*}
By convexity of $N_E$, the term involving $\rmD^2 N_E$ is non-negative, so
it remains to show 
\begin{equation}
  \label{eq:N.ddot}
  \pl_\rho N_E(\rho,\gamma) \ddot \rho + \pl_\gamma N_E(\rho,\gamma) \ddot
\gamma\geq 0. 
\end{equation}
To establish this, we use first that the scaling property
$N_E(s^{1+d/2}\rho,s \gamma)=s^2 N_E(\rho,\gamma)$ for all $s>0$ (which follows
from the definition of $N_E$ via $E$) and the convexity of $N_E$ imply
\begin{equation}
\label{eq:N.monotonicity}
(1{-}4/d^2) \;\! \rho\;\! \pl_\rho N_E(\rho,\gamma)+ \gamma\;\!
 \pl_\gamma N_E(\rho,\gamma)\geq 0,
\end{equation}
see Proposition \ref{pr:N.E.monot}. Second, we rely on a 
nontrivial curvature estimate for $(\rho,\gamma)$, namely  
\begin{align}
 \label{eq:I.curvature}
\frac{\ddot\gamma(t,x)}{\gamma(t,x)} \geq 0 \quad \text{and} \quad 
 \frac{\ddot\rho(t,x)} {\rho(t,x)} \leq \big(1-\frac4d\big)\!\;
\frac{\ddot\gamma (t,x)} {\gamma(t,x)}.
\end{align}
Estimates \eqref{eq:I.curvature} are provided
in Proposition \ref{pr:Est.alpha.delta} and strongly rely on the explicit
representation and the regularity
properties of the geodesics developed in Sections \ref{sec:HJ}
and~\ref{se:GeodCurves}.%

Combining \eqref{eq:I.curvature}  with $\pl_\rho N_E(\rho,\gamma)\leq 0$, the
desired relation \eqref{eq:N.ddot} easily follows, see Section
\ref{se:Geod.Conv}.  Finally, a simple integration over $\OOmega$ provides the
convexity of $t \mapsto \calE(\mu(t))$.  Note that we have indeed the larger factor
$(1{-}4/d^2)$ in \eqref{eq:N.monotonicity} while the curvature estimate in
\eqref{eq:I.curvature} has the smaller and hence ``better'' factor
$(1{-}4/d)$. 

As a consequence, we find that the power functionals $\calE_m$ with
$E_m(c)=c^m$ with $m>1$ are all geodesically convex, see Corollary
\ref{co:m.q}. This result was already exploited in
\cite[Thm.\,2.14]{DimChi20TGMH}. We can study the
discontinuous ``Hele--Shaw case'' $E(c)=-\lambda c$ for $c\in [0,1]$ and
$E(c)=\infty$ for $c>1$.  Moreover, in dimensions $d=1$ or $2$ the
densities $E_q(c)=- c^q$ with $q\in [\frac{d}{d+2},\frac12]$ also lead
to geodesically convex functionals $\calE_q$, see again Corollary \ref{co:m.q}.

\GGG Two important differences with the balanced Kantorovich--Wasserstein case
are worth noting: First, the Boltzmann logarithmic entropy functional
corresponding to $E(c)=c\log c$ is not geodesically $\lambda$-convex for any
value of $\lambda$, see Example \ref{ex:BoltzNotGLCvx}.  Second, if the space
dimension $d$ is larger than or equal to $3$, then there are no geodesically
convex power functionals of the form $E(x)=-c^m$ with exponent $m<1$, see
Example \ref{ex:q.le.1}. \AAA Some of these statements follow easily by
observing that $\mu_t = t^2 \mu_1$ is the unique geodesic connecting $\mu_0=0$
and $\mu_1$. \EEE

\subsection{Applications and outlook} 
\label{su:IntroAppl}

In \cite{Flei21MMAC, LasMie22?EVIH}, the JKO scheme (minimizing movement scheme) for a gradient system $(\MMM(\Omega),\HK_{\alpha,\beta},\calE)$ is considered, i.e., for $\tau>0$ we iteratively define
\begin{equation}
    \label{eq:RDE}
  \mu_{k\tau} \in \mathrm{Arg\,Min} \Bigset{\frac1{2\tau} \HK^2_{\alpha,\beta} 
  ( \mu_{(k-1)\tau}, \mu) + \calE(\mu)}
  { \mu \in \MMM(\Omega)} 
\end{equation}
and consider the limit $\tau\downarrow 0$ (along subsequences) to obtain \emph{generalized minimizing movements} (GMM) (cf.\:\cite{AmGiSa05GFMS}). Under suitable conditions, including the assumption $ \calE(\mu)= \int_\Omega \big( E(c) + c V\big) \dd x $ with $\mu= c\calL^d$ and  $E$ superlinear,  it is shown in \cite[Thm.\,3.4]{Flei21MMAC} that all GMM $\mu$ have the form $\mu(t)=c(t) \calL^d$, and the density $c$ is a weak solution of the reaction-diffusion equation  
\[
\pl_t c = \alpha \DIV\big( c \nabla (E'(c){+}V) \big)  - \beta \,u\,\big(E'(c){+}V \big) \ \text{ in }\Omega, 
  \quad \   
  c \nabla (E'(c){+}V)\cdot \rmn=0 \ \text{ on }\pl\Omega.
\]
In \cite{Li06GESE}, the equation $u_t =0= \Delta u + a u \log u + bu$ is
studied, whose solutions are steady states for HK gradient flows for
$\calE(u)=\int_\OOmega u\log u \dd x$.  We also refer to \cite{PeQuVa14HSAM,
  DimChi20TGMH}, where equation \eqref{eq:RDE} was studied for
$E(c)=\frac1m c^m - \lambda c$ and $V\equiv 0$. The linear functional
$\Phi(\mu)=\int_\OOmega V(x)\dd \mu$ for a given potential
$V\in \rmC^0(\OOmega)$ can easily be added, as its geodesic $\lambda$-convexity
is characterized in \cite[Prop.\,20]{LiMiSa16OTCR}. Note that \GGG our main
convexity result, proved here for the first time, plays an important role in
the existence and uniqueness results of \cite{DimChi20TGMH}, cf.\ Thm.\,2.14
there.  \EEE

In \cite{LasMie22?EVIH} it is shown that the GMM for the gradient system $(\MMM(\Omega),\HK_{\alpha,\beta},\calE)$ are EVI$_\lambda$ solutions in the sense of \cite{MurSav20GFEV}. Again the main ingredient is the geodesic $\lambda$-convexity of $\calE$ in the form \eqref{eq:calE} contained in our main Theorem \ref{th:Geod.Cvx}.

\bigskip
\centerline{\large\bfseries Main notation}

\smallskip
\halign{\small$#$\hfil\ &\small#\hfil
\cr
\MMM(X),\ \MMM_2(X)&finite positive Borel measures on $X$ (with finite
quadratic moment)
\cr
\cP(X),\ \cP_2(X)&Borel probability measures on $X$ (with finite
quadratic moment)
\cr
T_\sharp\mu&push forward of $\mu\in \MMM(X)$ by a map $T:X\to Y$:
\eqref{eq:push_forw}
\cr 
\mu=c\Leb d+\mu^\perp&Lebesgue decomposition of a measure $\mu\MMM(\R^d)$
\cr 
\rmC_\rmb(X)&continuous and bounded real functions on $X$
\cr
\cos_a(r)& truncated function $\cos\big(\min\{a,r\}\big)$, $a>0$ (typically
$a=\pi/2$)
\cr
\mathsf W_X(\mu_1,\mu_2)&Kantorovich--Wasserstein distance in 
$\calP_2(X)$
\cr
\bfsin,\bftan, \bfarctan,\cdots&vector-valued version of the usual
scalar functions, see \eqref{eq:5}
\cr
\HK(\mu_1,\mu_2)&Hellinger--Kantorovich distance in $\MMM(X)$:
Section \ref{s:HK.definition}
\cr
(\mfC,\ \sfd_{a,\mfC}),\ \mfo&metric cone on $\R^d$ and its vertex, see Subsection
\ref{sss:conespace}
\cr
\sfW_{a,\mfC}&$L^2$-Kantorovich--Wasserstein distance on $\calP_2(\mfC)$
induced
by $\sfd_{a,\mfC}$ \cr
\sfx,\sfr&coordinate maps on $\mfC$, see Subsection \ref{sss:conespace}
\cr
\pi^0,\pi^1&coordinate maps on a Cartesian product $X_0\times X_1$, $\pi^i(x_0,x_1)=x_i$
\cr
\mfh&homogeneous projection from $\MMM_2(\mfC)$ to $\MMM(\R^d)$,
see \eqref{eq:frak.P}
\cr
S_i, S_i', S_i'', S_i^{\pi/2},\mu_i',\mu_i''&see
\eqref{eq:33}-\eqref{eq:43}
\cr
(\bfT,q)_\star & action of a \TRAGRO\ map,
Def.~\ref{def:dilation-transport}
\cr
\mathrm{AC}^p([0,1];X) & space of curves $\rmx:[0,1]\to
X$ with $p$-integrable metric speed
\cr
\varphi_0^{\LLL\to}, \varphi_1^{\shortleftarrow\LLL} &
forward and backward $\LLL$-transform for cost function $\LLL$, see~\eqref{eq:Def.Ell.Trafo}
\cr
D_i',D_i'' & domains of $\nabla\varphi_i$ and $\rmD^2\varphi_i$, see Theorem~\ref{thm:regularity0}
\cr
\xi_s=\HopfLax s\xi, \bar\xi_s=\RopfLax s\bar\xi& for- and backward solution of  Hamilton--Jacobi equation, \eqref{eq:23}, \eqref{eq:ForwardBackwardSolution}
\cr
\Xi_s& contact set of forward and backward solutions $\xi_s$, $\ol\xi_s$, see \eqref{eq:116}
\cr
(\bfT_{s\to t}, q_{s\to t}) & \TRAGRO\ map induced by for/backward solutions, Theorem~\ref{thm:HJ1}
\cr
}

\section{The Hellinger--Kantorovich distance}
\label{s:HK.definition}
In this section, we recall a few properties and equivalent
characterizations of the Hellinger--Kantorovich distance
from \cite{LiMiSa16OTCR,LiMiSa18OETP}, 
that will turn out to be crucial in the following.

First, we fix some notation that we will extensively use:
Let $(X,\sfd_X)$ be a complete and separable metric space. 
In the present paper $X$ will typically be $\R^d$ with the Euclidean distance,
a closed convex subset thereof, the cone space
$\mfC$ on $\R^d$ (see Subsection \ref{sss:conespace}),
product spaces of the latter two, etc.
We will denote by $\MMM(X)$ the space of all
non-negative and finite Borel measures on $X$ endowed with the 
weak topology induced by the duality with the continuous and bounded
functions of $\rmC_\rmb(X)$. The subset of measures with 
finite quadratic moment will be denoted by $\MMM_2(X)$. 
The spaces $\cP(X)$ and $\cP_2(X)$ are the corresponding subsets of probability
measures.

If $\mu\in \MMM(X)$ and $T:X\to Y$ is a Borel map
  with values in another metric space $Y$,  
then $T_\sharp \mu$
denotes the push-forward measure on $\MMM(Y)$, defined by
\begin{equation}\label{eq:push_forw}
  T_\sharp\mu(B):=\mu(T^{-1}(B))\quad
  \text{for every Borel set $B\subset Y$}.
\end{equation}
We will often denote elements of $X\times X$ by $(x_0,x_1)$ and
the canonical projections by $\pi^i:(x_0,x_1)\to x_i$, $i=0,1$.
A coupling on $X$ is a measure $\ggamma\in \MMM(X{\times} X)$ with marginals
$\gamma_i:=\pi^i_\sharp \ggamma$.

Given two measures $\mu_0,\mu_1\in \MMM_2(X)$ with equal mass $
\mu_0(X)=\mu_1(X)$, their
(quadratic) Kantorovich--Wasserstein distance $\sfW_X$ is defined by
\begin{equation}
  \label{eq:11}
  \begin{aligned}
    \sfW_{\kern-1pt X}(\mu_0,\mu_1)^2:=\min\bigg\{\iint &
    \sfd_X(x_0,x_1)^2\dd\ggamma(x_0,x_1)\,\Big|\, 
    \\&\ggamma\in \MMM(X{\times} X),\ \pi^i_\sharp\ggamma=\mu_i,\ i=0,1\bigg\}.
  \end{aligned}
\end{equation}
We refer to \cite{AmGiSa08GFMS} for a survey on the Kantorovich--Wasserstein
distance and related topics.

\subsection{Equivalent formulations of the Hellinger--Kantorovich
  distance}
\label{su:Form.HK}
The Hellinger--Kantorovich distance was introduced in 
\cite{LiMiSa18OETP,LiMiSa16OTCR} and independently 
in \cite{KoMoVo16NOTD} and \cite{CPSV18UOTG,CPSV18IDOT}.
It is a generalization of the Kantorovich--Wasserstein distance
to  arbitrary non-negative and finite measures by taking
creation and annihilation of mass into account. Indeed, the latter
can be associated with a different notion of distance, namely the
Hellinger--Kakutani distance, see \cite{Hell09NBTQ} and \cite{Schi97GAFV}.
In this sense, the Hellinger--Kantorovich distance should be viewed
as an infimal convolution of the Kantorovich--Wasserstein and the 
Hellinger--Kakutani distance, cf.\ \cite[Rem.\,8.19]{LiMiSa18OETP}.

In \cite{LiMiSa18OETP}, five different equivalent formulations of
the Hellinger--Kantorovich distance are given: (i) the dynamical
formulation, (ii) the cone space formulation, (iii) the optimal
entropy-transport problem, (iv) the dual formulation in terms of
Hellinger--Kantorovich potentials, and (v) the formulation using
Hamilton--Jacobi equations. We will present and briefly discuss each
of them below,  as all are useful for our analysis of geodesic
convexity. 

In the following, we consider the Hellinger--Kantorovich distance for
measures on the domain $\R^d$. However, it is easy to see that all arguments
also work in the case of a closed and convex domain $\Omega\subset\R^d$. 
In particular, the latter is a complete, geodesic space.

\subsubsection{Dynamic approach}
\label{subsec:HKdynamic}
A first approach to the Hellinger--Kantorovich distance
is related to the dynamic formulation, which naturally 
depends on two positive parameters
$\alpha,\beta>0$:
they control the relative
strength of the Kantorovich--Wasserstein part and of the
Hellinger-Kakutani part (see \cite[Section
  8.5]{LiMiSa18OETP}).
\begin{definition}[The dynamic formulation]
  \label{def:dynamic}
  For every $\mu_0,\mu_1\in \MMM(\OOmega)$ we set 
  \begin{equation}
    \begin{aligned}
      \HK_{\alpha,\beta}(\mu_0,\mu_1)^2= \min \bigg\{ &\int_{0}^1\!\!
      \int_\OOmega \big(\alpha\, |\Upsilon( t,x)|^2 {+}\beta\xi(t,x)^2\big)
      \dd\mu_t(x)\dd t\,\Big|\,  \\
      & \mu\in \rmC([0,1];\MMM(\R^d)), \ \mu_{t=i}=\mu_i,\ \text{\upshape(gCE)
        holds} \bigg\},
    \end{aligned}\label{eq:15}
  \end{equation}
  where the generalized continuity equation  
  for the Borel vector and scalar fields
  $\Upsilon: (0,1)\times \R^d\to\R^d$ and $\xi:(0,1)\times \R^d\to \R$ reads
  \[
    \text{\upshape (gCE) } \qquad \frac\partial{\partial t} \mu +\alpha
  \div (\mu\;\! \Upsilon )  =\beta\, \xi \mu\quad \text{in } 
    \mathcal D'( (0,1)\ti \OOmega).
  \]
\end{definition}

Notice that \eqref{eq:15} yields in particular that $\mu\;\!\GGG \Upsilon\EEE$ and
$\xi\mu$ are (vector and scalar) measures with finite total mass, so
that the canonical formulation of (gCE) in $\mathcal D'( (0,1)\ti
\OOmega)$ makes sense.  For optimal solutions one has
$\Upsilon(t,x)=\nabla \xi(t,x)$ and the dual potential solves the 
generalized Hamilton--Jacobi equation 
\begin{equation}
\label{eq:genHJeqn}
 \pl_t\xi + \frac\alpha2 |\nabla \xi|^2 +\frac\beta2 \xi^2=0 
\end{equation}
in a suitable sense \cite[Theorem 8.20]{LiMiSa18OETP}.

A simple rescaling technique shows that it is sufficient to restrict 
ourselves to a specific choice of the parameters $\alpha$ and $\beta$.
In fact, it is easy to see that for every $\theta>0$ we have 
\begin{equation*}
  \HK_{\alpha,\beta}(\mu_0,\mu_1)^2
  =\theta\HK_{\theta\alpha,\theta\beta}(\mu_0,\mu_1)^2.
\end{equation*}
Moreover, if $\lambda>0$ \GGG and we consider the \AAA spatial \GGG dilation
$H:x\mapsto \lambda x$ in $\OOmega$, \EEE
we find 
\begin{equation*}
  \HK_{\alpha,\beta}(\mu_0,\mu_1)^2=
  \HK_{\alpha/\lambda^2,\beta}(H_\sharp\mu_0,H_\sharp\mu_1)^2.
\end{equation*}
Choosing $\lambda:=\sqrt{4\alpha/\beta}$, $\theta= 4/\beta$, and setting
$\HK:=\HK_{1,4}$ we get
\begin{equation*}
  \HK_{\alpha,\beta}(\mu_0,\mu_1)^2=
  \frac 4\beta
  \HK_{4\alpha/\beta,4}(\mu_0,\mu_1)^2=
  \frac 4\beta\HK(H_\sharp\mu_0,H_\sharp\mu_1)^2.
\end{equation*}
Therefore, in order to keep simpler notation,
in the remaining paper we will mainly consider the case
$\alpha=1$ and $\beta=4$.

\subsubsection{Cone space formulation}
\label{sss:conespace} 

There is a second characterization that connects $\HK$ with the
classic Kantorovich--Wasser\-stein distance on the extended cone
$\mfC:=(\OOmega\times[0,\infty[)/\!\sim$, where $\sim$ is the
equivalence relation which identifies all the points $(x,0)$ with the
vertex $\mfo$ of $\mfC$.  More precisely, we write $(x_0,r_0)\sim
(x_1,r_1)$ if and only if $x_0=x_1$ and $r_0 = r_1$ or $r_0=r_1=0$
and introduce the notation $[x,r]$ to denote the equivalence class
associated with $(x,r)\in\OOmega\times[0,\infty[$.  The cone $\mfC$ is
a complete metric space endowed with the
 cone distances
\begin{equation}
  \label{eq:7}
  \sfd_{a,\mfC}(\mfz_0,\mfz_1)^2:=
  r_0^2+r_1^2-2r_0r_1  \tcos{a}{ |x_1{-}x_0|},
  \quad  \mfz_i=[x_i,r_i],\ a\in (0,\pi],
\end{equation}
see e.g.\ \cite[Sect.\,3.6.2]{BuBuIv01CMG}, where we use the abbreviation
$\cos_a(r):=\cos\big(\!\min\{a,r\}\big)$. Notice that the projection map
$(x,r)\mapsto [x,r]$ is bijective from $\R^d\times (0,\infty)$ to
$\mfC_\STAR:=\mfC\setminus\{\mfo\}$; we will denote by $(\sfx,\sfr)$ its
inverse, which we extend to $\mfo$ by setting $\sfx(\mfo)=0,\ \sfr(\mfo)=0$.

The most natural choice of the parameter $a$ in
\eqref{eq:7} is $a:=\pi$: in this case the cone $(\mfC,\sfd_{\pi,\mfC})$ is
a geodesic space, i.e.,
given $\mfz_i=[x_i,r_i]$, $i=0,1$, there exists a curve
$\mfz_t=[x_t,r_t]=\geo t{\mfz_0}{\mfz_1}$, $t \in [0,1]$, connecting $\mfz_0$ to
$\mfz_1$ and satisfying
\begin{equation}\forall\,0\leq s,t\leq 1:\quad
  \label{eq:126}
  \sfd_{\pi,\mfC}(\mfz_s,\mfz_t)=|t{-}s|\,\sfd_{\pi,\mfC}(\mfz_0,\mfz_1).
\end{equation}
If one of the two points coincides with $\mfo$, e.g.\ for $\mfz_0=\mfo,$
it is immediate to check that $\mfz_t=[x_1,tr_1]$.
If $r_0,r_1>0$ and $|x_1 {-}x_0|< \pi/2 $ 
then the unique geodesic curve
reads
(recall the convention in \eqref{eq:5})
\begin{equation}
  \label{eq:128}
\begin{aligned}  
  &r_t:=r_0\Big((1{+}t u)^2+t^2|\bfv|^2\Big)^{1/2}
  ,\quad
  x_t:=x_0+
  \bfarctan\big(\frac{t\bfv}{1{+}tu}\Big) ,
\\
& \text{where } u=\frac{r_1}{r_0}\cos(|x_1{-}x_0|)-1 \text{ and } 
    \bfv:=\frac {r_1}{r_0}\,\bfsin(x_1{-}x_0).
  \end{aligned}
\end{equation}
For example, if we operate the same construction starting from the
one-dimensional set $\Omega=[0,L]\subset \R$ with $0<L \le \pi $ we can
isometrically identify the cone space over $\Omega$ with the two-dimensional
sector $\Sigma_\Omega = \bigset{y=(r\,\cos x,\ r\,\sin x)\in \R^2}{
  r\geq0,~x\in[0,L]}$ endowed with the Euclidean distance. For
$L\in {]\pi,2\pi[}$ the identification with the sector still holds, but the
sector $\Sigma_\Omega$ is no more convex and for $x_0,x_1\in\Omega$ with
$|x_0{-}x_1|\geq \pi$ the cone distance corresponds to the geodesic distance on
the sector $\Sigma_\Omega$, i.e.\,the length of the shortest path in
$\Sigma_\Omega$ connecting two points.

On the one hand, we can define a homogeneous projection
$\mfh: \MMM_2(\mfC)\to \MMM(\OOmega)$, via
\begin{equation}
  \label{eq:frak.P}
  \mfh\lambda:=
  \sfx_\sharp(\sfr^2\lambda)=
  \int_{r=0}^\infty r^2\,\lambda(\cdot,\rmd r),
\end{equation}
i.e.\ for every $\lambda\in \MMM_2(\mfC)$ and $\zeta\in
\rmC_b(\OOmega)$ we have
$$\int_\OOmega \zeta(x)\dd (\mfh\lambda )=
\int_\mfC r^2 \zeta(x)\dd\lambda(x,r).$$
On the other hand, measures in $\MMM(\OOmega)$ can be ``lifted'' to
measures in $\MMM_2(\mfC)$, e.g.\ by considering the measure $\mu
\otimes \delta_{1} $ for $\mu\in\MMM(\OOmega)$.
More generally, for every Borel map $r:\R^d\to \left]0,\infty\right[$ and
constant $m_0\ge0$, the measure
$\lambda = m_0 \delta_{\mfo} + \mu \otimes \frac{1}{r(\cdot)^2}
\delta_{r(\cdot)}$ gives $\mfh\lambda= \mu$.

Now, the cone space formulation of the Hellinger--Kantorovich distance
between two measures $\mu_0$, $\mu_1\in\MMM(\OOmega)$ is given as
follows, see \cite[Sec.\,3]{LiMiSa16OTCR}.

\begin{theorem}[Optimal transport formulation on the cone]
  \label{thm:OTcone}
  For $\mu_0,\mu_1\in \MMM(\R^d)$ we have
  \begin{subequations}\label{eq:HK.ConeSpace.pi}
    \begin{align}    \label{eq:6}
    \HK(\mu_0,\mu_1)^2 &= \min\Big\{\sfW_{\pi,\mfC}(\lambda_0,\lambda_1)^2\,\Big|\,\lambda_i\in \mathcal
    P_2(\mfC),\ \mfh\lambda_i = \mu_i\Big\}\\    \label{eq:70}
    &=\min\Big\{\iint_{\mfC\times\mfC}
    \msd_{\pi,\mfC}(z_0,z_1)^2\rmd\bflambda(z_0,z_1)\,\Big|\,
    \mfh_i\bflambda=\mu_i\Big\},
  \end{align}
  \end{subequations}
  where $\mfh$ is defined in
  \eqref{eq:frak.P} and $\mfh_i\bflambda:=\mfh(\pi^i_\sharp\bflambda)$
  for $\bflambda\in \MMM_2(\mfC{\times}\mfC)$ and $i=0,1$.
\end{theorem}
The cone space formulation is reminiscent of classical optimal transport
problems. Here, however, the marginals $\lambda_i$ of the transport plan
$\bflambda\in\MMM(\mfC\times\mfC)$ are not fixed, and it is part of the problem
to find an optimal pair of measures $\lambda_i$ satisfying the constraints
$\mfh\lambda_i=\mu_i$ and having minimal Kantorovich--Wasserstein distance on
the cone space.

\begin{remark}[Hellinger--Kantorovich space as
  cone] \label{rem:HK.SHK} In
  {\upshape \cite{LasMie19GPCA}}
  it is shown that
  the space $(\MMM(\R^d);\HK)$ can be understood as a cone space
  over the geodesic space $(\cP(\R^d),\mathsf S\!\HK)$
  where
  the \emph{spherical Hellinger--Kantorovich distance} in $\cP(\R^d)$ 
  reads $\mathsf S\!\HK(\nu_0,\nu_1):= \arccos\big( 1-\frac12
  \HK(\nu_0,\nu_1)^2\big)$.  It would be interesting to analyze
  geodesic convexity properties of functionals $\calE$ as in
  \eqref{eq:calE} on this space; see {\upshape \cite{LasMie22?EVIH}} for a first result. 
\end{remark}

The cone space formulation in \eqref{eq:HK.ConeSpace.pi} 
reveals many interesting geometric properties of
the Hellinger--Kantorovich distance, e.g.\ Hellinger--Kantorovich geodesics
are directly connected to geodesic curves in the cone space $\mfC$, see below.
Moreover, it can be deduced that a sharp threshold exists, which 
distinguishes between transport of mass and pure \AAA growth \EEE (i.e.\ creation or
destruction) of mass.

\begin{remark}
\label{rm:dynamical}
The link between the dynamical formulation in \eqref{eq:15} and the
cone-space formulation in \eqref{eq:HK.ConeSpace.pi} of the Hellinger--Kantorovich
distance can be best seen from a Lagrangian point of view. Let\/ 
$\mathsf{Lag}_{\alpha,\beta}( X, r;V,\varrho)=
\frac{r^2}{\alpha}|V|^2+
\frac{4}{\beta} \varrho^2$
denote the rescaled Lagrangian in the definition of the 
dynamical functional \eqref{eq:15} corresponding to a curve of the form
$\mu_t:=r^2(t)\delta_{X(t)}$ and consider for fixed $r_0$, $r_1>0$ and $x_0$,
$x_1\in\OOmega$ the minimization problem
\begin{multline*}
  M_{\alpha,\beta} (x_0,r_0;x_1,r_1) :=
\min \Big\{ \!\int_0^1 \!\!\mathsf{Lag}_{\alpha,\beta}\big(X(s),r(s);\dot
X(s),\dot r(s)\big)\dd s\, 
  \Big|\\
  (X,r)\in\rmC^1\big([0,1];\OOmega\ti\R_+\big),
          X(i)=x_i,~r(i)=r_i,\Big\}.
\end{multline*}
It is not hard to check {\em \cite[Sec.\,3.1]{LiMiSa16OTCR}} that   we obtain
for $(\alpha,\beta)=(1,4)$ the explicit formula 
\[
 \HK(\mu_0,\mu_1)^2=M_{1,4} (x_0,r_0;x_1,r_1) = \msd_{\pi,\mfC}([x_0,r_0],[x_1,r_1])^2,
\]
which is the Hellinger--Kantorovich distance
of the two Dirac measures $\mu_0=r_0^2\delta_{x_0}$ and $\mu_1 = r_1^2\delta_{x_1}$
in the case that $|x_0{-}x_1|<\pi/2$.

When $|x_0{-}x_1|\ge \pi/2$,
one can always connect $\mu_0$ to $\mu_1$ by the curve
$\mu_t:=\big((1{-}t)r_0\big)^2\delta_{x_0}+t^2r_1^2\delta_{x_1}$ (whose
support is no longer concentrated on a single point) obtaining
\begin{displaymath}
  \HK(\mu_0,\mu_1)^2=2=\msd_{\pi/2,\mfC}([x_0,r_0],[x_1,r_1])^2,
\end{displaymath}
and showing the role of the threshold $\pi/2$ instead of $\pi$ in
the computation of $\HK$.
\end{remark}

The explicit computation of the previous remark is in fact a
particular case of a general result \cite[Lem.\,7.9+7.19]{LiMiSa18OETP}.
\begin{theorem}[Effective $\pi/2$-threshold in the cone distance] 
  \label{thm:effective}
  Let $\mu_0,\mu_1\in \MMM(\OOmega)$, if 
  $\bflambda \in\MMM_2(\mfC{\times}\mfC)$ is an optimal plan for the
  cone-space formulation \eqref{eq:HK.ConeSpace.pi} then
  $\bflambda\res (\mfC\ti\mfC)\setminus \{(\mfo,\mfo)\}$ is still
  optimal and 
  \begin{equation}
    \label{eq:42}
    \bflambda\Big(\big\{([x_0,r_0],[x_1,r_1])\in \mfC \ti \mfC \,\Big|
    \, r_0r_1>0 \text{ and } |x_0{-}x_1|>\frac\pi2\big\}\Big)=0,
  \end{equation}
  so that 
  \begin{subequations}\label{eq:HK.ConeSpace.pi2}
    \begin{align}    \label{eq:6bis}
    \HK(\mu_0,\mu_1)^2 &= \min\Big\{\sfW_{\pi/2,\mfC}(\lambda_0,\lambda_1)^2\,\Big|\,\lambda_i\in \mathcal
    P_2(\mfC),\ \mfh\lambda_i = \mu_i\Big\}\\    \label{eq:70bis}
    &=\min\Big\{\iint_{\mfC\times\mfC}
    \msd_{\pi/2,\mfC}(z_0,z_1)^2\rmd\bflambda(z_0,z_1)\,\Big|\,
    \mfh_i\bflambda=\mu_i\Big\}.
  \end{align}
\end{subequations}
Moreover, setting for 
$i=0,\one$
\begin{equation}
  \label{eq:33}
\begin{aligned} 
 & S_i:=\supp(\mu_i),\quad 
  S^\pit_i:= 
  \bigset{x\in { \OOmega } }{ \dist(x,S_i) <\pit },
\\
&  S_i':=S_i \cap S^\pit_{1-i}, \quad \text{ and }
  S_i'':= S_i \setminus S^\pit_{1-i} , 
\end{aligned}
\end{equation}
(see Figure \ref{fig:A'i}) with the related decomposition 
\begin{equation}
  \label{eq:43}
  \mu_i:=\mu_i'+\mu_i'',\quad
  \mu_i':=\mu_i\res S_i'=\mu_i\res S^\pit_{1-i}, \quad
  \text{ and }\mu_i'':=\mu_i\res S_i'',
\end{equation}
then we have that
\begin{subequations}
\begin{gather}
  \label{eq:48}
  \HK(\mu_0,\mu_\one)^2= \HK(\mu_0',\mu_\one')^2 + \HK(\mu_0'',\mu_\one'')^2,\\
    \label{eq:61}
    \HK(\mu_0'',\mu_\one'')^2=\mu_0''(\OOmega)+\mu_\one''(\OOmega)=
    \mu_0(\OOmega\setminus S_0')+
    \mu_\one(\OOmega\setminus S_\one').
\end{gather}
\end{subequations}
  \end{theorem}
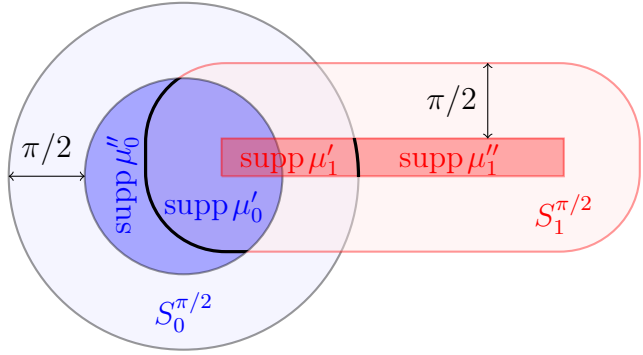
\begin{figure}[ht]
\centering \begin{tikzpicture}
 \draw[fill=blue!10, thick,opacity=0.4] (0,0) circle (2.3); 
 \draw[red, fill=red!10, thick,rounded corners=30,opacity=0.4]  (-0.5,-1)--
 (6,-1)--(6,1.5) -- (-0.5,1.5)--cycle;
 \draw[fill=blue!80, thick,opacity=0.4] (0,0) circle (1.3); 
 \draw[red, fill=red!80, thick,opacity=0.4]  (0.5,0)-- (5,0)--(5,0.5)-- (0.5,0.5)--(0.5,0);

\node[red] at (1.4,0.22) {supp$\!\;\mu'_1$};
\node[red] at (3.5,0.22) {supp$\!\;\mu''_1$};
\node[blue] at (0.4,-0.40) {supp$\!\;\mu'_0$};
\node[blue,rotate=90] at (-0.8,-0.10) {supp$\!\;\mu''_0$};

\begin{scope} 
\clip (0,0) circle (1.3);  
 \draw[rounded corners=30,very thick]  (-0.5,-1)--
 (6,-1)--(6,1.5) -- (-0.5,1.5)--cycle;
\end{scope}
\begin{scope} 
\clip  (0.5,0)-- (5,0)--(5,0.5)-- (0.5,0.5)--(0.5,0);
 \draw[rounded corners=30,very thick]  (-0.5,-1)--
 (6,-1)--(6,1.5) -- (-0.5,1.5)--cycle;
\draw[very thick] (0,0) circle (2.3); 
\end{scope}
\draw[<->] (-2.3,0)-- (-1.3,0) node[pos=0.5, above]{$\pi/2$};
\draw[<->] (4,0.5)-- (4,1.5) node[pos=0.5, left]{$\pi/2$};

\node[blue] at (0.0,-1.8) {$S^\pit_0$};
\node[red] at (5.0,-0.5) {$S^\pit_1$};
\end{tikzpicture}
\caption{The decomposition of the closed supports $S_i=\supp\mu_i$ of the measures
  $\mu_i=\mu'_i+\mu''_i$ as given in \eqref{eq:43} with cut-off at $\pi/2$. The
  open sets $S^\pit_0$ and $S^\pit_1$ denote the $\pi/2$-neighborhoods of
  the supports $S_1$ and $S_0$, respectively, and
  $\mu_i' = \mu_i\res (S^\pit_{1-i} \cap S_i)$, $\mu_i'' = \mu_i\res
  (S_i\setminus S^\pit_{1-i} )$
  are the corresponding restrictions of the measures $\mu_i$.  }
\label{fig:A'i}
\end{figure}

Note that \eqref{eq:48} shows that the decomposition in \eqref{eq:43} is
extremal with respect to the subadditivity property in Lemma~7.8 of
\cite{LiMiSa18OETP}, and \eqref{eq:61} shows that the computation of $\HK^2$
between $\mu_0''$ and $\mu_1''$ is trivial, so that
no information is lost if one restricts the evaluation of $\HK^2$
to $\mu'_0=\mu_0\res S_0'$ and $\mu'_1=\mu_1\res S_1'$.
Motivated by the above
properties, we introduce the following definition of reduced
pairs, which will play a crucial role in our analysis of geodesic curves.

\begin{definition}[Reduced pairs]
\label{def:reduced}
A pair $(\mu_0,\mu_1)\in \MMM(\OOmega)^2$ is called \emph{reduced}
(resp.\ \emph{strongly reduced}) if $\mu_i(S_i'')=0$, i.e.\,$\mu_i=\mu_i'$ for
$i=0$ and $1$ (resp.\,if $S_i\subset S^\pit_{1-i}$).
\end{definition}

By definition the sets $S_i=\mafo{supp}(\mu_i)$ are closed and $S^\pit_i$ are
open, so that $S''_i=S_i\setminus S^\pit_{1-i}$ is closed as well, but
$S'_i=S_i\cap S^\pit_{1-i}$ may be neither closed nor open. In the strongly
reduced case the condition $S_i\subset  S^\pit_{1-i}$ means that, at least
locally, the closed set $S_i$ has a positive distance to the boundary of the
open set $ S^\pit_{1-i}$. 

Notice that for every $(\mu_0,\mu_1)\in \MMM(\OOmega)^2$ the corresponding
pair $(\mu_0',\mu_1')$ defined according to \eqref{eq:33}--\eqref{eq:43} is
reduced by construction.  In fact, if $ x\in S_0'$ then there exists
$y\in \supp(\mu_1)$ with $|x{-}y|<\pi/2$: clearly $y\in S_1'$ so that
$ \dist(x,\supp(\mu_1'))\le \dist(x,S_1') < \pi/2$.
 

\subsubsection{\TTRAGRO\ maps}
\label{sss:dtm}

It is useful to express \eqref{eq:70bis} in an equivalent way, which extends
the notion of transport maps to the unbalanced case.  It relies on special
families of plans in $\bflambda\in \MMM_2(\mfC^2)$ with
$\mfh_i\bflambda = \mu_i$ generated by \TRAGRO\ systems.

\begin{definition}[\TTRAGRO\ maps]
  \label{def:dilation-transport}
  Let $\nu\in \MMM(Y)$, where $Y$ is some Polish space.  A \TRAGRO\ 
  map is a $\nu$-measurable map $(\bfT,q):Y\to X\times [0,\infty)$ with
  $q\in L^2(Y,\nu)$. It acts on $\nu$ according to this rule:
  \begin{equation}
    \label{eq:3}
    (\bfT,q)_\star \nu:=\bfT_\sharp(q^2\nu)
    =\mfh((\bfT,q)_\sharp\nu),
  \end{equation}
  where the last identity involves the obvious generalization of
  the definition \eqref{eq:frak.P} of homogeneous projection $\mfh$
  from $\MMM_2(X\times [0,\infty))$ to $\MMM(X)$.
\end{definition}
We notice that \TRAGRO\ maps obey the composition rule
\begin{equation}
  \label{eq:4}
  (\bfT_2,q_2)_\star (\bfT_1,q_1)_\star \nu=(\bfT,q)_\star \nu\quad
  \text{where}\quad
  \bfT:=\bfT_2\circ\bfT_1,\ q:=(q_2\circ \bfT_1)q_1.
\end{equation}
\TTRAGRO\ maps provide useful upper bounds for the $\HK$
metric, playing a similar role of transport maps for the
Kantorovich--Wasserstein distance.  In fact, for every choice of
 maps $(\bfT_i, q_i):Y\to \R^d\times[0,\infty)$, $i=0,1$,
associated with the measure $\nu\in \MMM(Y),$ we have
\begin{equation}
  \label{eq:95}
  \HK^2(\mu_0,\mu_1)\le 
  \int_Y  \Big(q_0^2+q_1^2-2q_0q_1
           \tcos{\pi/2}{|\bfT_0{-}\bfT_1|} \Big)\dd\nu\quad
  \mu_i:=(\bfT_i,q_i)_\star \nu.
\end{equation}
%
In order to show \eqref{eq:95} it is sufficient to check that the measure
$\bflambda\in\MMM_2(\mfC^2)$ defined by
\[
  \bflambda:=(\bfT_0,q_0;\bfT_1,q_1)_\sharp \nu,\quad
  \text{  satisfies}\quad
  \mfh_i\bflambda=\mu_i
\]
so that \eqref{eq:95} follows from
 \eqref{eq:70bis}
and the identity
\begin{equation}
  \label{eq:95bis}
  \int_{\mfC^2}
   \msd_{\pi/2,\mfC}(z_0,z_1)^2
   \dd\bflambda
  = \int_Y 
  \Big(q_0^2+q_1^2-2q_0q_1
   \tcos{\pi/2}{|\bfT_0{-}\bfT_1|}\Big) \dd\nu.
\end{equation}

On the other hand, choosing $Y=\mfC \ti\mfC$ and an optimal plan
$\nu=\bflambda\in \MMM_2(\mfC \ti\mfC)$  for
 \eqref{eq:70bis}

and setting 
$\bfT_i([x_0,r_0],[x_1,r_1]):=x_i$ and $q_i([x_0,r_0],[x_1,r_1])=r_i$,
we immediately find 
\begin{equation}
  \label{eq:95ter} 
  \HK^2(\mu_0,\mu_1)=
  \int_{\mfC \ti\mfC} 
  \Big(q_0^2+q_1^2-2q_0q_1
  \tcos{\pi/2}{|\bfT_0{-}\bfT_1|}
  \Big)\dd\bflambda,\quad
  \mu_i:=(\bfT_i,q_i)_\star \bflambda,
\end{equation}
and therefore equality holds in \eqref{eq:95}.

\begin{corollary}[$\HK$ via \TRAGRO\ maps]
\label{cor:HK-dt}
  For every $\mu_0,\mu_1\in \MMM(\R^d)$ we have
  \begin{equation}
    \label{eq:101}
    \begin{aligned}
      \HK^2(\mu_0,\mu_1)= 
      \min\Big\{&\int_{\mfC \ti\mfC} 
      \Big(q_0^2+q_1^2 -2q_0q_1
      \tcos{\pi/2}{|\bfT_0{-}\bfT_1|}\Big)
      \dd\bflambda
      \Big|\, \bflambda\in \MMM( Y), 
      \\
      &Y\text{Polish,\ }(\bfT_i,q_i):Y 
      \to \R^d\ti[0,+\infty),\ 
      \mu_i:=(\bfT_i,q_i)_\star \bflambda \Big\};
    \end{aligned}
  \end{equation}
  moreover, it is not restrictive to choose $Y=\mfC \ti\mfC$ in \eqref{eq:101}.
\end{corollary}

Inspired by the so-called Monge formulation of Optimal Transport, it is natural
to look for similar improvement of \eqref{eq:101}, when $Y=\R^d$, $\nu=\mu_0$,
$\bfT_0(x_0)=x_0$ is the identity map, and $q(x_0)\equiv 1$.
\begin{problem}[Monge formulation of $\HK$ problem]
\label{prob:Monge}
Given $\mu_0,\mu_1\in \MMM(\OOmega)$
such that $\mu_1=\mu_1'$, $\mu_1''=0$
(recall \eqref{eq:33} and \eqref{eq:43}),
find an optimal \TRAGRO\ pair
$(\bfT,q):\R^d\to \R^d\times[0,\infty)$ minimizing the cost
\begin{equation}
  \label{eq:24}
  \Mcost \bfT q{\mu_0}:=\int_{\OOmega} 
  \Big(1+q^2(x)-2q(x)
  \tcos{\pit}{|\bfT(x){-}x|}\Big)\dd\mu_0(x)
\end{equation}
among all the \TRAGRO\ maps satisfying $(\bfT,q)_\star \mu_0=\mu_1$
\end{problem}
By \eqref{eq:95} we have the  bound
\begin{equation}
  \label{eq:9}
  \begin{aligned}
    \HK(\mu_0,\mu_1)^2 \leq    \inf
    \Big\{\Mcost \bfT q{\mu_0}\Big|
  (\bfT,q)_\star \mu_0=\mu_1\Big\} .
  \end{aligned}
\end{equation}
When $\mu_0 \ll \Leb d$ and the support of $\mu_1$ 
is contained in the closed neighborhood of radius 
$\pi/2$ of the support of $\mu_0$,
the results of the next section (cf.\ Corollary
\ref{cor:Monge}), which are a consequence of the 
optimality conditions in Theorem~\ref{thm:optimality.cond},
show that the minimum of Problem \ref{prob:Monge} is attained and realizes the
equality in \eqref{eq:9}.%

\subsubsection{Entropy-transport problem}
\label{suu:ET.Prob} 

A third point of view, typical of optimal transport problems,
characterizes the Hellinger--Kantorovich distance via the static
Logarithmic Entropy Transport (LET) variational formulation.

We define the logarithmic entropy density
$F:\left[0,\infty\right[\to[0,\infty[$ via
\[
F(s) := s\log s - s + 1\quad\text{for $s>0$}\quad \text{and} \quad  F(0):=1,
\]
and the cost function $\Ell_1:\R^d\to[0,\infty]$ via
\begin{equation}\label{eq:Ell}
  \Ell_1(x) := \frac12\ell(|x|),\quad
  \ell(r):=
  \begin{cases}
    \displaystyle
    -\log(\cos^2(r))=\log\big(1{+}\tan^2(r)\big)&\text{for }r <\pi/2,\\
      +\infty&\text{otherwise.}
  \end{cases}
\end{equation}
For given $\mu_0,\mu_1 \in \MMM(\OOmega)$ the entropy-transport
functional $\calET(\,\cdot\,;\mu_0,\mu_1):
\MMM(\OOmega\ti\OOmega)\to[0,\infty]$ reads
\begin{equation}
  \label{eq:ET.sig.eta}
  \calET(\bfeta;\mu_0,\mu_1):=
  \int_\OOmega F(\sigma_0)\dd\mu_0+ \int_\OOmega F(\sigma_1)\dd\mu_1+
   \iint_{\OOmega\times\OOmega}2\Ell_1(x_0{-}x_1)\dd\bfeta 
\end{equation}
with $(\pi_i)_\sharp \bfeta=\sigma_i\mu_i\ll\mu_i$.  As usual, we set
$\calET(\bfeta;\mu_0,\mu_1):=+\infty$ if one of the marginals
$(\pi_i)_\sharp \bfeta$ of $\bfeta$ is not absolutely continuous with respect
to $\mu_i$. With this definition, the equivalent formulation of the
Hellinger--Kantorovich distance as entropy-transport problem reads as follows.

\begin{theorem}[LET formulation]
  \label{thm:LET}
  For every $\mu_0,\mu_1\in \MMM(\OOmega)$ we have
  \begin{equation}\label{eq:HK.ET}
\HK(\mu_0,\mu_1)^2=\min\Big\{\calET(\bfeta;\mu_0,\mu_1)\,\big|\,
\bfeta\in\MMM(\OOmega\times\OOmega)
\Big\}.
\end{equation}
Moreover, recalling the decomposition \eqref{eq:33}--\eqref{eq:43},
\begin{enumerate}[{\upshape(1)}]
\item
  \label{thm:LET.lab1}
      the pairs $(\mu_0,\mu_\one) $ and $(\mu'_0,\mu'_\one)$ 
share
      the same optimal plans $\bfeta$
\item
  \label{thm:LET.lab2}
  if we set $g_0(x_0):=([x_0,1],\mfo)$ and
  $g_1(x_1):=(\mfo,[x_1,1])$,
  every optimal plan $\bfeta\in \MMM(\OOmega\ti\OOmega)$ for the
  entropy-transport formulation in \eqref{eq:HK.ET} induces optimal
  plans $\bfbeta$ (resp.~$\bfbeta'$) in $\MMM( \mfC\ti\mfC)$ for the
  pair $(\mu_0,\mu_1)$ (resp.\ the reduced pair $(\mu_0',\mu_1')$) via
  \begin{equation}
    \label{eq:65}
    \bfbeta':=(x_0,\sigma_0^{-1/2};x_1,\sigma_1^{-1/2})_\sharp
    \bfeta,\quad
    \bfbeta:=\bfbeta'+(g_0)_\sharp\, \mu''_0
    +(g_1)_\sharp \, 
    \mu_1''.
  \end{equation}
\end{enumerate}
\end{theorem}

An optimal transport plan $\bfeta$, which always exists, gives the effective
transport of mass. Note, in particular, that the finiteness of $\calET$ only
requires $(\pi_i)_\sharp \bfeta=\eta_i\ll\mu_i$ (which is considerably weaker
than the usual transport constraint $(\pi_i)_\sharp \bfeta=\mu_i $) and the
cost of a deviation of $\eta_i$ from $\mu_i$ is given by the entropy
functionals associated with $F$. Moreover, the cost function $\ell$ is finite
in the case $|x_0{-}x_1|<\pi/2$, which highlights the sharp threshold between
transport and pure creation/destruction.  Notice that we could equivalently use
the truncated function $\tcosq{\pit}r=\cos^2(\min\{r ,\pit\})$ instead of
$\cos^2(r)$ in \eqref{eq:Ell}. As we have already seen,
the function $r \mapsto \tcosq{\pit}r$
plays an important role in many formulae.

In general, optimal entropy-transport plans
$\bfeta\in\MMM(\OOmega\ti\OOmega)$ are not unique. However, due to the
strict convexity of $F$, their marginals $\eta_i$ are unique so that the
non-uniqueness of the plan $\bfeta$ is solely a property of the optimal
transport problem associated with the cost function $(x_0,x_1) \mapsto
2\Ell_1(x_1{-}x_0)= \ell\big(|x_1{-}x_0|\big)$.

\begin{remark}
\label{rm:Perspective}
Besides \eqref{eq:65},
the connection between the cone-space formation and the logarithmic
entropy-transport problem is given by the homogeneous marginal perspective
function, namely
\begin{align*}
  \msd_{\pi/2,\mfC}
  ([x_0,  r_0],[x_1,r_1])^2  =\inf \big\{ 
  r_0^2 F(\tfrac{\theta}{r_0^2}) + r_1^2F\big(\tfrac{\theta}{r_1^2}\big) 
  +2\theta \Ell_1(x_0{-}x_1) \,\big|\, \theta > 0\big\},
\end{align*} 
where $r_i^2$ plays the role of the reverse densities $1/\sigma_i$ and $\theta$
is a scaling parameter, see \cite[Sec.\,5]{LiMiSa18OETP}.

We highlight that the logarithmic entropy-transport formulation
\eqref{eq:HK.ET} can be easily generalized by considering convex and lower
semi-continuous functions $F_0$ and $F_1$ and cost functions $\ell$, see
\cite[Part I]{LiMiSa18OETP}.
\end{remark}
Applying the previous Theorem~\ref{thm:LET} we can refine formula
\eqref{eq:95bis} by providing an optimal pair of \TRAGRO\ maps solving
\eqref{eq:101} in the restricted set
$Y=S_0\times S_1\subset \OOmega\ti\OOmega$. Indeed, we can choose arbitrary
points $\bar x_i\in S_i$ and
\begin{equation}
  \begin{aligned}
    \bfnu:={}&\bfeta+ \mu_0'' \oti
    \delta_{\bar x_1}+\delta_{\bar x_0}\oti\mu_1'',\\
    \bfT_i(x_0,x_1):={}&x_i,\quad q_i(x_0,x_1):=
    \begin{cases}
      \sigma_i^{-1/2}(x_i)&\text{if }(x_0,x_1)\in S_0'\ti S_1',\\
      1&\text{if }(x_0,x_1)\in (S_0\ti S_1)\setminus (S_0'\ti S_1')
    \end{cases}
  \end{aligned}
  \label{eq:98}
\end{equation}
which satisfies
\begin{equation}
  \label{eq:92}
  (\bfT_i,q_i)_\star \bfnu =\mu_i,\quad  \HK^2(\mu_0,\mu_1)=
  \int_Y \Big(q_0^2+q_1^2-2q_0q_1\tcos{\pit}{|\bfT_0{-}\bfT_1|}\Big)\dd\bfnu.
\end{equation}

\subsubsection{Dual formulation with Hellinger--Kantorovich potentials}
\label{ss:dualformulation}

In analogy to the Kantorovich--Wasserstein distance, we can give a dual
formulation in terms of Hellinger--Kantorovich potentials.  We slightly modify
the notation of \cite{LiMiSa18OETP}, in order to be more consistent with the
approach by the Hamilton--Jacobi equations (and the related Hopf--Lax solutions)
of Section \ref{sec:HJ} and to deal with rescaled distances. As we will
study segments of constant-speed geodesics $t\to \mu_t$ of length $\tau= t{-}s$ for $0\leq
s < t \leq 1$, it will be convenient to introduce a scaling parameter
$\tau>0$ that in certain parts will be replaced by $1$, namely if we consider a whole
geodesic. With this parameter, we set
\begin{equation}
  \label{eq:30}
  F_\tau(s):=\frac 1{2\tau} F(s),\quad
  \Ell_\tau(x)=\frac 1{2\tau} \ell(|x|),
  \quad
  \calET_\tau(\bfeta;\mu_0,\mu_\one)=\frac 1{2\tau}\calET(\bfeta;\mu_0,\mu_\one)
\end{equation}
and the corresponding
\begin{equation}
  \label{eq:HK.ETtau}
  \frac 1{2\tau}\HK^2(\mu_0,\mu_\one)=\min
  \Big\{\calET_\tau(\bfeta;\mu_0,\mu_\one)
  \,\big|\, \bfeta \in \MMM(\OOmega \ti \OOmega)
  \Big\}.
\end{equation}
It is clear that minimizers $\bfeta$ of \eqref{eq:HK.ETtau} are
independent of the coefficient $\frac 1{2\tau}$ in front of $\HK$
and coincide with solutions to \eqref{eq:HK.ET} if
$\mu_\tau=\mu_1$. The role of $\tau$ just affects the rescaling of the
potentials $\varphi$ and $\xi$ we will introduce below.

We also introduce the Legendre transform of $F_\tau$
\begin{equation}
  \label{eq:20}
  \chGtrafo_\tau(\varphi):= F_\tau^* (\varphi)
  =  \sup_{s>0} \varphi s-F_\tau(s)=
  \frac{\ee^{2\tau
      \varphi}-1}{2\tau}
 ,\quad
  \Gtrafo_\tau(\varphi):=\frac{1-\ee^{-2\tau
      \varphi}}{2\tau}=-\chGtrafo_\tau(-\varphi),
\end{equation}
extended to $[-\infty,+\infty]$ by
\begin{equation}
  \label{eq:76}
  \Gtrafo_\tau(+\infty)=-\chGtrafo_\tau(-\infty)=\frac 1{2\tau},\quad
  \Gtrafo_\tau(-\infty)=-\chGtrafo_\tau(+\infty)=+\infty,
\end{equation}
and their inverses
\begin{equation}
  \label{eq:21}
  \cIGtrafo_\tau(\psxi):=\frac1{2\tau}\log (1{+}2\tau\psxi),\quad
  \IGtrafo_\tau(\psxi):=
  -\frac 1{2\tau}\log(1{-}2\tau\psxi)=-\cIGtrafo_\tau(-\psxi),
\end{equation}
defined for $\psxi\in [-\frac1{2\tau},+\infty]$ and
$\psxi\in[-\infty,\frac 1{2\tau}]$ respectively, with the obvious convention
induced by \eqref{eq:76}.  
With Theorem 6.3 in \cite{LiMiSa18OETP} (see also Section 4 therein), we have the equivalent characterization of $\HK$
via the dual formulation
\begin{subequations}
\label{eq:HK.dual}
\begin{align}
  \nonumber
  \frac 1{2\tau}
  \HK(\mu_0,\mu_\one)^2&=\sup\Big\{\int_\OOmega \Gtrafo_\tau(\varphi_\one)\dd\mu_\one
  -\int_\OOmega \chGtrafo_\tau(\varphi_0)
  \dd\mu_0\,\Big|
  \,\\&\hspace{5em}\varphi_0,\varphi_\one\in\rmC_b(\OOmega),~
  \varphi_\one(x_\one){-}\varphi_0(x_0)\leq \Ell_\tau(x_\one{-}x_0)\Big\}
  \label{eq:HK.dual.a}
  \\
  \nonumber&=
  \sup\Big\{\int_\OOmega
  \psxi_\one\dd\mu_\one-\int_\OOmega 
  \psxi_0\dd\mu_0\,\Big|\,
  \psxi_i\in\rmC_b(\OOmega),~\sup_\OOmega\psxi_\one<\frac1{2\tau},\
  \inf_\OOmega\psxi_0>- \frac1{2\tau}
  \\
  &\hspace{7em}
  \big(1{-} 2\tau\psxi_\one(x_\one)\big)
  \big(1{+}2\tau\psxi_0(x_0)\big)
  \geq \tcosq{\pit}{|x_0{-}x_\one|}
  \Big\}.
  \label{eq:HK.dual.b}
\end{align}
\end{subequations}
Note that the formulations in \eqref{eq:HK.dual.a} and \eqref{eq:HK.dual.b} are
connected by the transformation
$ \psxi_\one= \Gtrafo_\tau(\varphi_\one),\ \psxi_0=\chGtrafo_\tau(\varphi_0)$ and the last
condition in \eqref{eq:HK.dual.b} is equivalent to
\begin{equation}
  \label{eq:22}
  \IGtrafo_\tau
  \big(\psxi_\one(x_\one)\big){-}
  \cIGtrafo_\tau
  \big(\psxi_0(x_0)\big)\leq
  \Ell_\tau(x_\one{-}x_0).
\end{equation}
It is not difficult to check that one can also consider Borel
functions in \eqref{eq:HK.dual.a} and \eqref{eq:HK.dual.b}, e.g.\ for
all Borel functions $\varphi_i:\R^d \to [-\infty,+\infty]$ with
\begin{equation}
  \label{eq:condphi}
  \begin{gathered}
    \int_{\R^d}\rme^{-2\tau \varphi_\one}\dd\mu_\one<\infty,\quad
    \int_{\R^d}\rme^{2\tau \varphi_0}\dd\mu_0<\infty,\\
    \varphi_\one(x_1) \leq \Ell_\tau(x_1{-}x_0) +\varphi_0(x_0)
    \quad\text{for all }
    x_0,x_\one\in \R^d\text{ with }|x_0{-}x_\one|<\pi/2,
  \end{gathered}
\end{equation}
we have
\begin{equation}
  \label{eq:73}
  \frac 1{2\tau}\HK (\mu_0,\mu_\one)^2
  \ge \int_\OOmega \Gtrafo_\tau(\varphi_\one)\dd\mu_\one
  -\int_\OOmega \chGtrafo_\tau(\varphi_0)\dd\mu_0.
\end{equation}
  If we allow extended valued Borel functions, the supremum in
\eqref{eq:HK.dual.a} and \eqref{eq:HK.dual.b} are attained.

\begin{theorem}
[Existence of optimal dual pairs]%
\label{thm:Borel-duality}\!\!%
For all $\mu_0,\mu_\one\in \MMM(\R^d)$ and $\tau>0$ there exists an optimal
pair of Borel potentials $\varphi_0,\varphi_\one:\R^d\to[-\infty,+\infty]$
which is admissible according to \eqref{eq:condphi} and realizes equality
in \eqref{eq:73}, namely 
\begin{equation}
  \label{eq:73.equality}
  \frac 1{2\tau}\HK (\mu_0,\mu_\one)^2
  =\int_\OOmega \Gtrafo_\tau(\varphi_\one)\dd\mu_\one-
  \int_\OOmega \chGtrafo_\tau(\varphi_0)\dd\mu_0.
\end{equation}
The transformations
$\psxi_0:=\chGtrafo_\tau(\varphi_0):\R^d\to [ -1/(2\tau), +\infty]$, and 
$\psxi_\one:=\Gtrafo_\tau(\varphi_\one):\R^d\to [-\infty,1/(2\tau) 
]$, give an optimal pair for \eqref{eq:HK.dual.b} (dropping $\xi_i\in
\rmC_b(\R^d)$) satisfying 
\begin{align}
  \label{eq:8bis}
  \int_{\R^d}|\xi_i|\dd\mu_i&<\infty,\quad i=0,\one,\\
  \label{eq:8}
  (1{-}2\tau\psxi_\one(x_\one))( {1{+}2\tau\psxi_0(x_0)})&\ge
  \tcosq{\pit}{|x_0{-}x_\one|}
  \quad\text{if }\psxi_0(x_0)<\infty, 
  \ \psxi_\one(x_\one) > - \infty,\\
  \label{eq:10}
  \frac 1{2\tau} \HK(\mu_0,\mu_\one)^2&=\int_{\OOmega}\psxi_\one\dd\mu_\one-
  \int_{\OOmega} \psxi_0\,\dd\mu_0.
\end{align}
\end{theorem}

\begin{remark}
\label{rem:tedious}\upshape 
Denoting by $S_i: =\supp(\mu_i)$ the support of $\mu_i$ for $i=0$ and $1$, 
we remark that it is always sufficient to find Borel potentials
$\varphi_i:S_i\to [-\infty,+\infty]$ satisfying \eqref{eq:condphi} on
$S_0\ti S_1$ instead of $\R^d \ti \R^d$. By setting $\tilde \varphi_1:=-\infty$
in $\R^d\setminus S_1$ and $\tilde \varphi_0:=+\infty$ in $\R^d\setminus S_0$
we obtain a pair still satisfying \eqref{eq:condphi} and
\eqref{eq:73.equality}. This freedom will be useful in Theorem~\ref{thm:optimality.cond} below. 
\end{remark}

Moreover, notice that \eqref{eq:HK.dual.b} can be rewritten as
\[
 \frac 1{2\tau}\HK (\mu_0,\mu_\one)^2
  =
      \sup\bigg\{
      \int_\OOmega \mathscr P_{\kern-2pt\tau}\xi_0 \dd\mu_\one  - 
   \int_\OOmega\xi_0\dd\mu_0\,\Big|\,
      \xi_0\in \rmC_b(\OOmega),\ \xi_0 > - \frac1{2\tau} 
      \bigg\},          
    \]
where $\HopfLax\tau\xi$ is defined in \eqref{eq:23intro}. 
In particular, the operator $\mathscr P_\tau$ is directly connected to the dynamical formulation in \eqref{eq:15},
and we will thoroughly study its properties in Section \ref{sec:HJ}.

\subsection{\texorpdfstring{First order optimality for $\HK$}
{First order optimality for HK}}
\label{su:1st.Optim}

From the above discussion, we have already seen that there is never any
transport over distances larger than $\pi/2$. This transport bound
will also be seen in the following optimality conditions for the
marginal densities $\sigma_i$ defined in \eqref{eq:ET.sig.eta}. 

\begin{theorem}[Optimality conditions {\cite[Thm.\ 6.3]{LiMiSa18OETP}}]
\label{thm:optimality.cond} 
Let $\mu_0,\mu_\one\in\MMM(\OOmega)$
and let $S_i,S_i',S_i'',\mu_i'$ be defined as in
\eqref{eq:33}--\eqref{eq:43}.
The following holds:
\begin{enumerate}[{\upshape(1)}]
\item  \label{thm:optim.cond.lab2} \label{th:OptiCond.labelii} A plan $\bfeta\in\MMM(\OOmega\times\OOmega)$
  is optimal for the logarithmic entropy-transport problem in
  \eqref{eq:HK.ETtau} if and only if
  \begin{itemize}
  \item[-] $\iint\ell\dd\bfeta<\infty$

  \item[-] its marginals $\eta_i$ are absolutely continuous with respect to
    $\mu_i'$ (equivalently, $\eta_i$ are absolutely continuous with respect to
    $\mu_i$ and $\eta_i(S_i'')=0$),

  \item[-] there exist Borel densities $\sigma_i: \R^d\to [0,\infty]$ such that
    $ \eta_i=\sigma_i \mu_i'$ and
    \begin{subequations}
      \begin{align}
        \label{eq:1}
        \sigma_i=0\quad &\text{on }S_i'',
        \\
        \label{eq:sigma.g.0}     
        0<\sigma_i<\infty\quad&\text{on }S_i',\\
        \sigma_i=+\infty\quad&\text{on }\R^d\setminus S_i,\\
        \sigma_0(x_0)\sigma_\one(x_\one)\geq
        \tcosq{\pit}{|x_0  {-}x_\one|}\quad&\text{on }  S_0\times S_\one,
        \\
        \label{eq:34}
        \sigma_0(x_0)\sigma_\one(x_\one)= \tcosq{\pit}{|x_0 {-}x_\one|}
        \quad&\bfeta\text{-a.e.\,on }
        S_0\times S_\one.
      \end{align}
   \end{subequations}
   \end{itemize}
   In particular, the marginals $\eta_i$ are unique and the densities
   $\sigma_i$ are unique $\mu'_i$-a.e. 

 \item \label{thm:optim.cond.lab3} If $\bfeta$ is optimal and
   $S_i,\,S_i',\,S_i''$ and $\sigma_i$ are defined as above, the pairs of
   potentials defined by
\begin{align}
  \label{eq:19}
  \varphi_\one:={}&\begin{cases}
    -\frac 1{2\tau}\log \sigma_\one&\text{in }S_\one',\\
    +\infty&\text{in }S_\one'',\\
    -\infty&\text{in }\OOmega\setminus S_\one;
  \end{cases}
  & \varphi_0:={}&\begin{cases} \frac 1{2\tau} \log\sigma_0
    &\text{in }S_0',\\  -\infty&\text{in }S_0'',\\
    +\infty&\text{in }\OOmega\setminus S_0;
    \end{cases}
    \\   \label{eq:19xi}
    \psxi_\one:={}&  \begin{cases} \frac{1{-}\sigma_\one}{2\tau}
      &\text{in }S_\one',\\  \frac 1{2\tau}&\text{in }S_\one'',\\
      -\infty&\text{in }\OOmega\setminus S_\one;
    \end{cases}
    &    \psxi_0:={}& \begin{cases}
      \frac{\sigma_0{-}1}{2\tau} 
      &\text{in }S_0',\\
      -\frac 1{2\tau}&\text{in }S_0'',\\
      +\infty&\text{in }\OOmega\setminus S_0;
    \end{cases}
\end{align}
are optimal in the respective dual relaxed characterizations 
of Theorem~\ref{thm:Borel-duality} 
and satisfy $\bfeta$-a.e.\,in $\R^d\times \R^d$
\begin{subequations}
  \begin{align}
    \varphi_i(x_i)&\in \R, &
     \varphi_\one(x_\one)-\varphi_0(x_0)&= \Ell_\tau(x_\one{-}x_0),
                                 \label{eq:2}\\
    -\psxi_0(x_0),\psxi_\one(x_\one)&\in \big(\frac1{2\tau},\infty\big),\kern-6pt
  &  (1{+}2\tau\psxi_0(x_0))(1{-}2\tau\psxi_\one(x_\one)) 
   &=  \tcosq{\pit}{|x_0{-}x_\one|}.
   \label{eq:2bis}
 \end{align}
\end{subequations}

\item \label{thm:optim.cond.lab4} Conversely, if $\bfeta$ is optimal and
  $(\varphi_0,\varphi_\one)$ (resp.\,$(\xi_0,\xi_\one)$) is an optimal pair
  according to Theorem~\ref{thm:Borel-duality}, then \eqref{eq:2}
  (resp.\,\eqref{eq:2bis}) holds $\bfeta$-a.e.\,and
\begin{equation}
    \label{eq:18}
  \begin{aligned}
    \sigma_\one={}&\ee^{-2\tau\varphi_\one}=1{-}2\tau\xi_\one
    \kern-4pt&\text{$\mu_\one$-a.e.\,in }S_\one',
   &&
    \varphi_\one={}&+\infty, \ \xi_\one=\frac1{2\tau}&\text{$\mu_\one$-a.e.\,in }S_\one'',
    \\
    \sigma_0={}&\ee^{2\tau\varphi_0}=1{+}2\tau\xi_0
        \kern-4pt&\text{$\mu_0$-a.e.\,in }S_0',
    &&
    \varphi_0={}&-\infty, \ \xi_0=-\frac{1}{2\tau}&\text{$\mu_0$-a.e.\,in }S_0''.
  \end{aligned}  
\end{equation}
\end{enumerate}
\end{theorem}

\section{Regularity of static 
\texorpdfstring{$\HK$}{HK} potentials 
\texorpdfstring{$\varphi_0$}{phi0} and
\texorpdfstring{$\varphi_1$}{phi1}} 
\label{se:2nd.Optim}

In this section, we will carefully study the regularity of a pair
$(\varphi_0,\varphi_1)$ of optimal $\HK$ potentials arising in \eqref{eq:19} of
Theorem~\ref{thm:optimality.cond}.  We will improve the previous approximate
differentiability result of \cite[Thm.\ 6.6(iii)]{LiMiSa18OETP} (see also
\cite[Thm.\,6.2.7]{AmGiSa08GFMS}) by adapting the argument of
\cite{FigGig11LSKP} and extending the classical result of \cite{GanMcc96GOT} to
the $\HK$ setting. In fact, this section is largely independent of the specific
$\HK$ setting but relies purely on the theory of $\LLL$-transforms. As we are
interested in the special case of \GGG continuous, extended values cost
functions \EEE $\LLL = \Ell_\tau=\frac1{\tau} \Ell_1:\R^d\to[0,+\infty]$ which
attain the value $+\infty$ outside a ball, we cannot rely on existing results
and have to provide a careful analysis of this case (but see also \GGG
\cite{GO07,MM09,JS12,BP13,BPP18} \EEE for different situations of discontinuous
costs taking the value $+\infty$).

We will use the notion of locally semi-concave and semi-convex functions; recall
that a function $\varphi:U\to \R$ defined in some open set $U$ of $\R^d$ is
locally semi-concave if for every point $\bar x\in U$ there exists $\rho>0$ and
a constant $C>0$ with
\begin{equation}
  \label{eq:89}
  x\mapsto \varphi(x)-\frac C2 |x|^2 \quad\text{is concave in
  }B_\rho(\bar x).
\end{equation}
A function $\varphi$ is locally semi-convex if $-\varphi$ is locally
semi-concave.
Let us recall that 
locally semi-concave functions are locally Lipschitz and thus
differentiable almost everywhere. We will denote by $\mathrm{dom}(\nabla
\varphi)$ the domain of their differential. By Alexandrov's Theorem
(see~\cite[Thm.\ 5.5.4]{AmGiSa08GFMS}),  there
exists for almost every $x\in \mathrm{dom}(\nabla\varphi)$ a symmetric matrix $\mathsf A=:\rmD^2 \varphi( x)$ such that
\begin{subequations}
 \label{eq:Hess.sfA}
\begin{align}
  \label{eq:72}
  &\lim_{y\to x}\frac{\varphi(y)-\varphi(x)-\langle \nabla 
   \varphi(x),y{-}x\rangle-\frac
    12 \langle \mathsf A(y{-}x),y{-}x\rangle}{|y{-}x|^2}=0,
\\[0.5em]
&\text{and} \quad
  \label{eq:72bis}
  \lim_{\underset{y\in \mathrm{dom}(\nabla\varphi)}{y\to x}}\frac{\nabla\varphi(y)-\nabla\varphi(x)-
    \mathsf A(y{-}x)}{|y{-}x|}=0.
\end{align}
\end{subequations}
We will denote by $\mathrm{dom}(\rmD^2\varphi)$ the subset of density points in
$\mathrm{dom}(\nabla\varphi)$ where \eqref{eq:72} and \eqref{eq:72bis}
hold.

As the optimality of potential pairs $(\varphi_0,\varphi_1)$ is closely related
to the theory of $\LLL$-transforms, we give the basic definitions first and then
derive the associated regularity properties under additional smoothness
assumptions.

For simplicity, we restrict the analysis of the remaining text to continuous functions
$\LLL : \R^d \to [0,\infty]$ satisfying
$\mathrm{dom}(\LLL) = \bigset{ z \in \R^d}{ \LLL(z)\in \R}
=B_{\RRR}(0)$ for some ${\RRR}>0$,
i.e.\ $\LLL(z)<\infty$
for $|z|<\RRR$ and $\LLL(z)=+\infty$ for $|z|\geq \RRR$.
By continuity of
$\LLL$ this behavior implies $\LLL(z_k)\to +\infty$ if
$\liminf_{k\to\infty} |z_k| \geq \RRR$.

We define the \emph{forward  $\LLL$-transform} $\varphi_0^{\LLL\to}$ of a
l.s.c.\ function $\varphi_0$ and the \emph{backward $\LLL$-transform}
$\varphi_1^{\shortleftarrow\LLL}$ of an u.s.c.\ function $\varphi_1$ via 
\begin{equation}
  \label{eq:Def.Ell.Trafo}
  \begin{aligned} 
    \varphi_0^{\LLL\to}(x_1) &:= \inf_{x_0\in B_\RRR(x_1)
    }
    \varphi_0(x_0) + \LLL(x_1{-}x_0) \quad \text{and}
    \\
    \varphi_1^{\shortleftarrow\LLL}(x_0) &:= \sup_{x_1 \in B_\RRR(x_0)
    } \varphi_1(x_1) -\LLL(x_1{-}x_0),
  \end{aligned}
\end{equation}
where the restriction of the infimum and supremum in
\eqref{eq:Def.Ell.Trafo} to 
the balls $B_\RRR(x_i)$, corresponding to the shifted proper domain of
$\LLL$,
is important to avoid
the expression ``$\infty- \infty$''. It will turn out that
$ \varphi_0^{\LLL\to}$ is u.s.c.\ and $\varphi_1^{\shortleftarrow\LLL}$ is
l.s.c. Of course, these transformations are related by
\begin{equation}
  \label{eq:156}
    \varphi^{\LLL\to}_0(x)=- 
({-}\varphi_0)^{\shortleftarrow\LLL}(x),
\end{equation}
and for arbitrary functions $\psi_i:\R^d \to [-\infty,+\infty]$ we have
the general relations  
\begin{equation}
  \label{eq:forw.backw.Ell}
  \psi_0^{\LLL\to} = \Big(\big(\psi_0^{\LLL\to}\big)^{\shortleftarrow\LLL}
  \Big)^{\LLL\to} \quad \text{and} \quad  \psi_1^{\shortleftarrow\LLL} = 
 \Big(\big(\psi_1^{\shortleftarrow\LLL}\big)^{\LLL\to}\Big)^{\shortleftarrow\LLL},
\end{equation}
see \cite[Ch.\,5]{Vill09OTON}.  For later usage, we consider the following
elementary example.

\begin{example}[Forward and backward $\LLL$-transform]
\slshape\label{ex:EllTransform}
We consider the potentials 
\[
\varphi_0(x_0) = \begin{cases} a_0 &\text{for } x_0=y_0, \\ +\infty& \text{otherwise},
\end{cases} \quad \text{ and }\quad  \varphi_1(x_1)=  
\begin{cases} a_1 &\text{for } x_1=y_1, \\ -\infty& \text{otherwise},
\end{cases}
\]
where $-\infty \leq a_0 <+\infty$, $-\infty < a_1 \leq +\infty$
and $y_0,y_1\in \R^d$ are fixed. For
$a_0,a_1 \in \R$ we find the transforms 
\begin{align*}
\varphi^{\LLL\to}_0(x_1) &= \begin{cases}\! a_0{+}\LLL(x_1{-}y_0)
  &\text{for } x_1 \in  B_\RRR (y_0), \\ 
  +\infty& \text{otherwise},
\end{cases} \\
\varphi^{\shortleftarrow\LLL}_1(x_0) &= \begin{cases}\! a_1{-}\LLL(y_1{-}x_0)
  &\text{for }x_0 \in B_\RRR(y_1), \\ 
  -\infty& \text{otherwise}.
\end{cases}
\end{align*}
For $a_0=-\infty $ and $a_1 = + \infty$, we obtain the transforms 
\[
\varphi^{\LLL\to}_0(x_1) = \begin{cases}\! -\infty
  &\text{for }x_1  \in B_\RRR (y_0), \\ +\infty& \text{otherwise},
\end{cases}\  \text{ and } 
\varphi^{\shortleftarrow\LLL}_1(x_0) = \begin{cases}\! +\infty
  &\text{for } x_0 \in B_\RRR (y_1), \\ -\infty& \text{otherwise},
\end{cases}
\]
As $B_\RRR (y_i)$ is open, we see that
$\varphi^{\LLL\to}_0$ is u.s.c.\ and $\varphi^{\shortleftarrow\LLL}_1$ is
l.s.c.
  Moreover, observe that
$\big(\varphi^{\LLL\to}_0\big)^{\shortleftarrow\LLL} = \varphi_0$ and
$\big(\varphi^{\shortleftarrow\LLL}_1\big)^{\LLL\to} = \varphi_1$, so that
\eqref{eq:forw.backw.Ell} is true for $\psi_0\in \big\{ \varphi_0,
\varphi_1^{\shortleftarrow\LLL} \big\}$  and
$\psi_1 \in \big\{\varphi_1, \varphi_0^{\LLL\to} \big\}$, respectively. 
\end{example}

For $\RRR>0$ and sets $S\subset \R^d$, we
introduce the notation
\begin{equation}
  \label{eq:108}
  \begin{gathered}
    S^{\RRR}:=\bigset{x\in \R^d }{ \dist(x,S)< \RRR },\\[0.3em]
    \Ext{\RRR}{S}:=\bigcup_{x:\ \dist(x,S)> \RRR}
    B_{\RRR}(x) ,\quad \Bdry{\RRR}S:=\partial S\cap \partial\big(\Ext{\RRR}S\big).
  \end{gathered}
\end{equation}
In particular, $\Ext{\RRR}S$ is the open subset of $\R^d\setminus S$ obtained by
taking the union of all the open balls of radius ${\RRR}$ that do not intersect
$S$.  If $S$ is closed and satisfies an exterior sphere condition of radius
${\RRR}$ at every point of its boundary (e.g.\,if $S$ is convex) then
$\Ext{\RRR}S$ coincides with $\R^d\setminus S$ and $\Bdry{\RRR}S=\partial S$.

\begin{figure}
\begin{tikzpicture}[scale=1.2, transform shape, > = triangle 45]
\draw[fill=red!20, color=red!20] (0,0.3)--(5,0.3)--(5,1)--
            (3,1)--(3,2)--(2,2)--(2,1.2) --(1.8,1) -- (0,1)--(0,0.3);
\draw[fill=red!20, color=red!20] (5.0,0.65) circle (0.35);

\node[color=red] at (3.5,0.65) {$S$};

\draw[fill=white, color=white] (1,0)--(1,0.3)--(1.2,0.7)--(1.4,0.3) 
    -- (1.6,0.7)--(1.8,0.3) --(2.0,0.7)--(2.2,0.3) 
    -- (2.4,0.7)--(2.6,0.3)--(2.6,0)--(1,0);

\draw[color=black!40, thin, shift={(3.71,1.71)}] (0,0) circle (0.71);
\draw[color=black!40, ->, shift={(3.71,1.71)}] (0,0) -- 
  node[pos=0.3, right]{{\smaller[2]$\RRR$}} (0.5,0.5); 

\draw[color=black!40, thin, shift={(1.29,1.71)}] (0,0) circle (0.71);
\draw[color=black!40, ->, shift={(1.29,1.71)}] (0,0) -- 
  node[pos=0.3, left]{{\smaller[2]$\RRR$}} (-0.5,0.5); 

\draw[color=black!40, thin, shift={(5.6,1.53)}] (0,0) circle (0.71);
\draw[color=black!40, ->, shift={(5.6,1.53)}] (0,0) -- 
  node[pos=0.3, right]{{\smaller[2]$\RRR$}} (0.5,0.5); 

\draw[color=blue, very thick] (2,1.71) --(2,2) --(3,2)--(3,1.71)
                   (1.29,1) --(0,1)--(0,0.3)--(1,0.3)
                   (3.71,1) --(5,1)
                   (2.6,0.3) --(5,0.3)
                   (5,0.3) arc (-90:90:0.35) ;

\fill[color=blue] (1.4,0.3) circle (0.04)
                  (1.8,0.3) circle (0.04)
                  (2.2,0.3) circle (0.04);

\end{tikzpicture}\hfill
\begin{minipage}[b]{0.45\textwidth}
\caption{Visualization of $\Bdry{\RRR}S$ (thick) as subset of the boundary
  $\partial S$ of the set $S$ (light red).}
\end{minipage}
\label{fig:bdry}
\end{figure}
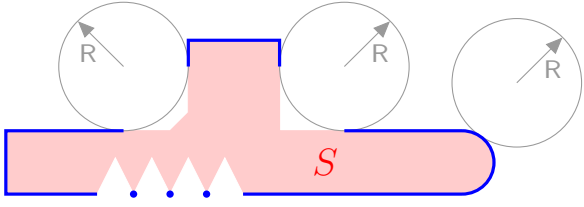

In general, $\Bdry{\RRR}S$ is a subset of the boundary of $S$, precisely made by
all points of $\partial S$ satisfying an exterior sphere condition of radius
${\RRR}$ with respect to $S$:
\begin{equation}
  \label{eq:109}
  x\in \Bdry{\RRR} S\quad\Longleftrightarrow\quad
  x\in \partial S\text{ and }\exists\,y\in \R^d:
  \ |x{-}y|={\RRR},\ B_{\RRR}(y)\cap S=\emptyset.
\end{equation}
In fact, if $x\in \Bdry{\RRR}S$ then there exist sequences $x_n, y_n$ such that
$x_n\to x$, $|x_n{-}y_n|<\RRR$ and $B_{\RRR}(y_n)\cap S=\emptyset$. Possibly
extracting a subsequence, we can assume that $y_n\to y$,
$B_{\RRR} (y)\cap S=\emptyset$, and $|x{-}y|\le \RRR$.  Since $x\in \partial S$,
it is not possible that $|x{-}y|<{\RRR}$, so that the left-to-right implication
of \eqref{eq:109} holds.  On the other hand, if $x\in \partial S$,
$|x{-}y|=\RRR$, and $B_\RRR (y)\cap S=\emptyset$, it is immediate to check
that $x\in \partial(\Ext\RRR S)$, see also Figure \ref{fig:bdry}.  

In Theorem~\ref{thm:regularity}(\ref{th:reg.label2}) we will use that for
arbitrary sets $S$ the boundary part $\Bdry\RRR S$ is countably
$(d{-}1)$-rectifiable, see \cite[Th.\,10.48(ii)]{Vill09OTON}, and hence
has $\Leb d$ measure $0$. 

The following result shows how the properties of $\LLL$ provide
regularity of the backward transform $\varphi^{\shortleftarrow\LLL}$. Of
course, an analogous statement holds for the forward transform using
\eqref{eq:156}. The important fact is that the \emph{upper bounds} on the
second derivatives of $\LLL$ generate semi-convexity of $\varphi_0$ (i.e.\
\emph{lower bounds} on $\rmD^2\varphi_0$), see Assertions \ref{Assert5} and
\ref{Assert6}. As $\rmD^2 \LLL(z)$ blows up at the boundary of $B_\RRR(0)$, it
is essential to use the fact that $\LLL(z_k)\to +\infty$ for
$ |z_k|\uparrow \RRR$. 

\begin{theorem}[Regularity of the $\LLL$-transform]
\label{thm:regularity0}
Let $\LLL: \R^d\to [0,+\infty]$ satisfy 
\begin{subequations}
  \label{eq:Cond.LLL}
  \begin{align}
\label{eq:Cond.LLL.1}
   & \LLL : \R^d\to [0,+\infty] \text{ is continuous and } \LLL(0)=0, 
\\ 
\label{eq:Cond.LLL.2}
  & \LLL\big|_{B_\RRR(0)} \in
   \rmC^2(B_\RRR(0)) \text{ and } \LLL(z)=+\infty \text{ if } |z|\geq \RRR,
\\ 
\label{eq:Cond.LLL.3}
  & \LLL \text{ is uniformly convex, i.e.\ } \exists\, \lambda_*>0 \ \forall \,
  z \in B_\RRR(0): \ \rmD^2\LLL(z) \geq \lambda_* I. 
  \end{align}
\end{subequations}
For an u.s.c.\ function $\varphi_1:\R^d\to[-\infty,+\infty]$, we consider the 
 backward $\LLL$-transform $\varphi_0=\varphi_1^{\shortleftarrow\LLL}$ and set
\begin{equation}
\label{eq:setsO_iQ_i}
\begin{aligned}&
O_0 = \{\varphi_0>-\infty\}, \quad Q_0=\{\varphi_0<+\infty\},
\\
&O_1 = \{\varphi_1<+\infty\}, \quad   Q_1=\{\varphi_1>-\infty\},\quad \text{and} \quad
\Omega_0 = O_0\cap \mathrm{int}(Q_0).
\end{aligned}
\end{equation}
Then, the following assertions hold: 
\begin{enumerate}[\upshape(1)]
\item \label{item:thm:regularity0:lsc} 
The function $\varphi_0$ is l.s.c.\ and satisfies
  \begin{align}
    \label{eq:154}
   & \inf \varphi_0\ge \inf \varphi_1 
    \quad \text{and} \quad
    \sup \varphi_0\le \sup \varphi_1,
  \\
   &\hspace*{-2em}(Q_1)^{\RRR} \subset O_0 \quad \text{and} \quad 
  \Big(\varphi_0(x_0)=-\infty\ \Leftrightarrow\ 
    B_\RRR(x_0)\subset \Ext\RRR{Q_1}\subset \{\varphi_1=-\infty\} \Big).
    \label{eq:155}
  \end{align}
  The sets $O_0,\ O_1$, and $\Omega_0$ are open. 
\item \label{item:thm:regularity0:sphere} 
The set $Q_0$ satisfies an external sphere condition of radius $\RRR$,
  namely
  \begin{equation}
    \label{eq:107}
    \R^d\setminus \Cl(Q_0)=\Ext\RRR{Q_0}\quad \text{and} \quad
    \partial Q_0=\Bdry\RRR{Q_0},
  \end{equation}
  so that the topological boundary of $Q_0$ is countably $(d{-}1)$-rectifiable.

\item \label{item:thm:regularity0:contact}
The ``contact set'' $M := M_{-\infty}\cup M_{+\infty} \cup M_\mathrm{fin}\subset 
 \R^d\times\R^d$ defined  via 
  \begin{equation}
    \label{eq:115}
    \begin{aligned}
      M_\mathrm{fin}:={}&\bigset{(x_0,x_1) }{  \varphi_i(x_i)\in \R,
      ~\varphi_1(x_1)=\LLL(x_0{-}x_1)  +\varphi_0(x_0) },
     \\[0.5em]
      M_{-\infty}:={}&\bigset{(x_0,x_1) }{\varphi_0(x_0)=-\infty,\
      |x_1{-}x_0|\ge \RRR },
     \\[0.5em]
      M_{+\infty}:={}&\bigset{(x_0,x_1) }{ \varphi_1(x_1)=+\infty,\
      |x_1{-}x_0|\ge \RRR },
    \end{aligned}
   \end{equation}
   is closed.

 \item \label{item:thm:regularity0:contactsec}
 For every $\bar x_0\in \Omega_0$, the section
   $M_{0\to1}[\bar x_0] := \bigset{x_1}{(\bar x_0,x_1)\in M_{\mathrm{fin}}}$ of
   $M_{\mathrm{fin}}$ is nonempty, compact, and included in $Q_1$. 
     Moreover, for every compact
   $K\subset \Omega_0$ there exists $\theta \in (0,\RRR) $ and
   $a',a''\in \R$ such that%
   \begin{equation}
     \label{eq:91}
     |x_1{-}\bar x_0|\le \theta \text{ and } a'\le \varphi_1(x_1)\le a'' 
     \quad\text{whenever } \bar x_0\in K 
     \text{ and } x_1\in M_{0\to1}[\bar x_0]. 
   \end{equation}

 \item \label{item:thm:regularity0:semiconvex}
 \label{Assert5} The restriction of $\varphi_0$ to the open set 
   $\Omega_0$ is locally semi-convex, and in particular locally Lipschitz and
   thus continuous.

 \item \label{item:thm:regularity0:map}
 \label{Assert6} If
   $D_0':=\mathrm{dom}(\nabla\varphi_0)\subset \Omega_0$,
   $D_0''=\mathrm{dom}(\rmD^2\varphi_0) \subset D_0'$, then $D''_0$ has
   full Lebesgue measure in $\Omega_0$. For  every $x\in D_0'$, the set
   $M_{0\to1}[x]$ contains a unique point
   $y=\bfT_{0\to1}(x)$. The induced map $\bfT_{0\to1}:D_0'\to \R^d$ is
   differentiable according to \eqref{eq:72bis} in $D_0''$ and satisfies the
   following properties:
   \begin{align}
     \label{eq:45}
     \text{(a) }& |x{-}\bfT_{0\to1}(x)|<\RRR \text{ and }
     \nabla\varphi_0(x) = (\nabla \LLL )\big(x{-}\bfT_{0\to1}(x)\big)
     \text{ for all }x\in D_0',
\\
      \label{eq:46}
      \text{(b) }& \rmD^2 \varphi_0(x) \geq
      -\rmD^2\LLL\,\big(x{-}\bfT_0(x)\big) \quad\text{for all } x\in D_0'',
\\
      \label{eq:47}
      \text{(c) }& \rmD \bfT_{0\to1}(x) \text{ is diagonalizable with nonnegative
        eigenvalues on }  D_0''. 
    \end{align}
  \end{enumerate}
\end{theorem}
\begin{proof}
We divide the proof in various steps, corresponding to each assertion.
\medskip

\noindent \underline{Assertion (\ref{item:thm:regularity0:lsc}).}
To check that $\varphi_0$ is
l.s.c.\ we assume $\varphi_0(x_0)>a$ for some $a\in [-\infty,+\infty)$, then
there exists $y \in B_\RRR(x_0)$ such that
$\varphi_1(y)-\LLL(y{-}x_0)>a$. As $\LLL $ is continuous, 
we can find $\delta\in (0,\RRR -|y{-}x_0|)$ such that
$\varphi_1(y)-\LLL(y{-}x)> a$  for every $x\in
B_\delta(x_0)$. By definition of $\varphi_0$ this estimate implies $\varphi_0(x)>a$ on $B_\delta(x_0)$, and 
lower semi-continuity is shown. 

The estimates in \eqref{eq:154} are elementary following from $\LLL(0)=0$
and $\LLL(z)\geq 0$, respectively. The relation in \eqref{eq:155}
follows from the fact that $\varphi_0(x_0)=-\infty$ implies
$\varphi_1(y)\equiv -\infty$ in $B_\RRR(x_0)$. The openness of $O_0$ and
$O_1$ follows because $\varphi_0$ is l.s.c.\ and $\varphi_1$ is u.s.c.\ This property in turn
implies that $\Omega_0= O_0 \cap \mafo{int}(Q_0)$ is open. 
\medskip

\noindent
\underline{Assertion \eqref{item:thm:regularity0:sphere}.}  Recalling $Q_0=\{\varphi_0<+\infty\}$ it is
sufficient to notice that 
\begin{subequations}
      \label{eq:86new}
  \begin{align}
         \label{eq:86new.a}
   &\bar x \in \R^d\setminus Q_0 \ \Leftrightarrow \ 
    \varphi_0(\bar x)=+\infty \ \Rightarrow \ 
    \exists\,\bar y: |\bar x{-}\bar y|\le \RRR \text{ and } \varphi_1(\bar
      y)=+\infty,
\\
    \intertext{where we used $\LLL\geq 0$
    and the upper semicontinuity of $\varphi_1$.
    However, using
  $\mathrm{dom}(\LLL)=B_\RRR(0)$ we obtain}
&      \label{eq:86new.b}
   \varphi_1(\bar y)=+\infty \quad \Rightarrow \quad 
\varphi_0(x)=+\infty \text{ for all } x\in B_{\RRR}(\bar y). 
  \end{align}
\end{subequations}
This implication means that if $\bar x\in \R^d\setminus Q_0$ then $\bar x\in
\Cl(\Ext\RRR {Q_0})$,
so that $\partial Q_0=\partial (\R^d\setminus Q_0)= \partial \Cl( \Ext\RRR
{Q_0})=\partial \Ext\RRR {Q_0}$.
  
\medskip\noindent \underline{Assertion (\ref{item:thm:regularity0:contact}).}  The closedness of $M_{\pm \infty}$
follows easily by the semi-continuities of $\varphi_i$. For $M_\text{fin}$
we consider a sequence $(x_{0,n},x_{1,n})\in M_{\mathrm{fin}}$ to $(x_0,x_1)$. 
If 
$|x_0{-}x_1|<\GGG \RRR$, then we have $\varphi_1(x_1)\ge \LLL(x_1{-}x_0) + 
\varphi_0(x_0)$ by the semi-continuities. As the opposite
inequality is always satisfied,  
we obtain the equality.  We can also exclude that
$\varphi_0(x_0)=\varphi_1(x_1)=+\infty$ (resp.\,$-\infty$), since
otherwise
$\varphi_0(x)\equiv +\infty$ in $B_\RRR(x_1)$ by
\eqref{eq:86new.b}
which contains a neighborhood of $x_0$
(resp.\,$\varphi_1(x)\equiv -\infty$ in $B_\RRR(x_0)$ by
\eqref{eq:155}, which contains a neighborhood of $x_1$),
so that
$(x_0,x_1)\in M_{\mathrm{fin}}$.
If $|x_1{-}x_0|\ge \GGG\RRR $ and $(x_0,x_1)$ does
not belong to $M_{-\infty}$ then we have
$\liminf_{n\to\infty} \varphi_0(x_{0,n})\ge \varphi_0(x_0)>-\infty$ so that
\begin{displaymath}
  \varphi_1(x_1)\ge \limsup_{n\to\infty}\varphi_1(x_{1,n})=
  \limsup_{n\to\infty}\LLL(x_{1,n}-x_{0,n})+\varphi_1(x_{0,n})=
  +\infty
\end{displaymath}
and
$(x_0,x_1)\in M_{+\infty}$. Hence, $M=M_\mathrm{fin}\cup
M_{+\infty} \cup M_{-\infty}$ is closed.\medskip 

\noindent
\underline{Assertion (\ref{item:thm:regularity0:contactsec}).}  Let us first show that \emph{$\varphi_0$ is locally bounded
  from above in the interior of $Q_0$, i.e.\,the open set
  $Q_0\setminus \partial Q_0$.}  In fact, if a sequence $x_n$ is converging to
$\bar x\in Q_0\setminus \partial Q_0$ with $ \varphi_0(x_n) \uparrow +\infty$,
by arguing as before and using 
$\varphi_0(x_n)=\sup_{y\in B_{\RRR}(x_n)}\varphi_1(y)-\LLL(y {-}x_n)$,
we find $\bar y\in \overline{B_{\RRR}(\bar x)}$ with $\varphi_1(\bar
y)=+\infty$. Now \eqref{eq:86new.b} gives 
$\varphi_0(x)=+\infty$ for all $x\in B_\RRR(\bar y)$, which  
contradicts the fact that $ \varphi_0(x) < +\infty$ in a neighborhood of $\bar
x$, because of $|\bar x{-}\bar y|\leq \RRR$. 
  
We fix now a compact subset $K$ of the open set $\Omega_0$, a point
$\bar x\in K$, and consider the section $M_{0\to1}[\bar x]$ of the
contact set $M_{\mathrm{fin}}$.  Let $\eta>0$ be sufficiently small so that
$K_\eta:=\bigset{x\in \R^d}{\dist(x,K)\le \eta }\subset \Omega_0$ and let
$\ol a :=\sup_{K_\eta} \varphi_0$, where $a<+\infty$ by the previous
claim. By l.s.c.\ of $\varphi_0$, we also have $\underline{a} := \inf_{K_\eta}
\varphi_0 >-\infty$. 

By the definition of $\varphi_0=\varphi_1^{\shortleftarrow\LLL}$, for every
$\eps\in (0,1]$  the sets  
\begin{equation}
  M^\eps(\bar x):=\Bigset{y\in B_{\RRR}(\bar x) } {
    \varphi_1(y)
    \geq
   \LLL ( y {-} \bar x)    +\varphi_0(\bar x)
   \MAT-\eps\EEE  },\label{eq:82} 
\end{equation}
are non-empty. We choose $y\in M^1(\bar x)$ and set 
$x_\vartheta:=\vartheta\bar x+ (1 {-}\vartheta) y$ with $\vartheta= 1-\eta/\RRR$, which
implies  $|x_\vartheta{-}\bar x| \leq \eta$, and hence $x_\vartheta\in
K_\eta$.
Moreover, we have $|x_\vartheta{-}y|\le \RRR-\eta$. 
Therefore, for $y \in M^1(\bar x) \subset B_\RRR(\bar x)$ we find
\begin{equation}
  \label{eq:88}
  \begin{aligned}
    \varphi_1(y)&\leq \LLL(y{-}  x_\vartheta)+\varphi_0(x_\vartheta) \leq  
    a'':=  \ol a + \wh\ell(\RRR{-}\eta)  < \infty,
\\
\varphi_1(y)&\geq \varphi_0(\bar x) + \LLL(y{-}\bar x)
-1\geq a':=
     \underline a >-\infty, 
  \end{aligned}
\end{equation}
where $\wh\ell(\varrho):= \sup_{z\in B_{\varrho}(0)} \LLL(z)$.  Combining the last two estimates we additionally find 
\begin{equation}
  \label{eq:83}
  \LLL(y{-}\bar x)\leq \varphi_1(y)-\varphi_0(\bar x) \leq a''-\underline a =:
  \wh\ell(\theta) \quad
  \text{with }\quad \theta\in (0, \RRR).
\end{equation}
Hence, all elements $y \in M^1(\bar x)$ satisfies $|\bar x{-}y|\le \theta$
and \eqref{eq:88}.

We now consider a sequence $y_\eps\in M^\eps(\bar x)\subset M^1(\bar x)$, then a
standard compactness argument and the upper semi-continuity of $\varphi_1$ show
that any limit point $\bar y$ is an element of $M_{0\to1}[\bar x]$, which is
therefore not empty.  The compactness of $M_{0\to1}[\bar x]$ and \eqref{eq:91}
again follow by \eqref{eq:83}

\medskip\noindent
\underline{Assertion (\ref{item:thm:regularity0:semiconvex}).}
Let us now fix $\bar x_0\in \Omega_0$ and $\delta>0$ such that
$K:=\overline{B_\delta(\bar x_0)}\subset \Omega_0$. The previous
assertion yields  $\theta<\RRR$ and $a',a''\in \R$ such that
$|x'{-}x|\leq \theta$ and $a'\leq \varphi_1(x')\leq a''$
whenever $x\in K$ and $x'\in M_{0\to1}[x]$. By possibly
reducing $\delta$, we can also assume that $3\delta+\theta<\RRR$.
For every $x\in K$, we now have by construction 
\begin{equation}
  \label{eq:MaximExists}
  \varphi_0(x)=\max_{x'\in \overline {B_{\delta+\theta}(\bar x_0)}}\varphi_1(x')-\LLL(x'{-}x)
\end{equation}
which is bounded and semi-convex in $K$ because it is a supremum over a
family of uniformly semi-convex functions, where we use $| x'{-}x|\leq
|x'{-}\bar x_0| + |\bar x_0{-}x|\leq 2\delta {+} \theta$ and that $-\LLL$ is
semi-convex on $\ol{B_{2\delta {+} \theta}(\bar x_0)}$ by
\eqref{eq:Cond.LLL.2}. 

\medskip\noindent \underline{Assertion (\ref{item:thm:regularity0:map}).} 
This assertion follows in the standard way 
by using the extremality conditions in the contact set, see
e.g.\,\cite[Thm.\,6.2.4 and 6.2.7]{AmGiSa08GFMS}. We give the main
argument to show how the assumptions in \eqref{eq:Cond.LLL} enter.  By
Alexandrov's theorem and Assertion (\ref{item:thm:regularity0:semiconvex}) the set $D''_0$ has full
Lebesgue
measure. To obtain the optimality conditions, we fix $x_0\in Q_0\cap D''_0$ and know from
\eqref{eq:MaximExists} that there exists $\bar x_1$ such that 
$\varphi_0(x_0)=\varphi_1(\bar x_1)- \LLL(\bar x_1{-}x_0)$. However, for all
$x\in B_\delta(x_0)$ we have
$\varphi_0(x)+\LLL(\bar x_1{-}x)\geq \varphi_1(\bar x_1) $ with equality for
$x=x_0$. Thus, we obtain the optimality conditions 
\[
  \nabla\varphi_0(x_0)- \nabla \LLL(\bar x_1{-}x_0)=0 \text{ in } \R^d \quad \text{and}
  \quad \rmD^2 \varphi_0(x_0) + \rmD^2\LLL(\bar x_1{-}x_0)\geq 0 \text{ in }
  \R^{d\ti d}_\text{sym}.
\]
This result gives the conditions (a) to (c), if we observe that $\bar x_1$ is
unique. But this property follows from the first optimality condition by using
\eqref{eq:Cond.LLL.3} which allows us to write 
\[
\bar x_1 = \bfT_{0\to1}(x_0):= x_0+ \big( \nabla\LLL\big)^{-1}(\nabla
\varphi_0(x_0)),
\]
i.e.\ $\bar x_1$ is uniquely determined by $x_0$. Moreover, $\rmD \bfT_{0\to1}(x_0)$
exists and satisfies $\rmD^2 \varphi_0(x_0)= (\rmD^2
\LLL)(\bfT_{0\to1}(x_0){-}x_0) \big(\rmD \bfT_{0\to1}(x_0){-}I\big)$, which implies
the diagonalization result. 
\end{proof}

The previous result can now be applied to the solution of the LET problem in
Theorem~\ref{thm:LET} using $\LLL=\Ell_1$;
thus in this case $\RRR=\pi/2$.
Using the notations for
$\mafo{supp}(\mu_i)=S_i=S'_i+S''_i$ and $\mu_i=\mu'_i+ \mu''_i$ from Theorem~\ref{thm:effective} we can compare these to the sets $O_i$, $Q_i$,
$D'_i$, and $D''_i$ defined for an optimal pair $(\varphi_0,\varphi_1)$ as in
Theorem~\ref{thm:regularity0}. So far we constructed optimal pairs
$(\varphi_0,\varphi_1)$ satisfying 
\begin{equation}
  \label{eq:OptimPairs}
\varphi_0\geq \varphi_1^{\shortleftarrow\Ell_1} \text{on }\R^d, \quad 
 \varphi_0^{\Ell_1\to}\geq \varphi_1\text{ on } \R^d, \quad
\varphi_0= \varphi_1^{\shortleftarrow\Ell_1} \ \mu_0\text{-a.e.},\quad 
 \varphi_0^{\Ell_1\to}\MAT= \EEE \varphi_1 \  \mu_1\text{-a.e.}
\end{equation}
However, following 
\cite[Ch.\,5]{Vill09OTON}, we will show that it is possible to restrict to ``tight
optimal pairs'' satisfying $\varphi_0=\varphi_1^{\shortleftarrow \Ell_1}$ and
$\varphi_1= \varphi_0^{\Ell_1\to}$, which implies that $\varphi_0$ is l.s.c.\
and $\varphi_1$ is u.s.c.\ This possibility leads to the following refinement of 
the results in \cite[Thm.\,6.6(iii)]{LiMiSa18OETP}.

\begin{theorem}[Regularity of optimal $\HK$ potentials]
\label{thm:regularity}
Let $\mu_0,\mu_1$ be nontrivial measures in $\MMM(\OOmega)$ with
decompositions given by \eqref{eq:33}--\eqref{eq:43}. 
\begin{enumerate}[{\upshape(1)}]

\item \label{th:reg.label1} There exists an optimal pair of potentials
  $\varphi_0,\varphi_1:\R^d\to[-\infty,+\infty]$ with $\varphi_0$ being
  l.s.c.  and $\varphi_1$ u.s.c., solving the dual problem of Theorem
  \ref{thm:Borel-duality} and 
  \begin{gather}
    \label{eq:77pre}
   \varphi_0=\varphi_1^{\shortleftarrow \Ell_1}\quad \text{ and } \quad 
    \varphi_1= \varphi_0^{\Ell_1\to} \quad\text{on }\R^d,
    \\
    \label{eq:104}
     S_i\subset Q_i, \quad
    S_0'\subset S^\pit_1 \subset O_0,\quad
    S_1'\subset S^\pit_0  \subset O_1,
    \\
    \varphi_0=-\infty     \text{ on }   S_0'',\quad\text{and}\quad
    \varphi_1=+\infty\text{ on }  S_1'',
    \label{eq:105}
  \end{gather}
  where the sets $O_i$ and $Q_i$ are as in \eqref{eq:setsO_iQ_i}.

\item \label{th:reg.label2}
  If $\bfeta$ is
  an optimal solution of the LET problem \eqref{eq:HK.ET}, the functions
  $\sigma_0:=\mathrm e^{2\varphi_0}$ and $\sigma_1:=\mathrm e^{-2\varphi_1}$
  provide lower semi-continuous representatives of the densities of the
  marginals $ \eta_i = \pi^i_\sharp\bfeta$ with respect to $\mu_i$, i.e.,
  $\eta_i=\sigma_i\mu_i$, and $\bfeta$ is concentrated on the contact
  set $M_{\rm fin}$ so that 
  $\supp(\bfeta)\subset M$ (see Theorem~\ref{thm:regularity0}).
  The marginals $\eta_i$ are concentrated on the open sets $O_i$.  
  
  Conversely,
  if $\wt\bfeta$ satisfies $\supp(\wt\bfeta)\subset M$ and $\wt\eta_i=\sigma_i\mu_i$, then $\wt\bfeta$
  is an optimal solution of the LET problem \eqref{eq:HK.ET}.

\item \label{th:reg.label2bis}
  If $\mu_0$ (resp.~$\mu_0'$) does not charge $(d{-}1)$-rectifiable
  sets, e.g.\ in the case that $\mu_0\ll\Leb d$ or if $\mu_0(\Bdry\pit { S_0})=0$
  (resp.~$\mu_0'(\Bdry\pit { S_0})=0$), then for
  every optimal pair $(\varphi_0,\varphi_1)$ with
  $\varphi_0=\varphi_1^{\shortleftarrow\Ell_1}$ and $\varphi_1$ u.s.c., the
  measure $\mu_0$ is concentrated on the open set $\Int(Q_0)$
  (resp.~$\mu_0'$ is concentrated on the open set $\Omega_0$).
\item \label{th:reg.label3} If $\mu_0'$ is concentrated on
  $D_0'=\mathrm{dom}(\nabla\varphi_0)$ (in particular if $\mu_0'\ll \Leb d$)
  then the optimal transport plan $\bfeta$ solving the LET formulation is
  unique, it is concentrated on $D_0'\times S^\pit_0$, and it is induced
  by the graph of $\bfT_{0\to1}$,
  i.e.\,$\bfeta=(\mathrm{Id},\bfT_{0\to1})_\sharp \eta_0$ with
  $\bfT_{0\to1}$ from Theorem~\ref{thm:regularity0}(\ref{Assert6}). 

\item \label{th:reg.label4} If $\mu_0',\mu_1'\ll\Leb d$ then $\mu_0'$ is
  concentrated on $D_0''\cap \bfT_{0\to1}^{-1}(D_1'')$, where
  $D_i''=\mathrm{dom}(\rmD^2\varphi_i)$, and $\bfT_{0\to1}$ is
  $\mu_0'$-essentially injective with $\det \rmD\bfT_{0\to1}>0$ \ $\mu_0$-a.e.\,in
  $D_0''$.
\end{enumerate}
\end{theorem}
\begin{proof}
\noindent\underline{Assertion (\ref{th:reg.label1}).}
Let $(\phi_0,\phi_1)$ be an optimal Borel pair according to
Theorem~\ref{thm:optimality.cond}(\ref{thm:optim.cond.lab3}), see
\eqref{eq:19}, satisfying
\begin{equation}
  \label{eq:111}
  \phi_i\in \R\ \text{ $\mu_i$-a.e.\,in }S_i',\quad
  \phi_0=-\infty\ \text{ $\mu_0$-a.e.\,in }S_0'',\ 
  \phi_1=+\infty\ \text{ $\mu_1$-a.e.\,in }S_1''.
\end{equation}
With this pair, we set $\varphi_0:=\phi_1^{\shortleftarrow\Ell_1}$, and recalling
\eqref{eq:Def.Ell.Trafo} we easily obtain
\begin{equation}
  \label{eq:84}
  \varphi_0\le \phi_0\quad\text{in }\OOmega,\quad
  \phi_1(x_1)   \le \Ell_1(x_1{-}x_0)+\varphi_0(x_0)
  \quad \text{if }x_0,x_1\in \R^d,\ |x_1{-}x_0|<\pi/2.
\end{equation}
Looking at the dual problem \eqref{eq:HK.dual.a} with the more general
admissible set of Borel pairs as described in  \eqref{eq:condphi}, we see that
$ (\varphi_0,\phi_1)$ is still optimal. 

Repeating the argument, we can set $\varphi_1= \varphi_0^{\Ell_1\to}$ to find a
new optimal pair satisfying $\varphi_1 \geq \phi_1$. However, exploiting
\eqref{eq:forw.backw.Ell} we see that the tightness relation \eqref{eq:77pre}
holds for the optimal pair $(\varphi_0,\varphi_1)$. This fact implies that
$\varphi_0$ is l.s.c.\ and $\varphi_1$ is u.s.c. 

By the construction of $\phi_i$ in Theorem 
\ref{thm:optimality.cond}(\ref{thm:optim.cond.lab3}) we have 
\[
\{\phi_i \in \R\} = S'_i, \quad \{ \phi_0=-\infty\} = S''_0, \quad 
\text{and }\  \{ \phi_1=+\infty\} = S''_1.
\]
Together with $\phi_0 \geq \varphi_0 $ and $\phi_1 \leq \varphi_1$ we find 
\begin{align*}
&S''_0  = \{ \phi_0=-\infty\} \subset \{ \varphi_0=-\infty\} \ \text{ and }
\  S_0 = \{ \phi_0 <+\infty\} \subset  \{ \varphi_0 <+\infty\}=Q_0,
\\
&S''_1  = \{ \phi_1=+\infty\} \subset \{ \varphi_1=+\infty\} \ \text{ and }
\  S_1 = \{ \phi_1 >-\infty\} \subset  \{ \varphi_1>-\infty\} = Q_1. 
\end{align*}
Clearly, $S'_0 =S_0\cap S^\pit_1\subset S^\pit_1$. Moreover, for $x_0 \in
S^\pit_1$ we find $y_1\in S_1$ with
$|y_1{-}x_0| <\pit$, i.e.\ $\Ell_1(y_1{-}x_0)<\infty$. With this we have
$\varphi_0(x_0)= \varphi_1^{\shortleftarrow\Ell_1}(x_0) \geq \varphi_1(y_1) -
\Ell_1(y_1{-}x_0)>- \infty$ and conclude $x_0 \in O_0$. Thus, $S'_0\subset
S^\pit_1 \subset O_0$ is shown and $S'_1\subset
S^\pit_0 \subset O_1$ follows similarly. Hence, \eqref{eq:104} and
\eqref{eq:105} are established.

\medskip\noindent \underline{Assertion (\ref{th:reg.label2}).} The claim follows immediately from Theorem~\ref{thm:optimality.cond}.

\medskip\noindent \underline{Assertion (\ref{th:reg.label2bis}).}
We just consider the case of $\mu_0$, since the argument for
$\mu_0'$ is completely analogous and eventually uses the fact that
$\Omega_0=O_0\cap \Int(Q_0)$ and 
$\mu_0'$ is also concentrated on $O_0$ by \eqref{eq:104}.

By Theorem~\ref{thm:regularity0} (cf.\ \eqref{eq:107}) we know that 
$\partial Q_0=\Bdry\pit{Q_0}$. Since $\partial Q_0$ is
$(d{-}1)$-rectifiable and $\mu_0$ does not charge $(d{-}1)$-rectifiable
sets, we conclude $\mu_0(\partial Q_0)=0$.  

If $\mu_0(\Bdry\pit{S_0})=0$, we also obtain  $\mu_0(\partial Q_0)=0$
via the following arguments: By \eqref{eq:104} we have $S_0\subset Q_0$, which
implies that a point $x \in \partial S_0\cap \Bdry\pit{Q_0}$ also lies
$\Bdry\pit{S_0}$. Using $\partial Q_0=\Bdry\pit{Q_0}$ we obtain $\partial S_0
\cap \partial Q_0 \subset \Bdry\pit{S_0}$ and find  
\[
\mu_0(\partial Q_0) \overset{\text{(i)}}= \mu_0(\pl Q_0 \cap S_0) 
\overset{\text{(ii)}}= \mu_0(\pl Q_0 \cap \pl S_0) \leq 
\mu_0\big(\Bdry\pit{S_0}\big)= 0,
\]
where we used $S_0= \mafo{sppt}(\mu_0)$ in $\overset{\text{(i)}}=$ and $S_0
\subset Q_0$ in $\overset{\text{(ii)}}=$. Thus, we have shown that $\mu_0$ is
concentrated on $\Int(Q_0)$.

\medskip\noindent \underline{Assertion (\ref{th:reg.label3}).}
If $\mu_0'\ll\Leb d$ then $\mu_0'$ is concentrated on $\Omega_0$ by
Claim \ref{th:reg.label2bis} and
$\mu_0(\Omega_0\setminus D_0')=0$ by \ref{thm:regularity0}(\ref{Assert6}). 
By the previous claim \ref{th:reg.label2}, we know that the first marginal $\eta_0$ of
$\bfeta$ is given by $\rme^{2\varphi_0}\mu_0=
\rme^{2\varphi_0}\mu_0'=\rme^{2\varphi_0}\mu_0'\res D_0'$
(in particular $\eta_0(\R^d\setminus D_0')=0$) so that
$\bfeta$ is concentrated on 
$M_{\rm fin}\cap (D_0'\times \R^d)$ which is the graph of the map $\bfT_{0\to1}$
given by Theorem~\ref{thm:regularity0}(\ref{Assert6}). 
\medskip

\noindent \underline{Assertion (\ref{th:reg.label4}).}
Let us first recall that for $i=0,1$ the marginal $\eta_i$ of $\bfeta$ and the
measure $\mu_i'$ are
mutually absolutely continuous.
Since $\mu_i'\ll \Leb d$ we know by Theorem~\ref{thm:regularity0}(\ref{Assert6}) and the third claim 
that $\mu_i'(\R^d\setminus D_i'')=\mu_i'(\Omega_i\setminus D_i'')=0$,
so that $\eta_i(\R^d\setminus D_i'')=0$ and
$\eta_0(\bfT^{-1}_{0\to1}(\R^d\setminus D_1''))=
\eta_1(\R^d\setminus D_1'')=0$; we deduce that $\eta_0$ and $\mu_0'$ are
concentrated on $D_0''\cap \bfT^{-1}_{0\to1}(D_1'').$

We can apply Theorem \ref{thm:regularity0}(\ref{Assert6}),
inverting the order of the pair $(\varphi_0,\varphi_1)$ and obtaining
that for every $x_1\in D_1'$ there is a unique element $x_0\in \R^d$
in the section $M_{1\to 0}(x_1)$, i.e.~such that
$(x_0,x_1)\in M_{\rm fin}$. This result precisely shows that the restriction
of $\bfT_{0\to1}$ to $D_0'\cap \bfT^{-1}_{0\to1}(D_1')\supset
D_0''\cap \bfT^{-1}_{0\to1}(D_1'')$ is
injective. Since $(\bfT_{0\to1})_\sharp \eta_0=\eta_1\ll\Leb d$, we
can eventually apply \cite[Lemma 5.5.3]{AmGiSa08GFMS} which shows that
$\det \rmD \bfT_{0\to1}>0$ $\mu_0$-a.e.~in $D_0''$.
\end{proof}

It is important to realize that the tightness condition \eqref{eq:77pre} is
strictly stronger than the optimality conditions \eqref{eq:OptimPairs}.
However, 
even for tight optimal pairs there is some freedom outside the supports of the measures $\mu_0$ and $\mu_1$, 
as is seen in the following simple case.

\begin{example}[Tight optimal pairs for two Diracs]
\label{ex:TightDirac}
This example lies in-between Examples \ref{ex:EllTransform} and
\ref{ex:ContactDirac}. For two points $z_0,z_1\in \R^d$ with
$\varrho=|z_1{-}z_0|= \pi/3$, such that $\tcos{\pit}\varrho=1/2$. We consider
two measures $\mu_i=\delta_{z_i}$. With $s_i=S'_i =\{z_i\}$ we easily find the
two optimal potential $(\phi_0,\phi_1)$ according to Theorem~\ref{thm:optimality.cond}, see
\eqref{eq:19}:
\[
\phi_0(x_0) = \begin{cases}-\frac{\log 2}2 &\text{for } x_0=z_0, \\
 +\infty& \text{otherwise}, \end{cases} 
\quad \text{and} \quad 
\phi_1(x_1) = \begin{cases}\frac{\log 2}2 &\text{for } x_1=z_1, \\
 -\infty& \text{otherwise}, \end{cases} 
\] 
In particular, we have $\phi_1(z_1)-\phi_0(z_0)=\log 2 =
\Ell_1(z_1{-}z_0)=\frac12\ell(\varrho)$.  

Proceeding as in Step 1 of the above proof with $\varphi_0=
\phi_1^{\shortleftarrow \Ell_1}$ and taking into account the
calculations of Example \ref{ex:EllTransform}, we obtain a first tight optimal pair 
\[
(\varphi^{(1)}_0,\varphi^{(1)}_1) \quad \text{with } \varphi^{(1)}_0(x_0)= \begin{cases}
  \frac{\log 2}2 - \Ell_1(z_1{-}x_0) & \text{for } x_0 \in B_\pit(z_1), \\ 
-\infty& \text{otherwise}, 
\end{cases} \quad \text{and } \ \varphi^{(1)}_1= \phi_1. 
\]
Interchanging the roles of $\phi_0$ and $\phi_1$ we arrive at a second  tight
optimal pair  
\[
(\varphi^{(2)}_0,\varphi^{(2)}_1) \quad \text{with } \varphi^{(2)}_0=\phi_0
 \quad \text{and } \ \varphi^{(2)}_1(x_1)= \begin{cases}
  -\frac{\log 2}2 + \Ell_1(x_1{-}z_0) & \text{for } x_1 \in B_\pit(z_0), \\ 
 \infty& \text{otherwise}. \end{cases} 
\]

A third case is obtained by choosing $z_1^*\neq z_1$ and considering an optimal
pair $(\phi_0,\wt\phi_1)$ with $\phi_0$ from above and   
\[
\wt\phi_1(x_1)= \begin{cases} \frac{\log 2}2 & \text{for }
  x_1=z_1 , \\  a_1 &\text{for } x_1=z_1^*, \\ -\infty &\text{otherwise}.
\end{cases} \quad \text{where } a_1\leq - \frac{\log 2}2 + \Ell_1(z_1^*{-}z_0).
\]
We obtain $\varphi^{(3)}_0: x_0\mapsto \max\{\frac{\log 2}2 -
\Ell_1(z_1{-}x_0),\, a_1 - \Ell_1(z_1^*{-}x_0) \} $ and the tight optimal pair
$\big(\varphi^{(3)}_0 ,( \varphi^{(3)}_0)^{\Ell_1\to}\big)$. 

With the notation of Theorem~\ref{thm:regularity0} we have
$O_0^{(3)}=\{\varphi_0^{(3)} >- \infty\} = B_\pit(z_1)\cup B_\pit(z_1^*) = \wt Q_1^\pit$,
since $\wt Q_1=\{\wt\phi_1>-\infty\} = \{z_1,z_1^*\}$, i.e.\ \eqref{eq:155}
holds. Because of $Q_0^{(3)} = \{\varphi_1^{(3)} <+\infty\}=\R^2$, also \eqref{eq:107}
is true. 
\end{example}

The following corollary shows that in the case of an absolutely
continuous reduced pair $(\mu_0,\mu_1)$ the density of $\mu_1$ can be
written in terms of the optimal pair $(\sigma_0,\sigma_1)$, the transport map
$\bfT$, and the density of $\mu_0$, and vice versa.

\begin{corollary}[Monge solutions]
\label{cor:Monge}
Let $\mu_0,\mu_1\in \MMM(\OOmega)^2$  with  $\mu_1''=0$, and let
$(\varphi_0,\varphi_1)$ be a tight optimal pair of potentials
according to Theorem~\ref{thm:regularity}.
If $\mu_0'$  is
concentrated on $D_0'=\mathrm{dom}(\nabla\varphi_0)$
(in particular  if $\mu_0'\ll\Leb d$), then there exists
a ``unique'' 
(up to  $\mu_0$-negligible sets) optimal \TRAGRO\ pair $(\bfT,q)$ attaining the minimum for the Monge Problem \ref{prob:Monge}, namely 
\begin{equation}
  \label{eq:27bis}
  (\bfT,q)_\star \mu_0=\mu_1 \quad \text{and } \quad 
  \Mcost q\bfT{\mu_0}=\HK^2(\mu_0,\mu_1).
\end{equation}
If $\sigma_i,\varphi_i,D_i',D_i'',\bfeta$, $\bfT_{0\to1},\bfT_{1\to0}$
are given as in Theorem~\ref{thm:regularity0} and \ref{thm:regularity}, the
pair $(\bfT,q)$ can be obtained in the following way:
\begin{enumerate}[{\upshape(1)}]

\item The restriction of $\bfT$ to $D_0'$ coincides with the map
  $\bfT_{0\to1}$ (and the plan $\bfeta$) as in Theorem
  \ref{thm:regularity0}, whereas $\bfT(x):=x$ for every
  $x\in \R^d\setminus D_0'$ (in particular in $S_0''$).

\item $q(x)\equiv0 $ for $ x\in \R^d\setminus D_0'$ (in particular in
  $S_0''$) and 
  \begin{equation}
    q^2(x) = \frac{\sigma_0(x)}{\sigma_1(\bfT_{0\to1}(x))}=
    \sigma^2_0(x)+\frac 14|\nabla\sigma_0(x)|^2
    \text{ for  }x\in D_0'.\label{eq:114}
  \end{equation}
\end{enumerate}
Moreover, $\bfT$ satisfies
\begin{equation}
  \label{eq:52}
  |\bfT(x){-}x|<\pi/2\quad\text{and}\quad
  \sigma_0(x)\sigma_1(\bfT(x))=\cos(|x{-}\bfT(x)|)^2
  \quad\text{in } D_0'.
\end{equation}
If $\mu_0\ll\Leb d$, then $\mu_1\ll \Leb d$ if and only if
$\det \rmD \bfT(x)>0$ for $\mu_0$-a.e.\,$x\in D_0''$.  In this case, setting
$\mu_i=c_i\Leb d\ll \Leb d$ we have
\begin{equation}
  \label{eq:49}
  c_1=\Big(c_0\frac{q^2}{\det \rmD \bfT}\Big)\circ
  \bfT^{-1}\quad\text{$\Leb d$-a.e.\,in
    $\bfT(D_0'')\subset \Omega_1$}.
  \end{equation}
\end{corollary}

To obtain the second identity in \eqref{eq:114}, we exploit the first-order
optimality \eqref{eq:45} and $\sigma_0=\ee^{2\varphi_0}$ giving
$\frac1{2\sigma_0} \nabla \sigma_0= \nabla \varphi_0(x) = \bftan(x{-} \bfT(x))$
by \eqref{eq:45}.  Thus, using the optimality condition \eqref{eq:52} (coming
from \eqref{eq:34}) we find 
\begin{equation}
  \label{eq:148}
  q^2(x)=\frac{\sigma_0^2(x)}{\cos^2(|x{-}\bfT(x)|)} =
  \sigma_0^2(x)(1 {+}\tan^2(|x{-}\bfT(x)|))=
  \sigma_0^2(x)+\frac 14|\nabla\sigma_0(x)|^2.
\end{equation}
We can also rephrase the above results in terms of the optimal Kantorovich
potential $\psxi_0$ in \eqref{eq:HK.dual.b}. This potential, which satisfies
the relations $\psxi_0 = \frac12(\sigma_0{-}1) =  \Gtrafo_1(\varphi_0) = 
\frac1{2}(\ee^{2\varphi_0}{-}1)$, will be the best choice for characterizing the
densities of the Hellinger--Kantorovich geodesic curves. Indeed, the transport
map $\bfT$ on $D_0''$ takes the form
\begin{equation}
  \label{eq:transport.map}
\begin{aligned}
  \bfT(x) &= x  +
  \bfarctan\Big(\frac{\nabla\psxi_0(x)}{1{+}2\psxi_0(x)}\Big)=
  \GGG x  +
  \bfarctan\Big({\nabla\varphi_0(x)}\Big),
  \\
  q^2(x)&=(1{+} 2 \xi_0(x))^2+  |\nabla\xi_0(x)|^2.    
\end{aligned}    
\end{equation}
If $\mu_0,\mu_1$ have full support $S_0=S_1=\R^d$, then Theorem~\ref{thm:regularity} immediately yields
$\Omega_i = O_i\cap \mafo{int}(Q_i)=\R^d$, so that $\varphi_0$ and $\varphi_1$
take values in $\R$, are locally Lipschitz, and locally semi-convex and
semi-concave, respectively.  Another important case where the properties of
$\varphi_0,\varphi_1$ can be considerably refined is when $\mu_0,\mu_1$ are
\emph{strongly reduced} (cf.\ Definition \ref{def:reduced}) and have compact support.

\begin{theorem}[Improved regularity in case of strongly reduced pairs]
\label{thm:extra-regularity}
Let us assume that the supports $S_0,S_1$ of $\mu_0,\mu_1$ are compact and
satisfy $S_i\subset S^\pit_{1-i}$, so that $\mu_0,\mu_1$ is a strongly
reduced pair (cf.\ Definition~\ref{def:reduced}).  Then it is possible to find
a pair of optimal potentials $\varphi_0,\varphi_1$ as in Theorem~\ref{thm:regularity} satisfying the additional properties:
  \begin{enumerate}[{\upshape(1)}]
  \item $\varphi_i$ are uniformly bounded
    (in particular $\Omega_i=\R^d$ and $M=M_{\mathrm{fin}}$):
    there exist
    constants $\phi_{\rm min}<\phi_{\rm max}\in \R$ such that
    \begin{equation}
      \label{eq:103}
      \phi_{\rm min}\le \varphi_i\le \phi_{\rm max}\quad\text{in }\R^d.
    \end{equation}

  \item If $\theta\in [0,\pi/2[$ satisfies $\cos^2(\theta)=\ee^{2(\phi_{\rm min}-\phi_{\rm max})}$ then for every
    $x_0,x_1\in \R^d$
    \begin{equation}
      \label{eq:113}
      (x_0,x_1)\in M\quad\Rightarrow\quad
      |x_1{-}x_0|\le \theta.
    \end{equation}

  \item 
    $\varphi_i$ are Lipschitz,    
    $\varphi_0$ is semi-convex, $\varphi_1$ is semi-concave.    
  \end{enumerate}
\end{theorem}
\begin{proof}
\underline{Assertion (1).} Let $\varphi_0',\varphi_1'$ be an optimal pair as in
Theorem~\ref{thm:regularity}.  Since $\varphi_1'$ is u.s.c~and
$\varphi_1'<+\infty$ on $S^\pit_0$, we have
$\phi_{\rm max}:=\max_{S_1}\varphi_1'<+\infty$.  We can then define
$\zeta_1:= \min\{\varphi_1' , \phi_{\rm max}\}$ observing that
$\zeta_1\le \phi_{\rm max}$ in $\R^d$ and $(\varphi_0',\zeta_1)$ is still
optimal since $\zeta_1=\varphi_1'$ on $S_1$.

Arguing as in the proof of Theorem~\ref{thm:regularity}, we define
$\zeta_0:=(\zeta_1)^{\shortleftarrow\Ell_1}$, observing that
$\zeta_0\le \phi_{\rm max}$ as well.  On the other hand, $\zeta_0$ is
l.s.c.\,and $\zeta_0>-\infty$ on $S^\pit_1 \supset S_0$, so that
$\phi_{\rm min}:=\min_{S_0}\zeta_0'>-\infty$.  Setting
$\zeta_0':= \max\{\zeta_0 , \phi_{\rm min}\} $ we obtain a new optimal
pair $(\zeta_0',\zeta_1)$ with $\phi_{\rm min}\le \zeta_0\le \phi_{\rm max}$.
Hence, with $\zeta_1':=(\zeta_0')^{\Ell_1\to}$ we get the
desired optimal pair $(\zeta_0',\zeta_1')$ satisfying
$\phi_{\rm min}\le \zeta_i'\le \phi_{\rm max}$ as well. 
 \smallskip

\noindent
  \underline{Assertion (2).} This assertion is now an easy consequence of the definition of
  contact set \eqref{eq:115} and the fact that
  $\varphi_1(x_1)-\varphi_0(x_0)\le \phi_{\rm max}-\phi_{\rm min}$.\\[0.2em]
  \underline{Assertion (3).} The last assertion follows as Theorem~\ref{thm:regularity0}(\ref{item:thm:regularity0:semiconvex}).  
\end{proof}

\section[Dynamic duality and regularity properties of the Hamilton--Jacobi equation]
{Dynamic duality and regularity properties of \\  
the Hamilton--Jacobi equation}
\label{sec:HJ}

\let\one\tau

In the previous section, the regularity properties of the optimal
$\HK$ pairs $(\varphi_0,\varphi_1)$ were studied, which can be understood via
the static formulations of $\HK$ as only the measures $\mu_0$ and $\mu_1$ are
involved. Now, we consider the dual potentials $\xi_t(x)=\xi(t,x)$ along
geodesics $(\mu_t)_{t\in [0,1]}$. At this stage, the present Section~\ref{sec:HJ} is completely independent of the previous Section~\ref{se:2nd.Optim}. Only in the upcoming Section~\ref{se:GeodCurves}, we will
combine the two results to derive the finer regularity properties of the
geodesics $\mu_t$. 

In \cite[Sect.\ 8.4]{LiMiSa18OETP}, it is shown that the optimal dual potentials
$\xi$ in the dynamic formulation in \eqref{eq:15} (but now for $\alpha=1$
and $\beta=4$) are subsolutions to a suitable Hamilton--Jacobi equation, namely
\begin{equation}
  \label{eq:16}
  \begin{aligned}
    \frac 1{2\tau} \HK(\mu_0,\mu_\tau)^2= \sup\bigg\{\int_\OOmega
    &\xi(\tau,\cdot)\dd\mu_\tau-\int_\OOmega\xi(0,\cdot)\dd\mu_0\,\Big|\, \xi\in
    \rmC^\infty_\rmc([0,\tau]\times \R^d),
    \\
    & \frac\partial{\partial t}\xi+
    \frac12 |\nabla \xi|^2+2 \xi^2\le0 \quad\text{in
    }[0,\tau]\times\R^d\bigg\}.
  \end{aligned}
\end{equation}
Theorem 8.11 in \cite{{LiMiSa18OETP}} shows that the maximal subsolutions of
the generalized Hamilton--Jacobi equation \eqref{eq:genHJeqn} for
$t\in (0,\tau)$ are given by the following \emph{generalized Hopf--Lax formula}
\begin{equation}
\xi_t(x)=\xi(t,x) = \big(\HopfLax t{ \xi_0}\big)(x)=
  \frac 1{t} \HopfLax 1{\big(t\xi_0(\cdot)\big)}(x)=\inf_{y\in\OOmega}
\frac{1} 
{2t}
\Big(1-\frac{\tcosq{\pit}{|x{-}y|}}{1+
2t\xi_0(y)}\Big),\label{eq:23}
\end{equation}
where $\xi_0\in\rmC^1(\OOmega)$ is fixed and such that
$\inf_{\OOmega}\xi_0(\cdot)>-\frac 1{2\tau} $, compare with \eqref{eq:23intro}.

In the spirit of the previous section, it is possible to derive  some semi-concavity properties of
$\xi_t$ from this
formula. However, these are not enough as we need
more precise
second order differentiability. To obtain
the latter, we use the fact that a geodesic curve is not oriented, meaning that
$t \mapsto \mu_{1-t}$ is still a geodesic, or in other words that $t\mapsto
\xi_{1-t}$ has to also solve a Hamilton--Jacobi equation. Thus, our strategy
will be the following: For an optimal pair $(\xi_0,\bar \xi_1)$ in \eqref{eq:HK.dual.b}, we construct a
forward solution $\xi_t$ starting from $\xi_0$ and backward solutions starting
from $\bar\xi_1$ via 
\begin{equation}
\label{eq:ForwardBackwardSolution}
\xi_t = \HopfLax t{ \xi_0} \text{ for }t\in (0,1]\quad \text{and} \quad
\bar\xi_t = \RopfLax t\bar\xi_1:=- \HopfLax t({-}\bar\xi_1) \text{ for } t\in [0,1).
\end{equation}
In Section \ref{se:GeodCurves}, optimality will be used to guarantee
that $\xi_t$ and $\bar \xi_t$ are essentially the same so that semi-concavity
of $\xi_t$ and semi-convexity of $\bar\xi_t$ provide the desired smoothness.

\subsection{Exploiting the generalized Hopf--Lax formula for
  regularity} 
\label{su:HopfLax.reg}
 
In this section, we study in detail the regularity properties of the function
$\xi_t$ arising in \eqref{eq:23}.  Assuming that
$\inf_{x\in \OOmega}\xi_0(x)\ge -
\frac 1{2\tau}$ we see that $\mathscr P_t\xi_0$
is well-defined for $t\in (0,\tau)$  and can be
equivalently characterized by
\begin{align}
\label{eq:23bis}
  \big(\HopfLax t{ \xi_0}\big)(x)
  &= \inf_{y \in B_\pit(x)}
\frac{1} {2t} \Big(1- 
 \frac{\tcosq{\pit}{|x{-}y|}}{1+ 2t\xi_0(y)}\Big).
\end{align}
We can extend \eqref{eq:23bis} at $t=\tau$ if we define the quotients
$a/0:=+\infty$, $a/(+\infty):=0$ for every $a>0$.
Moreover, since $t\mapsto \HopfLax t{\xi_0}(x)$ is decreasing, we
easily get
\begin{equation}
  \label{eq:63}
  \xi_t(x)=\big(\HopfLax t{ \xi_0}\big)(x)=\lim_{s\uparrow t }
  \big(\HopfLax s{\xi_0}\big)(x)\quad\text{for every
  }x\in \OOmega,\ t\in (0,\tau]
\end{equation}
so that many properties concerning the limiting case $t=\tau$ can be
easily derived by continuity as $t\uparrow\tau$.

If $\xi_0$ is l.s.c.\,and
$\big(\HopfLax t{ \xi_0}\big)(x)<\frac 1{2t}$,
the infimum in \eqref{eq:23bis} it attained at a compact set denoted
by
\begin{equation}
  \label{eq:119}
  \ArgM_t\xi_0(x):=\operatorname{argmin}_y
  \frac{1}{2t}
  \Big(1-\frac{\tcosq{\pit}{|x{-}y|}}{1+
    2t\xi_0(y)}\Big)\subset {B_{\pit}(x)} .
\end{equation}
Notice that $\big(\HopfLax t{ \xi_0}\big)(x)=\frac 1{2t}$ only if $\xi_0$ is
identically $+\infty$ in $B_{\pit}(x)$ and in this case any element
of $\overline {B_{\pit}(x)}$ is a minimizer. For later usage we also
define $\ArgM_0\xi(x)=\{x\}$. 

We also observe that if $\xi_0(x)=a$ is constant then $\HopfLax t\xi_0$
is constant in $x$, namely 
\begin{equation}
  \label{eq:147}
  \HopfLax t\xi_0(x)=  \HopfLax t{a }(x)=P_a(t):=\frac {a}{1+2at}, \quad
  \text{with }P_{\infty}(t):=\frac 1{2t}.
\end{equation}
A crucial property of \eqref{eq:23} is the link with
the classical Hopf--Lax formula on the cone $\mfC$ for a function
$\zeta:\mfC\to \R$ satisfying 
$\zeta([x,r])\ge -\frac 1{2\tau} r^2$. For $  t\in (0,\tau)$ the Hopf--Lax
formula on $\mfC$ reads 
\begin{equation}
  \label{eq:67}
  \mathscr Q_t\zeta([x,r]):=\inf_{[x',r'] \in \mfC}\zeta([x',r'])+\frac
  1{2t}\sfd^2_\mfC\big([x,r],[x',r']\big).
\end{equation}
For $\xi_0$ satisfying $\xi_0\ge -\frac 1{2\tau}$
and $t\in (0,\tau)$  we set $ \zeta([x,r]):= \xi_0(x)r^2$ and find 
(cf.\ \cite[Thm.\,8.11]{LiMiSa18OETP}) 
\begin{equation}
  \label{eq:71}
\xi_t=\HopfLax t{\xi_0}\quad
  \Longleftrightarrow\quad
  \xi_t(x)r^2=\mathscr Q_t\zeta([x,r])\quad \text{for all } x\in \R^d;
\end{equation}
Moreover, if $\xi_0$ is lower semi-continuous
the infimum in \eqref{eq:67} is attained and 
we have
\begin{equation}
  \label{eq:120}
  \xi_t(x)r^2=\zeta([x',r']) {+}  \frac1{2t} 
  \sfd^2_\mfC\big([x,r],[x',r']\big) \Longleftrightarrow
  \begin{cases}
    x'\in \ArgM_t\xi_0(x) \text{ and}\\
    (1{+}2t\xi_0(x')) (r')^2=(1{-}2t\xi_t(x))r^2
  \end{cases}
\end{equation}
(where $[x',r']=\mfo$ if $r'=0$, corresponding to the case $1{-}2t\xi_t(x)=0$).
From \eqref{eq:67} and \eqref{eq:71} we also deduce the estimate
\begin{equation}
  \label{eq:129}
  \big(1{-} 2  t \psxi_t(x)\big)r^2+ \big(  1{+} 2 t\psxi_0(x')\big)(r')^2
  \geq 2rr'\tcos{\pit}{|x{-}x'|}
\end{equation}
for every $x,x'\in \R^d$ and $r,r'\ge0$.  Optimizing with respect to $r,r'$ we
find
\begin{equation}
  \big(1{-}  2  t \psxi_t(x)\big)\big(  1{+}
  2 t \psxi_0(x')\big)
  \geq \tcosq{\pit}{|x{-}x'|}
  \quad\text{for every }x,x'\in \R^d\label{eq:130}
\end{equation}
and arrive at the following characterization:
For all $x\in \R^d$ with $1{-}2t\xi_t(x)>0$ we have
\begin{equation}
  \label{eq:131}
  x'\in \ArgM_t\xi_0(x)\quad\Longleftrightarrow\quad
  \big(1{-} 2 t \psxi_t(x)\big)\big(  1{+}2t \psxi_0(x')\big)
= \tcosq{\pit}{|x{-}x'|}.
\end{equation}

To treat the factor of $r$ and $r'$ in \eqref{eq:129} efficiently, we
define the function 
\begin{equation}
  \label{eq:124}
  Z_t(u',u):=\frac{1{-}2t u}{1{+}2t u'}\ \text{ for }1{+}2tu',\ 
  1{-}2tu\geq 0 \quad \text{and }  
  Z_t(+\infty,u)\equiv 0. 
\end{equation}
Using \eqref{eq:131}, the optimal $r'$ in \eqref{eq:120} can now be
equivalently characterized by 
\begin{equation}
  \label{eq:146}
  (r')^2=Z_t(\xi_0(x'),\xi_t(x))=\frac{(1{-}2t\xi_t(x))^2}{\tcosq{\pit}{|x{-}x'|}}=
  (1{-}2t\xi_t(x))^2(1+\tan^2(|x{-}x'|)).
\end{equation}

The following result collects the properties of $\HopfLax t$ that will be
needed in the sequel.

\begin{proposition}[Properties of the generalized Hopf--Lax operator $\HopfLax t$]
\label{prop:HJ1}
Let $\xi_0:\OOmega\to [a,b]$ with $-1/2\le a\le b\le +\infty$ be lower
semi-continuous and set $\xi_t:=\HopfLax t{ \xi_0}$ for $t\in [0,1]$.
\begin{enumerate}[{\upshape(1)}]
\item \label{prop:HJ1.label1} \emph{Lower/upper bounds.} The functions $\xi_t$ are well defined
  and satisfy (cf.\:\eqref{eq:147} for $P_a$)
  \begin{equation}
    \label{eq:25}
    -\frac{1}{2(1{-}t)}\le  P_a(t)\le \xi_t\le P_b(t)\le \frac 1{2t}\quad
    \text{for every }t\in (0,1),\ x\in \R^d.
  \end{equation}
  Moreover, it holds that
  \begin{equation}
    \label{eq:168}
    \xi_0(x)=-1/2\quad\Leftrightarrow\quad
    \xi_t(x)=-\frac1{2(1{-}t)}.
  \end{equation}

\item \label{prop:HJ1.label2}
  \emph{Semi-concavity.} Setting
  $\Lambda_a(t) :=\frac{1}{t(1+2at)}\le \frac 1{t(1{-}t)}$ the
  functions $\xi_t$ are $\Lambda_a(t)$-Lipschitz and $\Lambda_a(t)$
  semi-concave, i.e.\ $ x\mapsto \xi_t(x)-\frac {\Lambda_a(t)}2|x|^2$ is
  concave.

\item\label{prop:HJ1.label3}
 \emph{Semigroup property.} For every $0\le s<t\le 1$ we have  
  \begin{equation}
    \label{eq:39}
    \xi_t=\HopfLax {t-s}\xi_s
  \end{equation}

\item \label{prop:HJ1.label4}
 \emph{Concatenation of optimal points.}  For $s,t$ with $0\le s<t< 1$ and
  $x\in \R^d$ we define the set-valued function $\ArgM_{t\to s}$ via
  $\ArgM_{t\to s}(x) := \ArgM_{t-s}\xi_s(x)$.  For all $0\le t_0<t_1<t_2<1$ and
  all $x_0,x_1,x_2 \in \R^d$ we have:
  \begin{equation}
    \label{eq:121}
    \begin{aligned}
      \text{If }& x_{1}\in \ArgM_{t_2\to t_1}(x_2) \text{ and } x_{0}\in
      \ArgM_{t_1\to t_0}(x_{1}), \quad \text{then } x_{0}\in
      \ArgM_{t_2 \to t_0}(x_2) \text{ and}
      \\[0.2em]
      &   Z_{t_2-t_0}(\xi_{t_0}(x_0),\xi_{t_2}(x_2))=Z_{t_1-t_0}(\xi_{t_0}(x_{0}),
      \xi_{t_1}(x_{1}))\, Z_{t_2-t_1}(\xi_{t_1}(x_{1}),\xi_{t_2}(x_2)).
    \end{aligned}
  \end{equation}

\item\label{prop:HJ1.label5}
 \emph{Geodesics on $\mfC$.}  If $0\le t_0<t_1<t_2<1$,
  $x_0\in \ArgM_{t_2\to t_0}(x_2)$,
  $r_0=Z_{t_2-t_0}(\xi_{t_0}(x_0),\xi_{t_2}(x))r_2$, and
  $[x_1,r_1]=\geo{\theta}{[x_0,r_0]}{[x_2,r_2]}$ for
  $\theta=\frac{t_1-t_0}{t_2-t_0}$, then $ x_1\in \ArgM_{t_2\to t_1}(x_2).$

\item\label{prop:HJ1.label6} 
\emph{Characterization of optimality.} For all $x,y\in \R^d$ and
  $0\le s<t< 1$ with $\tau:=t{-}s$ we have
    \begin{gather}
      \label{eq:118}
      (1{-}2\tau\xi_t(x))(1{+}2\tau\xi_s(y))\ge
      \tcosq{\pit}{|x{-}y|},\\
      \label{eq:135}
      y\in \ArgM_{t\to s}(x),\ \xi_t(x)<\frac1{2\tau}\ \Leftrightarrow\
      (1{-}2\tau\xi_t(x))(1{+}2\tau\xi_s(y))= \tcosq{\pit}{|x{-}y|}.
      \end{gather}
    \end{enumerate}
\end{proposition}

\begin{proof}
  \noindent
\underline{Assertion (\ref{prop:HJ1.label1}).}
  The first assertion follows by the monotonicity property of $\mathscr P_t$ and
  \eqref{eq:147}. Note that \eqref{eq:168}
  is a simple consequence of the property
  \[
    \frac{1} 
    {2t}
\Big(1-\frac{\tcosq{\pit}{|x{-}y|}}{1+
  2t\xi(y)}\Big)\ge -\frac1{2(1{-}t)}
\]
with equality if and only if $x=y$ and $\xi(y)=-1/2$.
  
  \noindent
\underline{Assertion (\ref{prop:HJ1.label2}).}
  It is sufficient to observe that for every $y\in \R^d$ 
  \begin{equation}
    \label{eq:78}
    x\mapsto \tcosq{\pit}{|x{-}y|}
    \text{ is $2$-Lipschitz,}
    \quad
    x\mapsto \tcosq{\pit}{|x{-}y|}-|x|^2
    \text{ is concave},
  \end{equation}
  so that
  \begin{equation}
    \label{eq:78bis}
    x\mapsto \frac 1{2t}\Big(1-\frac{\tcosq{\pit}{|x{-}y|}}{1+2t\xi_0(y)}\Big)\quad
    \text{is $\Lambda_a(t)$-Lipschitz}
  \end{equation}
  and
  \begin{equation}
    \label{eq:78ter}
    x\mapsto \frac 1{2t}\Big(1-
   \frac{\tcosq{\pit}{|x{-}y|}}{1+2t\xi_0(y)}\Big) 
   -\frac{\Lambda_a(t)}2|x|^2\quad \text{is concave}.
  \end{equation}
  \noindent
\underline{Assertion (\ref{prop:HJ1.label3}).}
   If $t<1$ the semigroup property for $\mathscr P_t$ 
  can be derived by the link with the Hopf–Lax semigroup in $\mfC$
  given by \eqref{eq:71}
  and the fact that $(\mfC,\sfd_{\pi,\mfC})$ is a geodesic space.
  The case $t=1$ follows by approximation and \eqref{eq:63}.

  \noindent
\underline{Assertion (\ref{prop:HJ1.label4}).} We set $\tau_0:=t_1-t_0$, $\tau_1:=t_2-t_1$, $r>0$,
\[
  r_1=Z_{\tau_1}(\xi_{t_1}(x_1),\xi_{t_2}(x))r,\quad
  r_0=Z_{\tau_0}(\xi_{t_0}(x_0),\xi_{t_1}(x_1))r_1
\]
and use \eqref{eq:120} and \eqref{eq:39}:
\begin{align*}
  \xi_{t_2}(x)r^2
  &=
  \xi_{t_1}r_1^2+\frac{1}{2\tau_1}\sfd_{\pi,\mfC}^2
  ([x_1,r_1], [x,r])
  \\
  &= 
  \xi_{t_0}r_0^2+
  \frac{1}{2\tau_0}\sfd_{\pi,\mfC}^2
  ([x_0,r_0], [x_1,r_1])
  +\frac{1}{2\tau_1}\sfd_{\pi,\mfC}^2
  ([x_1,r_1], [x,r]).
\end{align*}
On the other hand
\begin{equation}
  \label{eq:133}
  \xi_{t_2}(x)r^2\le \xi_{t_0}r_0^2+
  \frac{1}{2\tau}
  \sfd_{\pi,\mfC}^2
  ([x_0,r_0], [x,r])
\end{equation}
so that we obtain
\begin{equation}
  \label{eq:132}
  \frac{1}{2\tau_0}\sfd_{\pi,\mfC}^2
  ([x_0,r_0], [x_1,r_1])  +\frac{1}{2\tau_1}\sfd_{\pi,\mfC}^2
  ([x_1,r_1], [x,r])\le   \frac{1}{2\tau}  \sfd_{\pi,\mfC}^2
  ([x_0,r_0], [x,r]);
\end{equation}
since $\tau=\tau_0+\tau_1$ the opposite inequality always hold in
\eqref{eq:132}, and we deduce the equality, which implies that the equality
holds in \eqref{eq:133} as well, showing \eqref{eq:121}
thanks to \eqref{eq:120}.
    
\noindent
\underline{Assertion (\ref{prop:HJ1.label5}).} We can argue as in the previous assertion, starting from the characterization
of $x_0,r_0$
\begin{equation}
  \label{eq:133bis}
  \xi_{t_2}(x)r^2= \xi_{t_0}r_0^2+  \frac{1}{2\tau}
  \sfd_{\pi,\mfC}^2  ([x_0,r_0], [x,r])
\end{equation}
and using the identity along the geodesic in $\mfC$ connecting $[x_0,r_0]$ to
$[x,r]$, namely 
\begin{equation}
  \label{eq:132bis}
  \frac{1}{2\tau_0}\sfd_{\pi,\mfC}^2  ([x_0,r_0], [x_1,r_1])
  +\frac{1}{2\tau_1}\sfd_{\pi,\mfC}^2  ([x_1,r_1], [x,r])= 
  \frac{1}{2\tau}   \sfd_{\pi,\mfC}^2   ([x_0,r_0], [x,r]).
\end{equation}
\noindent
\underline{Assertion (\ref{prop:HJ1.label6}).} The final assertion follows from \eqref{eq:130} and \eqref{eq:131}.
\end{proof}

\subsection{Backward generalized Hopf--Lax flow and contact sets}
\label{su:BackwardGHL}
 
Let us now consider the backward version of the generalized Hopf--Lax
semigroup. By the simple structure of the generalized Hamilton--Jacobi equation
\eqref{eq:genHJeqn}, we immediately see that time reversal leads to the same
effect  as the sign reversal $\xi \leadsto -\xi$. Hence, the
backward semigroup $\RopfLax t$ is defined for $\bar\xi$
with $\bar\xi \leq 1/(2\tau)$ via  
\begin{equation}
  \label{eq:44}
  \RopfLax t\bar \xi:=-\HopfLax t{(-\bar\xi)} \quad \text{for }
   t\in (0,\tau]. 
\end{equation}
The corresponding properties of $\RopfLax t$ follow easily from Proposition
\ref{prop:HJ1}, but observe that we use $\bar \xi_t = \RopfLax{1-t}{
\bar\xi_1}$ to go backward in time.  

\begin{corollary}[Properties of $\RopfLax t$] 
\label{cor:HJ1}
Let $\bar \xi_1:\OOmega\to[-\bar b,-\bar a]$ with
$-\infty\le -\bar b\le -\bar a\le 1/2$ be upper semi-continuous and set
\begin{equation}
  \label{eq:37}
  \bar \xi_t:=\RopfLax {1-t}{\bar \xi_1}\quad \text{for } t\in [0,1].
\end{equation}
\begin{enumerate}[{\upshape(1)}]

\item \emph{Lower/upper bounds.}  The functions $\bar \xi_t$ are well-defined and satisfy
  \begin{equation}
    \label{eq:25bis}
    -\frac 1{2t}\le P_{\,\bar b}(1{-}t)\le \bar\xi_t\le P_{\bar a}(1{-}t)
    \le \frac1{2(1{-}t)}\quad
    \text{for all }t\in (0,1),\ x\in \R^d.
  \end{equation}
  Moreover, we have the equivalence
  \begin{equation}
    \label{eq:168bis}
    \bar\xi_1(x)=1/2\quad\Leftrightarrow\quad
    \bar\xi_t(x)=\frac1{2t}.
  \end{equation}

\item \emph{Semi-convexity.}  The functions $\bar\xi_t$ are
  $\Lambda_{\bar a}(1{-}t)$-Lipschitz and $\Lambda_{\bar a}(1{-}t)$
  semi-convex, i.e.\ $x\mapsto \bar\xi_t(x)+\frac {\Lambda_{\bar a}(1{-}t)}2|x|^2$
  is convex (cf.\ Proposition {\upshape\ref{prop:HJ1}(\ref{prop:HJ1.label2})} for $\Lambda_a$). 

\item \emph{Time-reversed semigroup property.} For every $0\le s<t\le 1$
  we have
  \begin{equation}
    \label{eq:39bis}
    \bar\xi_s=\RopfLax{t-s}\bar \xi_t.
  \end{equation}

\item \emph{Concatenation of optimal points.} Setting $
 \ol\ArgM_{s\shortleftarrow t}(x) :=
\ArgM_{t-s}(-\bar \xi_t)(x)$ for every $0< s<t\le 1$ and $x\in \R^d$, the
set-valued function $\ol\ArgM_{s \shortleftarrow t}$
  satisfies the concatenation property for
$0<t_0<t_1<t_2\le 1$ and $x_0,x_1,x_2\in \R^d$:
  \begin{equation}
    \label{eq:121bis}
    \begin{aligned} \text{If }& x_{1}\in \ol\ArgM_{t_0\shortleftarrow t_1}(x_0) \text{ and }
x_{2}\in \ol\ArgM_{t_1\shortleftarrow t_2}(x_1), \quad \text{then } x_{2}\in
\ol\ArgM_{t_0\shortleftarrow t_2}(x_0)\text{ and} 
\\ &Z_{t_2-t_0}(\bar\xi_{t_0}(x_0),\bar\xi_{t_2}(x_2))=
Z_{t_1-t_0}(\bar\xi_{t_0}(x_{0}), \bar\xi_{t_1}(x_1))\cdot
Z_{t_2-t_1}(\bar\xi_{t_1}(x_{1}),\bar\xi_{t_2}(x_2)).
    \end{aligned}
  \end{equation}

  
\item \emph{Characterization of optimality.}  For all $x,y\in \R^d$ and
  $0< s<t\le 1$ with $\tau:=t-s$
  \begin{gather}
    \label{eq:118bis}
    (1{-}2\tau\bar \xi_t(x))\,(1{+}2\tau\bar\xi_s(y)) \geq
    \tcosq{\pit}{|x{-}y|},\\
    \label{eq:135bis}
    x\in \ArgM_{s\to t}(y),\ \bar\xi_s(y)>-\frac 1{2\tau}\ \Leftrightarrow\
    (1{-}2\tau\bar\xi_t(x))\,(1{+}2\tau\bar\xi_s(y))= \tcosq{\pit}{|x{-}y|}.
  \end{gather}
\end{enumerate}
\end{corollary}
\begin{proof}
We just observe that the second statement in \eqref{eq:121bis} follows by the
corresponding statement in \eqref{eq:121} which now reads as
\begin{equation}
  \label{eq:122pre}
  Z_{t_2-t_0}(-\bar\xi_{t_2}(x_2),-\bar\xi_{t_0}(x_0))=
  Z_{t_1-t_0}(-\xi_{t_1}(x_{1}),-\xi_{t_0}(x_0)      )
  \cdot Z_{t_2-t_1}(-\bar \xi_{t_2}(x_2),-\bar\xi_{t_1}(x_{1})),
\end{equation}
and the property $Z_{\tau}(-u',-u) =Z_{-\tau}(u',u)
=Z^{-1}_\tau(u,u')$. Equations \eqref{eq:118bis} and \eqref{eq:135bis} follow
by \eqref{eq:118} and \eqref{eq:135} changing $\xi_s(y)$ with $-\bar \xi_t(x)$
and $\xi_t(x)$ with $-\bar\xi_s(y)$.
\end{proof}

We are now in the position to compare the forward solution $\xi_t$ and the
backward solution $\bar\xi_t$. The main philosophy is that in general we only
have $\RopfLax{t}{} \HopfLax{t}{\xi_0} \leq \xi_0$ (cf.\ \eqref{eq:37bis}
below), but equality holds $\mu_t$-a.e.\ if $(\xi_0,\HopfLax 1\xi_0)$ is an
optimal pair.  In the following result, we still stay in the general case
comparing arbitrary forward solutions $\xi_t=\HopfLax t \xi_0$ and backward
solutions $\bar\xi_t=\RopfLax{1-t}\bar\xi_1$ only assuming
$\xi_1 \geq \bar\xi_1$.  Along the contact set $\Xi_t$ where $\xi_t$ and
$\bar\xi_t$ coincide, we can then derive differentiability and optimality properties of
$\xi_t$ and $\bar\xi_t$.

\begin{theorem}[Contact set $\Xi_t$]
\label{thm:HJ1}
Let $\xi_0:\OOmega\to [a,+\infty]$ be l.s.c.\ with $a\geq -1/2$ and
$\bar \xi_1:\OOmega\to [-\infty, -\bar a]$ u.s.c.\ with $ \bar a \leq
1/2$. Assume $ \HopfLax 1{\xi_0} \geq \bar \xi_1$ and set
\begin{equation}
  \label{eq:37bis}
  \xi_t:=\HopfLax t{ \xi_0} \ \text{ and } \ 
  \bar \xi_t:=\RopfLax {1-t}{\bar \xi_1}\quad
  \text{for } t\in [0,1].
\end{equation}
Then, the following assertions hold: 
\begin{enumerate}[{\upshape(1)}]
\item\label{thm:HJ1.label1} 
  For every $t\in [0,1]$ we have $ \xi_t\geq \bar\xi_t$ and the
  \emph{contact set} 
  \begin{equation}
    \label{eq:116}
    \Xi_t:=\big\{x\in \R^d:\bar \xi_t(t)=\xi_t(x)\big\}
    \quad\text{is closed}.
  \end{equation}
\item\label{thm:HJ1.label2} 
   For every $t\in (0,1)$ and $x\in \Xi_t$ there exists a unique
  $p=\GRAD_t(x)$  satisfying 
  \begin{equation}
    \label{eq:54}
    \xi_t(y)-\xi_t(x)-\frac1{2}\Lambda_{\bar a}(1{-}t)|x{-}y|^2\le
    \langle p,y-x\rangle\le
    \bar\xi_t(y)- \MAT \bar\xi_t \EEE (x)+\frac1{2}\Lambda_{a}(t)|x{-}y|^2
  \end{equation}
  so that in particular $\xi_t$ and $\bar \xi_t$ are differentiable at $x$ with
  gradient $\GRAD_t(x)$ (cf.\ Proposition {\upshape\ref{prop:HJ1}(\ref{prop:HJ1.label2})} for $\Lambda_a$).

\item\label{thm:HJ1.label3} 
 The map $x\mapsto \GRAD_t(x)$ is bounded and $C(t)$-Lipschitz with
  $C(t)\le 2(\Lambda_a(t)+\Lambda_{\bar a}(1{-}t))\le \frac{4}{t(1{-}t)}$ on
  $\Xi_t$. Moreover, the sets
  \begin{equation}
    \label{eq:117}
    \begin{aligned}
      \Xi^-_t:={}&\Bigset{x\in \R^d}{\xi_0=-\frac 12 }
         =\Bigset{x\in \Xi_t }{ \xi_t=\frac {-1}{2(1{-}t)} }
      \\
      \Xi^+_t:={}&\Bigset{x\in \R^d }{ \bar \xi_1=\frac12 } =\Bigset{x\in
      \Xi_t }{\bar\xi_t=\frac 1{2t} }
    \end{aligned}
  \end{equation}
      are independent of $t$, are contained in $\Xi_t$ for every $t\in [0,1]$, and
      \GGG the critical  set 
      $\Xi^0_t:=\{x\in \Xi_t:\GRAD_t(x)=0\}$ 
      of $\GRAD_t$ contains $\Xi^\pm_t$: \EEE
  \begin{equation}
    \label{eq:112}
   \GGG  \Xi^0_t\supset \Xi^-_t\cup\Xi^+_t
   \quad\text{for every }t\in (0,1).
  \end{equation}
\item\label{thm:HJ1.label4}  Let $s\in (0,1)$, $t\in [0,1]$, and $\tau:=t-s\neq 0$.  Then, for
  every $x_s\in \Xi_s$ with $1{+}2\tau\xi_s(x_s)>0$ the set
  $\ArgM_{s\to t}(x_s)$ consists of a unique element $x_t= :\bfT_{s\to t}(x_s)$
  satisfying
  \begin{equation}
    \label{eq:136}
    \begin{aligned}
      x_t&\in \Xi_t\ \text{and } x_s\in \ol\ArgM_{t \shortleftarrow
        s}(x_t),
      \\
      x_t&=\bfT_{s\to t}(x)=x+\bfarctan\Big(\frac{\tau
        \GRAD_s(x)}{1{+}2\tau\xi_s(x)}\Big), \\
      \big(1&{-}2\tau\xi_t(\bfT_{s\to t}(x))\big)\,\big(1{+}2\tau\xi_s(x)\big)
      = \tcosq\pit{\bigup|x{-}\bfT_{s\to t}(x)|}.
    \end{aligned}
  \end{equation}
\item \label{th:HJ.label5} For every $x\in 
\GGG \Xi^0_s\supset \Xi_s^{\pm}$ 
\EEE we have
  $\bfT_{s\to t}(x)=x$ (and thus we set $\bfT_{s\to t}(x):=x$ also for $t=0$ or
  $t=1$). Let $s\in (0,1)$ and define $\bfT_{s\to s}(x)=x$, then for all
  $x\in \Xi_s$ the mappings
  $t\mapsto \bfT_{s\to t}(x)$ are  analytic in $[0,1]$. \\
  For $s,t\in (0,1)$ the mappings $\bfT_{s\to t}:\Xi_s \to \R^d$ are 
  Lipschitz.  If $t=0$ (resp.\,$t=1$) then $\bfT_{s\to t}$ is locally Lipschitz
  in $\Xi_s\setminus \Xi_s^+$ (resp.\,in $\Xi_s\setminus \Xi_s^-$).

\item \label{th:HJ.label6} Setting
  \begin{equation}
    q^2_{s\to t}(x):=\frac{1+2\tau\xi_s(x)} { 1 {-}2\tau\xi_t 
      (\bfT_{s\to t}(x))}= (1 {+} 2\tau\xi_s(x))^2+\tau^2 |\GRAD_s(x)|^2
    \label{eq:137}
  \end{equation}
  for every $x\in \Xi_s$, the map $t\mapsto q_{s\to t}(x)$ is analytic in
  $[0,1]$, $q_{t\to s}$ is bounded and Lipschitz with respect to $x$, and
  $q_{s\to t}(x)>0$ for $t\in (0,1)$ or $t=0$ and $x\not \in \Xi_s^+$
  (resp.\,$t=1$ and $x\not\in \Xi_s^-$). Moreover, 
  $  q_{s\to t}(x)=1+2(t{-}s)\xi_s(x)\text{ for  }x\in \Xi_s^\pm$.

\item\label{thm:HJ1.label7}  For all $t_0, t_1\in (0,1)$, $t_2\in [0,1]$, the maps
  $\bfT_{t_i\to t_j}$ are Lipschitz on $\Xi_{t_i}$ for $i\in \{0,1\}$, and we have
  \begin{equation}
    \label{eq:26}
    \bfT_{t_1\to t_2}\circ \bfT_{t_0\to t_1}=\bfT_{t_0\to t_2},\quad
    q_{t_1\to t_2}(\bfT_{t_0\to t_1}(x))\cdot q_{t_0\to t_1}(x)=q_{t_0\to t_2}(x).
  \end{equation}
\end{enumerate}
\end{theorem}
\begin{proof}
  \noindent
\underline{Assertion (\ref{thm:HJ1.label1}).} The inequality
\begin{equation}
  \label{eq:81}
  \xi_s\ge \RopfLax{t-s}{\Big(\HopfLax{t-s}{\xi_s}\Big)}
  =
  \RopfLax{t-s}\xi_t\quad\text{for }  0<s<t<1
\end{equation}
can be derived by the link with the Hopf–Lax semigroup in $\mfC$ given by
\cite[Theorem 8.11]{LiMiSa18OETP} and arguing as in
\cite[Thm. 7.36]{Vill09OTON}.  We prove it by a direct computation as follows:
Set $\tau=t{-}s$, observe that
$\inf \HopfLax{\tau}{\xi_s}=\inf \xi_t\le \frac{1}{2t}<\frac1{2\tau}$, and use
$\xi_t = \HopfLax \tau \xi_s$ to obtain 
\begin{equation}
  \label{eq:87}
  \frac1{1{-}2\tau\xi_t(y)}=  \inf_{z\in B_{\pit}(y)}
  \frac{1+2\tau\xi_s(z)}{\tcosq{\pit}{|y{-}z|}}
  \leq  \frac{1+2\tau\xi_s(x)}{\tcosq{\pit}{|y{-}x|}}\quad
  \text{if }|x{-}y|<\pi/2.
\end{equation}
With this estimate, we find
\begin{align*}
  \RopfLax\tau{\xi_t}(x)
  & \overset{\text{def}}=
  \sup_{y\in B_{\pit}(x)}\frac
  1{2\tau}\Big(\frac{\tcosq{\pit}{|x{-}y|}}{1-2\tau\xi_s(y)}-1\Big)
  \\& \overset{\text{\eqref{eq:87}}} \leq 
  \sup_{y\in B_{\pit}(x)}
  \frac
  1{2\tau}\Big(\frac{\tcosq{\pit}{|x{-}y|}}{\tcosq{\pit}{|x{-}y|}}
  \,(1{+}2\tau\xi_s(x))-1\Big) \ = \ \xi_s(x).
\end{align*}
Using $\xi_t\ge \xi_1\ge \bar \xi_1$ we thus get \eqref{eq:81}.  
Passing to the limit as $t\uparrow 1$ in \eqref{eq:81}, we arrive at
$\xi_s\ge \RopfLax{1-s}{\bar\xi_1}=\bar \xi_s$.

The closedness of $\Xi_t$ follows from the semi-continuities of $\xi_t$ and
$\bar\xi_t$ and the estimate $\xi_t\geq \bar\xi_t$. Indeed, assume $y_k\to y$
with $y_k\in \Xi_t$, then we have $y\in \Xi_t$ because of 
\[ 
\xi_t(y) \overset{\text{l.s.c.}}\leq  \liminf_{k\to \infty} \xi_t(y_k) =
\liminf_{k\to \infty} \bar\xi_t(y_k) \leq \limsup_{k\to \infty} \bar\xi_t(y_k) 
\overset{\text{u.s.c.}}\leq \bar\xi_t(y) \leq \xi_t(y).
\]

\noindent
\underline{Assertion (\ref{thm:HJ1.label2}).} Let us fix $x\in \Xi_t$, $\Lambda:=\Lambda_a(t)$ and
$\bar \Lambda:=\Lambda_{\bar a}(1{-}t)$, and let $p$ (resp.\,$p'$) be an element of the
superdifferential of $x\mapsto \xi_t(x)-\frac12 \Lambda|x|^2$ (resp.\,of the
subdifferential of $x\mapsto \bar\xi_t(x)+\frac12\bar \Lambda|x|^2$). The
superdifferential (subdifferential) is not empty, since the function is
concave (convex) and finite everywhere.  For every $x,y\in \R^d$ with $x\in
\Xi_t$ we have 
\begin{align*}
    \langle p,y{-}x\rangle
    &\geq
      \xi_t(y)-\xi_t(x)-\frac12 \Lambda|x{-}y|^2 \text{ and }
    \langle p',y{-}x\rangle \leq \bar \xi_t(y)-\bar\xi_t(x)+\frac12
    \bar \Lambda|x{-}y|^2.
\end{align*}
Subtracting the two inequalities and using $\xi_t(x)=\bar\xi_t(x)$ and  $ \bar\xi_t(y)\le \xi_t(y)$ yields
\[
  \langle p'-p,y-x\rangle\le \bar\xi_t(y)-\xi_t(y) +\frac 12(\Lambda{+}\bar \Lambda) |y{-}x|^2 \leq \frac 12(\Lambda{+}\bar \Lambda) |y{-}x|^2\quad\text{for every   }y\in \R^d,
\]
so \EEE that $p=p'$ is uniquely determined and \eqref{eq:54} holds.
\medskip

\noindent
\underline{Assertion (\ref{thm:HJ1.label3}).} The fact that $\Xi^\pm_t$ are independent of $t$ and
contained in $\Xi_t$ follows from \eqref{eq:168} and \eqref{eq:168bis}. 
Moreover, \eqref{eq:112} follows easily since $\xi_t$ takes its minimum at $\Xi^-_t$ and
its maximum at $\Xi^+_t$.

Let us now fix $t\in (0,1)$, $x_0,x_1\in \Xi_t$, $p_i=\GRAD_t(x_i)+\bar \Lambda x_i$, and
set $\bar \zeta(x):=\bar \xi_t(x) + \frac 12 \bar \Lambda
|x|^2$, $\zeta(x):= \xi_t(x) +\frac 12 \bar \Lambda
|x|^2$.  Notice that $\bar\zeta(x)$ is convex and $\zeta(x)$ is $C=\Lambda+\bar \Lambda$
semi-concave with $\bar \zeta(x)\le \zeta(x)$.  We get
\begin{align*}
  \bar\zeta(x_0)
  &\le
    \bar\zeta(x)-\langle p_0,x-x_0\rangle
    \le \zeta(x)
    -\langle p_0,x-x_0\rangle
  \\&  \le
    \zeta(x_1)+\langle p_1,x-x_1\rangle
    -\langle p_0,x-x_0\rangle
    +
    \frac C2 |x{-}x_1|^2
    \\&=
  \bar\zeta(x_1)
    +\langle p_1-p_0,x-x_1\rangle-\langle
  p_0,x_1-x_0\rangle+ \frac C2 |x{-}x_1|^2.
\end{align*}
Minimizing with respect to $x$ 
we find
$  \bar\zeta(x_0)
  \le
   \bar\zeta(x_1)
    -\langle
  p_0,x_1{-}x_0\rangle-
  \frac 1{4C}|p_1{-}p_0|^2$.
Inverting the role of $x_0$ and $x_1$ and summing up gives 
$  \frac 1{2C}|p_1{-}p_0|^2\le
  \langle
  p_1{-}p_0,x_1{-}x_0\rangle $
and therefore
\begin{equation}
  \label{eq:100}
  |p_1{-}p_0| \leq 2C \,|x_1{-}x_0|.
\end{equation}
The boundedness of $\GRAD_t$ on $\Xi_t$ follows by the fact that $\xi_t$ is
Lipschitz.\smallskip

\noindent
\underline{Assertion (\ref{thm:HJ1.label4}).} Let us first consider the case $s>t$ with $\tau:=s{-}t$ and let
$y\in \ArgM_{s\to t}(x)$.  If $\xi_s(x)>-\frac 1{2\tau}$ then $y$ satisfies the
identity \eqref{eq:135}.  Since $\bar\xi_t(y)\le \xi_t(y)$, \eqref{eq:118bis}
and $\bar\xi_s(x)=\xi_s(x)$ yields $\bar \xi_t(y)=\xi_t(y)$ so that
$y\in \Xi_t$ as well with $x\in \ArgM_{t\shortleftarrow s}(y)$ since
$\bar \xi_t(y) \ge -\frac{1}{2(1{-}t)}>-\frac{1}{2\tau}$.

Since the function
$  x'\mapsto (1{+}2\tau \xi_s(x'))\,(1{-}2\tau\xi_t(y))
  -\tcosq{\pit}{|x'{-}y|}$
has a global minimizer at $x$, we arrive at the Euler--Lagrange equations
\[
  2\tau (1{-}2\tau\xi_t(y))\, \GRAD_s(x) +2\tcos{\pit}{|x{-}y|}\,\bfsin(x{-}y)=0
\]
Since we can assume $|x{-}y|<\pi/2$ we obtain
\begin{equation}
  \label{eq:134}
  x-y=-\bfarctan\Big(\frac{\tau \GRAD_s(x)}{ 1 {+} 2\tau\xi_s(x)}\Big),
\end{equation}
which characterizes $y$ uniquely and establishes \eqref{eq:136}. 

The case $t>s$ follows by the same arguments.

\noindent
\underline{Assertion (\ref{th:HJ.label5}).} This assertion is an immediate consequence of \eqref{eq:136} and \eqref{eq:112}.

\noindent
\underline{Assertion (\ref{th:HJ.label6}).} The claims are simple consequences of the identity \eqref{eq:146} and
the definition of $q_{s\to t}$ of \eqref{eq:137}.

\noindent
\underline{Assertion (\ref{thm:HJ1.label7}).} The final assertion follows by \eqref{eq:121} (and the corresponding
\eqref{eq:121bis}).
\end{proof}

\begin{remark}[Strongly reduced pairs]
  It is worth noticing that if $\inf \xi_0>-\frac 12$ and
  $\sup\bar \xi_1<\frac 12$, then the
  sets $\Xi^\pm_t$ in \eqref{eq:117} are empty
  and many properties of $\xi_t$, $\bfT_{s\to t}$ and $q_{s\to t}$
  become
  considerably simpler. This situation is, e.g., the case
  of the solution induced by a strongly reduced pair with compact support, see
  Theorem~\ref{thm:extra-regularity}.
\end{remark}

We close this subsection by giving a small example for $\xi_t$ and
$\bar\xi_t$ and their contact set $\Xi_t$ 
derived from an optimal pair $(\xi_0,\bar\xi_1)$ for the transport
between two Dirac measures. 

\begin{example}[The contact set for two Dirac measures]
\slshape\label{ex:ContactDirac} For points $z_0,z_1\in \R^d$ and $r_0,r_1>0$ we
consider the Dirac measures $\mu_j= r_j^2 \delta_{z_j}$. We have 
\[
\HK^2(\mu_0,\mu_1)=r_0^2 + r_1^2 - 2r_0r_1 \tcos\pit\varrho  \quad 
\text{with } \varrho= |z_1{-}z_0|,
\]
and all geodesic curves are known, see \cite[Sec.\,5.2]{LiMiSa16OTCR}. For
$\varrho< \pit$ we have a unique geodesic $\mu_t= r(t)^2 \delta_{z(t)}$ defined
by transport, and
for $ \varrho>\pit$ the unique geodesic $\mu_t=(1{-}t)^2 r_0^2\delta_{z_0} +
t^2r_1^2 \delta_{z_1}$ consists of \AAA growth \EEE (annihilation and decay) only. For
$\varrho=\pit$ there is an infinite-dimensional convex set of geodesics, and we
will see that this property is also reflected by  a larger contact set. 

Using the simple one-point supports of $\mu_j$ it is easy to calculate the 
optimal potentials and the transport plan $\bfeta$ in Theorem~\ref{thm:optimality.cond}(ii). We obtain 
\[
s_0:=\sigma_0(z_0)=\frac{r_1}{r_0}\,\tcos\pit\varrho, \quad 
s_1:=\sigma_1(z_1)=\frac{r_0}{r_1}\,\tcos\pit\varrho, \quad
\bfeta= r_0r_1\tcos\pit\varrho\,\delta_{(z_0,z_1)}.  
\] 
Thus, we will distinguish the case $\tcos\pit\varrho>0$ and
$\tcos\pit\varrho=0$.\medskip 

\underline{Case $\varrho<\pit$:} By \eqref{eq:19xi} the optimal pair
$(\xi_0,\bar\xi_1)$ reads
\[
\xi_0(x) = \begin{cases} \frac{s_0{-}1}2 &\text{for } x=z_0, \\
+\infty &\text{for } x\neq z_0;  \end{cases}
\quad \text{and} \quad 
\xi_1(x) = \begin{cases} \frac{1{-}s_1}2 &\text{for } x=z_1, \\
-\infty &\text{for } x\neq z_1.  \end{cases}
\]
From these identities, we obtain the forward and backward solutions $\xi_t = \HopfLax
t\xi_0$ and $\bar\xi_t = \RopfLax{1{-}t}\bar\xi_1$:
\begin{equation}
\label{eq:xi.barxi}
\xi_t(x) =\frac{1{-}t{+}ts_0- \tcosq\pit{|x{-}z_0|}} {2\,t\,(1{-}t{+}t s_0)} 
 \quad \text{and} \quad
\bar\xi_t(x) =\frac{\tcosq\pit{|x{-}z_1|} - t {-}(1{-}t) s_1} 
  {2\,(1{-}t)\,(t{+}(1{-}t) s_1)} . 
\end{equation}
The following optimality conditions can be checked by direct computation:
\begin{align*}
\text{(a) } \ &\xi_0\geq \bar\xi_0  \text{ and } \xi_1 \geq \bar\xi_1 \text{ \ on }\R^d 
\\
\text{(b) } \ &  \xi_0 = \bar\xi_0 \ \ \mu_0\text{-a.e.}\quad \text{and} \quad 
 \xi_1=\bar\xi_1 \ \ \mu_1\text{-a.e.} 
\end{align*}
As $\xi_0(x)=+\infty$ for $x\neq z_0$ and $\xi_1(x)=-\infty$ for $x\neq z_1$
statement (a) follows from (b). For (b) observe
\[
\bar\xi_0(z_0)= \frac{\cos^2\varrho -s_1}{2s_1} = \frac{ \cos^2\varrho
  - (r_0/r_1) \cos\varrho}{2(r_0/r_1)\cos\varrho} 
= \frac12\big(\frac{r_1}{r_0}\cos\varrho -1\big)=\frac{s_0{-}1}2= \xi_0(z_0).
\] 
Similarly, $\xi_1(z_1)=\bar\xi_1(z_1)$ follows, which provides a first result
on the contact sets $\Xi_t:=\bigset{x\in \R^d}{\xi_t(x)=\bar\xi_t(x)}$, 
namely $\Xi_0=\{z_0\}$ and $\Xi_1=\{z_1\}$. 
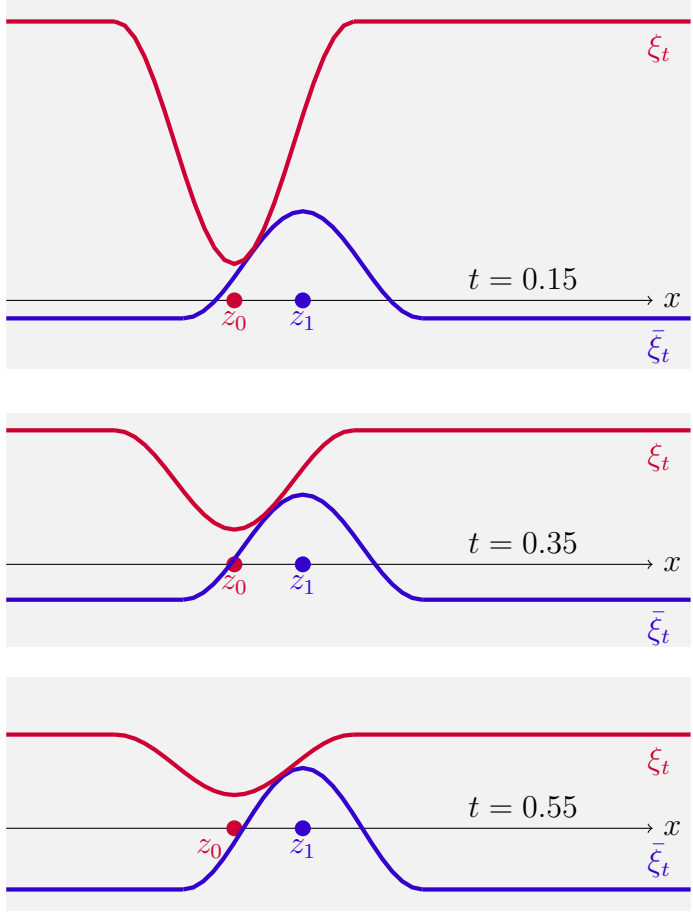
\begin{figure}
{\begin{tikzpicture}
\def\PITT{1.5708}
\def\rrrr{0.9} 
\def\dddd{2.0} 
\def\tttt{0.15}
\def\aaaa{3.0} 
\def\aaab{3.35}

\fill[color=gray!10] (-3,\aaaa-0.9) rectangle (6,\aaaa+4.0);
\draw[->] (-3,\aaaa)-- node[pos=0.8, above] {$t=\tttt$} (5.5,\aaaa) node[right]{$x$};
\fill[color=red!80!blue] (0,\aaaa) circle (3pt) node[below] {$z_0$};
\fill[color=blue!80!red] (\rrrr,\aaaa) circle (3pt) node[below] {$z_1$};

\draw[blue!80!red, ultra thick, domain=-3:\rrrr-\PITT] 
  plot (\x, {\aaab+(0.0 -\tttt -(1-\tttt)*cos(\rrrr r)/\dddd)/
    (2*(1-\tttt)*(\tttt + (1-\tttt)*cos(\rrrr r)/\dddd ) ) } );

\draw[blue!80!red, ultra thick, domain=\rrrr-\PITT:\rrrr+\PITT] 
  plot (\x, {\aaab+((cos((\x-\rrrr) r))^2 -\tttt -(1-\tttt)*cos(\rrrr r)/\dddd)/
    (2*(1-\tttt)*(\tttt + (1-\tttt)*cos(\rrrr r)/\dddd ) ) } );

\draw[blue!80!red, ultra thick, domain=\rrrr+\PITT:6] 
  plot (\x, {\aaab+(0.0 -\tttt -(1-\tttt)*cos(\rrrr r)/\dddd)/
    (2*(1-\tttt)*(\tttt + (1-\tttt)*cos(\rrrr r)/\dddd ) ) } ) 
    node[below] {$\bar\xi_t$\qquad\mbox{}}   ;

\draw[red!80!blue, ultra thick, domain=-3:-\PITT] 
  plot (\x, {\aaaa+(1.0 -\tttt +\tttt*cos(\rrrr)*\dddd)/
    (2*\tttt*(1-\tttt + \tttt*cos(\rrrr r)*\dddd ) ) } );

\draw[red!80!blue, ultra thick, domain=-\PITT:\PITT] 
  plot (\x, {\aaaa+(1-\tttt +\tttt*cos(\rrrr)*\dddd-(cos(\x r))^2)/
    (2*\tttt*(1-\tttt + \tttt*cos(\rrrr r)*\dddd ) ) } );

\draw[red!80!blue, ultra thick, domain=\PITT:6] 
  plot (\x, {\aaaa+(1-\tttt +\tttt*cos(\rrrr)*\dddd)/
    (2*\tttt*(1-\tttt + \tttt*cos(\rrrr r)*\dddd ) ) } )
   node[below]{$\xi_t$\qquad\mbox{}};
\def\tttt{0.35}
\def\aaaa{-0.5}
\def\aaab{-0.2}

\fill[color=gray!10] (-3,\aaaa-1.1) rectangle (6,\aaaa+2.0);
\draw[->] (-3,\aaaa)-- node[pos=0.8, above] {$t=\tttt$} (5.5,\aaaa) node[right]{$x$};
\fill[color=red!80!blue] (0,\aaaa) circle (3pt) node[below] {$z_0$};
\fill[color=blue!80!red] (\rrrr,\aaaa) circle (3pt) node[below] {$z_1$};

\draw[blue!80!red, ultra thick, domain=-3:\rrrr-\PITT] 
  plot (\x, {\aaab+(0.0 -\tttt -(1-\tttt)*cos(\rrrr r)/\dddd)/
    (2*(1-\tttt)*(\tttt + (1-\tttt)*cos(\rrrr r)/\dddd ) ) } );

\draw[blue!80!red, ultra thick, domain=\rrrr-\PITT:\rrrr+\PITT] 
  plot (\x, {\aaab+((cos((\x-\rrrr) r))^2 -\tttt -(1-\tttt)*cos(\rrrr r)/\dddd)/
    (2*(1-\tttt)*(\tttt + (1-\tttt)*cos(\rrrr r)/\dddd ) ) } );

\draw[blue!80!red, ultra thick, domain=\rrrr+\PITT:6] 
  plot (\x, {\aaab+(0.0 -\tttt -(1-\tttt)*cos(\rrrr r)/\dddd)/
    (2*(1-\tttt)*(\tttt + (1-\tttt)*cos(\rrrr r)/\dddd ) ) } ) 
    node[below] {$\bar\xi_t$\qquad\mbox{}}   ;

\draw[red!80!blue, ultra thick, domain=-3:-\PITT] 
  plot (\x, {\aaaa+(1.0 -\tttt +\tttt*cos(\rrrr)*\dddd)/
    (2*\tttt*(1-\tttt + \tttt*cos(\rrrr r)*\dddd ) ) } );  

\draw[red!80!blue, ultra thick, domain=-\PITT:\PITT] 
  plot (\x, {\aaaa+(1-\tttt +\tttt*cos(\rrrr)*\dddd-(cos(\x r))^2)/
    (2*\tttt*(1-\tttt + \tttt*cos(\rrrr r)*\dddd ) ) } ); 

\draw[red!80!blue, ultra thick, domain=\PITT:6] 
  plot (\x, {\aaaa+(1-\tttt +\tttt*cos(\rrrr)*\dddd)/
    (2*\tttt*(1-\tttt + \tttt*cos(\rrrr r)*\dddd ) ) } )
   node[below]{$\xi_t$\qquad\mbox{}};
\def\tttt{0.55}
\def\aaaa{-4.0}
\def\aaab{-3.7}

\fill[color=gray!10] (-3,\aaaa-1.1) rectangle (6,\aaaa+2.0);
\draw[->] (-3,\aaaa)-- node[pos=0.8, above] {$t=\tttt$} (5.5,\aaaa) node[right]{$x$};
\fill[color=red!80!blue] (0,\aaaa) circle (3pt) node[below left] {$z_0$};
\fill[color=blue!80!red] (\rrrr,\aaaa) circle (3pt) node[below] {$z_1$};

\draw[blue!80!red, ultra thick, domain=-3:\rrrr-\PITT] 
  plot (\x, {\aaab+(0.0 -\tttt -(1-\tttt)*cos(\rrrr r)/\dddd)/
    (2*(1-\tttt)*(\tttt + (1-\tttt)*cos(\rrrr r)/\dddd ) ) } );

\draw[blue!80!red, ultra thick, domain=\rrrr-\PITT:\rrrr+\PITT] 
  plot (\x, {\aaab+((cos((\x-\rrrr) r))^2 -\tttt -(1-\tttt)*cos(\rrrr r)/\dddd)/
    (2*(1-\tttt)*(\tttt + (1-\tttt)*cos(\rrrr r)/\dddd ) ) } );

\draw[blue!80!red, ultra thick, domain=\rrrr+\PITT:6] 
  plot (\x, {\aaab+(0.0 -\tttt -(1-\tttt)*cos(\rrrr r)/\dddd)/
    (2*(1-\tttt)*(\tttt + (1-\tttt)*cos(\rrrr r)/\dddd ) ) } ) 
    node[above] {$\bar\xi_t$\qquad\mbox{}}   ;

\draw[red!80!blue, ultra thick, domain=-3:-\PITT] 
  plot (\x, {\aaaa+(1.0 -\tttt +\tttt*cos(\rrrr)*\dddd)/
    (2*\tttt*(1-\tttt + \tttt*cos(\rrrr r)*\dddd ) ) } ); 

\draw[red!80!blue, ultra thick, domain=-\PITT:\PITT] 
  plot (\x, {\aaaa+(1-\tttt +\tttt*cos(\rrrr)*\dddd-(cos(\x r))^2)/
    (2*\tttt*(1-\tttt + \tttt*cos(\rrrr r)*\dddd ) ) } );  
\draw[red!80!blue, ultra thick, domain=\PITT:6] 
  plot (\x, {\aaaa+(1-\tttt +\tttt*cos(\rrrr)*\dddd)/
    (2*\tttt*(1-\tttt + \tttt*cos(\rrrr r)*\dddd ) ) } )
   node[below]{$\xi_t$\qquad\mbox{}};
\end{tikzpicture}}
\hfill
\begin{minipage}[b]{0.37\textwidth}
  \caption{For the case $\varrho=|z_1{-}z_0|=0.9<\pit$ the functions $\xi_t(x)$
    (red) and $\bar\xi_t(x)$ (blue) from 
    \eqref{eq:xi.barxi} are displayed for the different times $t=0.15,\ 0.35$,
    and $0.55$ (with parameter $r_1/r_0=2$). We
    always have $\xi_t(x) \geq \bar\xi_t(x) $ with equality at the one-point
    contact set $\Xi_t=\{z(t)\}$, where $z(t) =\bfT_{0\to t}(z_0)$ moves
    continuously from $z_0$ to $z_1$.}
\end{minipage}
\label{fig:ContactSets} 
\end{figure}

The general theory in Theorem~\ref{thm:HJ1}(i) guarantees $\xi_t\geq
\bar\xi_t$. A lengthy computation shows that $\Xi_t$ is a singleton also for
$t\in (0,1)$, i.e.\ $\Xi_t=\{a(t)\} $ from $\mu_t=r(t)^2 \delta_{z(t)}$ and
$\Xi^\pm=\emptyset$.  We refer
to Figure \ref{fig:ContactSets}, where $x\mapsto (\xi_t(x),\bar\xi_t(t))$ is
plotted.\medskip 

\begin{figure}
{\begin{tikzpicture}
\def\PITT{1.5708}
\def\tttt{0.35}
\def\aaaa{-0.5}
\def\aaab{-0.55}
\def\rrrr{1.5708} 

\fill[color=gray!10] (-3,\aaaa-1.5) rectangle (6,\aaaa+2.2);

\node [ fill=white] at (0.0,\aaaa+1.7) {$|z_1{-}z_0|=\pit$};

\draw[->] (-3,\aaaa)--node[pos=0.8, above] {$t=\tttt$} (5.5,\aaaa) node[right]{$x$};
\fill[color=red!80!blue] (0,\aaaa) circle (3pt) node[below] {$z_0$};
\fill[color=blue!80!red] (\rrrr,\aaaa) circle (3pt) node[below] {$z_1$};

\draw[blue!80!red, ultra thick, domain=-3:\rrrr-\PITT] 
  plot (\x, {\aaab+(0.0 -\tttt)/(2*(1-\tttt)*\tttt  ) } );

\draw[blue!80!red, ultra thick, domain=\rrrr-\PITT:\rrrr+\PITT] 
  plot (\x, {\aaab+((cos((\x-\rrrr) r))^2 -\tttt )/
    (2*(1-\tttt)*\tttt ) } );

\draw[blue!80!red, ultra thick, domain=\rrrr+\PITT:6] 
  plot (\x, {\aaab+(0.0 -\tttt )/(2*(1-\tttt)*\tttt  ) } ) 
    node[below] {$\bar\xi_t$\qquad\mbox{}}   ;

\draw[red!80!blue, ultra thick, domain=-3:-\PITT] 
  plot (\x, {\aaaa+(1.0 -\tttt)/(2*\tttt*(1-\tttt) ) } );  

\draw[red!80!blue, ultra thick, domain=-\PITT:\PITT] 
  plot (\x, {\aaaa+(1-\tttt -(cos(\x r))^2)/
    (2*\tttt*(1-\tttt ) ) } ); 

\draw[red!80!blue, ultra thick, domain=\PITT:6] 
  plot (\x, {\aaaa+(1-\tttt)/
    (2*\tttt*(1-\tttt ) ) } )
   node[below]{$\xi_t$\qquad\mbox{}};
\def\rrrr{2.5}
\def\aaaa{-4.5}
\def\aaab{-4.52}

\fill[color=gray!10] (-3,\aaaa-1.5) rectangle (6,\aaaa+2.2);
\node [ fill=white] at (0.0,\aaaa+1.7) {$|z_1{-}z_0|=2.5$};

\draw[->] (-3,\aaaa)--node[pos=0.87, above] {$t=\tttt$} (5.5,\aaaa) node[right]{$x$};
\fill[color=red!80!blue] (0,\aaaa) circle (3pt) node[below] {$z_0$};
\fill[color=blue!80!red] (\rrrr,\aaaa) circle (3pt) node[below] {$z_1$};

\draw[blue!80!red, ultra thick, domain=-3:\rrrr-\PITT] 
  plot (\x, {\aaab+(0.0 -\tttt)/
    (2*(1-\tttt)*\tttt  ) } );

\draw[blue!80!red, ultra thick, domain=\rrrr-\PITT:\rrrr+\PITT] 
  plot (\x, {\aaab+((cos((\x-\rrrr) r))^2 -\tttt )/
    (2*(1-\tttt)*\tttt ) } );

\draw[blue!80!red, ultra thick, domain=\rrrr+\PITT:6] 
  plot (\x, {\aaab+(0.0 -\tttt )/
    (2*(1-\tttt)*\tttt  ) } ) 
    node[below] {$\bar\xi_t$\qquad\mbox{}}   ;

\draw[red!80!blue, ultra thick, domain=-3:-\PITT] 
  plot (\x, {\aaaa+(1.0 -\tttt)/
    (2*\tttt*(1-\tttt) ) } );  

\draw[red!80!blue, ultra thick, domain=-\PITT:\PITT] 
  plot (\x, {\aaaa+(1-\tttt -(cos(\x r))^2)/
    (2*\tttt*(1-\tttt ) ) } ); 

\draw[red!80!blue, ultra thick, domain=\PITT:6] 
  plot (\x, {\aaaa+(1-\tttt)/
    (2*\tttt*(1-\tttt ) ) } )
   node[below]{$\xi_t$\qquad\mbox{}};
\end{tikzpicture}}
\hfill
\begin{minipage}[b]{0.37\textwidth}
  \caption{For $\varrho=|z_1{-}z_0|\geq \pit$ the contact set $\Xi_t$ for the
    functions $\xi_t(x)$ (red) and $\bar\xi_t(x)$ (blue) from
    \eqref{eq:xi.rho.ll} is no longer a singleton. For $\varrho=\pit$ (upper
    figure) we obtain $\Xi_t=[z_0,z_1]$. For $\varrho>\pit$ (lower figure), we
    have $\Xi_t=\Xi^+_t\cup\Xi^-$ with $\Xi^-_t=\{z_0\}$ and $\Xi^+_t=\{z_1\}$. }
\end{minipage}
\label{fig:ContactSets2} 
\end{figure}
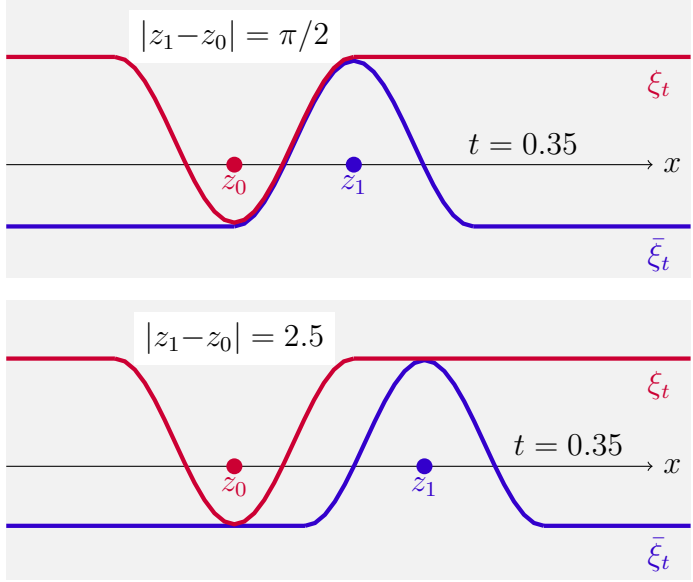

\underline{Case $\varrho \geq \pit$:} Now we have $s_0=s_1=0$ and $\xi_t$ and
$\bar\xi_t$ simplify accordingly:
\begin{equation}
  \label{eq:xi.rho.ll}
  \xi_t(x) = \frac{1-t - \tcosq\pit{|x{-}z_0|}}{ 2\,t\,(1{-}t)} 
\quad \text{and} \quad 
\bar\xi_t(x) = \frac{\tcosq\pit{|x{-}z_1|}-t }{ 2\,t\,(1{-}t)}.  
\end{equation}
The contact sets are easily found depending on $\varrho=\pit$ or
$\varrho>\pit$, namely 
\begin{align*}
\varrho>\pit:\quad & \Xi_t= \Xi^-\cup \Xi^+_t \ 
  \text{ with } \Xi^-=\{z_0\} \text{ and } \Xi^+_t=\{z_1\},
\\ 
\varrho=\pit: \quad & \Xi_t=[z_0,z_1] 
 \text{ and } \Xi^-=\{z_0\} \text{ and } \Xi^+=\{z_1\},
\end{align*}
where $[z_0,z_1]$ denotes the segment $\bigset{ (1{-}\theta)z_0{+}\theta z_1 
}{ \theta\in [0,1]}$, see Figure \ref{fig:ContactSets2}. 

The interesting fact that for $\varrho=|z_1{-}z_0|=\pit$ the contact set
$\Xi_t$ is constant and consists of a full segment reflects the observation in
\cite[Sec.\,5.2]{LiMiSa16OTCR} that $\mu_0$ and $\mu_1$ can be connected by 
geodesics satisfying $\mafo{sppt}(\mu_t)=[z_0,z_1]$ for all $t\in [0,1]$. 
\end{example}

\subsection{Geodesic flow and characteristics}
\label{subsec:flow}

Finally, we study the differentiability of $\GRAD_s=\nabla\xi_s$ and
$\bfT_{t\to s}$ on $\Xi_s$.  Let us denote by $\wt\Xi_t$ the subset of density
points of the contact set $\Xi_t$, which is closed by \eqref{eq:116}:
\begin{equation}
  \label{eq:143}
  x\in \wt\Xi_t\quad\Leftrightarrow\quad
  \lim_{\varrho\downarrow0}\frac{\mathcal L^d(\Xi_t\cap
    B_\varrho(x))}{\mathcal L^d(B_\varrho(x))}=1.
\end{equation}
\GGG Notice that $\wt\Xi_t$ is just the set
of Lebesgue points of the characteristic functions of $\Xi_t$,
so that \cite{AmFuPa00FBVF} $\mathcal L^d(\Xi_t\setminus \wt\Xi_t)=0$.  
\EEE 
By
\cite[Thm.\,1]{Bucz92DPBL}, the family of sets $(\wt\Xi_t)_{t\in (0,1)}$
is invariant with respect to the action of the bi-Lipschitz maps
$\bfT_{s\to t}$, i.e., $\bfT_{s\to t}(\wt\Xi_s)=\wt\Xi_t$ for every
$s,t\in (0,1)$.

Given a locally Lipschitz function $F:\Xi_t\to \R^d$ and $x\in \wt\Xi_t$,
we say that $F$ is differentiable at $x$ if there exists a matrix
$\mathsf A=\rmD F(x)\in \R^{d\times d}$ such that 
\begin{equation}
  \label{eq:138}
  |F(y)-F(x)-\mathsf A(y-x)|=o(|y-x|)\quad\text{as }y\to x,\ y\in \Xi_t.
\end{equation}
Since $x$ \GGG belongs to the set $\wt\Xi_t$ of density points of $\Xi_t$, \EEE 
the matrix $\sfA$ is unique and every
(locally) Lipschitz extension of $F$ is differentiable at $x$ with the same
differential $\sfA$ (e.g.\ one can argue as in the proof of
\cite[Thm.\,2.14]{AmFuPa00FBVF}).

We call $\mathrm{dom}_t (\rmD F)$ the set of differentiability points
$x\in \wt\Xi_t$ of $F$.  If $F$ is locally Lipschitz in $\Xi_t$, considering an
arbitrary Lipschitz extension of $F$ and applying Rademacher's theorem, we know
that $\mathcal L^d(\Xi_t\setminus \mathrm{dom}_t(\rmD F))=0$.  We will use the
simple chain-rule property that if $y=F(x)$ is a density point of
$F(\Xi_t)$ and $H:F(\Xi_t)\to \R^k$ is differentiable at $y$, then
\begin{equation}
  \label{eq:ChainRule}
  \rmD(H\circ F)(x)=\rmD H(F(x))\cdot \rmD F(x).
\end{equation}

In the proof of the following lemma we will denote by $\partial\xi_s$
the Fr\'echet subdifferential of $\xi_s$, which coincides with $\nabla\xi_s$
whenever $\xi_s$ is differentiable, in particular in $x \in \Xi_s$.

\begin{lemma}
  \label{le:nightmare}
  Let $s\in (0,1)$ and let $x\in \wt\Xi_s$ be a density point of $\Xi_s$ 
  where
  $\GRAD_s=\nabla\xi_s$ is differentiable in the sense of \eqref{eq:138}
  with $p=\GRAD_s(x)$ and $ \sfA=\rmD\nabla\xi_s(x)$. Then
\begin{subequations}
  \begin{gather}
    \label{eq:77}
    \sfA=\rmD\nabla\xi_s(x)\quad\text{is symmetric,}\\
    \label{eq:77.1}
    \sup_{z\in \partial\xi_s(y)} |z-p-\sfA(y{-}x)| = 
      o\,(|y{-}x|) \quad \text{as }y\to x,\\
    \label{eq:77.2}
      \xi_s(y)-\xi_s(x)-\langle p, y{-}x\rangle-\frac 12 \langle\sfA (y{-}x),
      y{-}x\rangle = o\,(|y{-}x|^2) \quad \text{as }y\to x,
    \end{gather}
\end{subequations}
Analogous results hold for $\bar \xi_s$.  We will denote $\rmD \nabla \xi_s$ by
$\rmD^2 \xi_s$.
\end{lemma}
Notice that the points $y$ in the limits in \eqref{eq:77.1} and
\eqref{eq:77.2} are not restricted to $\Xi_s$.\medskip

\noindent
\begin{proof}
We adapt some ideas of \cite{BiCoPu96SDCS, AlbAmb99GAMF} to our setting, and
we consider the case of $\bar \xi_s$ (to deal with a semi-convex function,
instead of semi-concave).  We will assume $x=0$ and will shortly write
$\bar \xi$ and $\Xi$ for $\bar\xi_s$ and $\Xi_s$ omitting the explicit
dependence on the parameter $s$.  For $h>0$ we define the blowup set
$\Xi^{h}:=h^{-1}\Xi$.  Up to an addition of a quadratic term, it is also not
restrictive to assume that $\bar\xi$ is convex.
   
For $h>0$ we set
$\omega_h(y):=\frac 1{h^2}\big(\bar \xi(hy)-\bar \xi(x)-h\langle
p,y\rangle\big)$ so that $\omega_h$ is a convex and nonnegative function. By
\eqref{eq:54} there exists a positive constant $C$ such that
\begin{equation}
  \label{eq:79}
  0\le \omega_h(y)\le C|y|^2\quad  \text{for every }y\in \Xi^h.
\end{equation}
Since $x=0$ is a density point of $\Xi$, $\Leb d(B_r(0)\setminus \Xi^h)\to 0$
as $h\downarrow0$ so that every point of $z\in B_r(0)$ is a limit of a sequence
in $z_h\in \Xi^h\cap B_r(0)$.  Therefore, for $h$ sufficiently small we can
find points $y_{h,i}\in \Xi^h\cap B_{4d}(0)$, $i=1,\cdots, 2d$, such that
$\overline {B_2(0)}\subset \mathrm{conv}(\set{y_{h,i}}{ i=1,\cdots, 2d})$. For
this it is sufficient to approximate the (rescaled) elements of the canonical
basis $\pm \mathsf e_i$, $i=1,\cdots , d$.  If $y\in B_2(0)$ we then find
coefficients $\alpha_{h,i}\ge 0$, $\sum_i \alpha_{h,i}=1$ such that
\[
  \omega_h(y)\leq  \sum_i \alpha_{h,i}\omega_h(y_{h,i})
  \le C \sum_i \alpha_{h,i}|y_{h,i}|^2\le 2dC
\]
so that $\omega_h$ is uniformly bounded in $B_2(0)$ and therefore is also
uniformly Lipschitz in $\overline {B_1(0)}$.  Every infinitesimal sequence
$h_n\downarrow0$ has a subsequence $m\mapsto h_{n(m)}$ such that
$\omega_{h_{n(m)}}$ is uniformly convergent to a nonnegative, convex Lipschitz
function $\omega:\overline {B_1(0)}\to \R$. We want to show that any limit
point $\omega$ coincides with the quadratic function induced by the
differential $\sfA$, namely 
$\omega(y)=\omega_{\sfA}(y)=\frac 12\langle \sfA y,y\rangle$

Let $\omega$ be the uniform limit of $\omega_h$ along a subsequence
$h_n\downarrow0$.  If $y_n\in \Xi^{h_n}\cap B_1(0)$ is converging to
$y\in B_1(0)$ we know that any limit point of $p_n=\nabla \omega_{h_n}(y_n)$
belongs to $\partial \omega(y)$.  On the other hand,
$p_n=\frac1{h_n}(\nabla\bar \xi(h_ny_n)-p)=\sfA y_n+o(1)$ thanks to the
differentiability assumption, so that $\sfA y\in \partial \omega(y)$.  Since we
can approximate every point of $B_1(0)$ we conclude that
$\sfA y\in \partial\omega(y)$ for every $y\in B_1(0)$.  On the other hand,
$\omega$ is Lipschitz, so that it is differentiable a.e.\,in $B_1(0)$ with
$\nabla\omega(y)=\sfA y$ and therefore the distributional differential of
$\nabla \omega$ coincides with $\sfA$. We conclude that $\sfA$ is symmetric and
$\omega(y)=\frac 12\langle \sfA y,y\rangle.$ The fact that $\omega_h$ uniformly
converges to $\omega$ eventually yields \eqref{eq:77.1} and \eqref{eq:77.2}.
\end{proof}

We now use the second-order differentiability of $\xi_s$ to derive
differentiability of $\bfT_{s\to t}$ by using the formula \eqref{eq:136} with
$\GRAD_s(x)=\nabla \xi_s(x)$. For $s \in (0,1)$ we define 
\begin{equation}
  \label{eq:def.Ds}
  \DDD_s = \mathrm{dom}_s(\rmD\nabla\xi_s)) \cap \wt\Xi_s =
\mathrm{dom}_s(\rmD^2\xi_s) \cap \wt\Xi_s. 
\end{equation}
As we already observed, since $\GRAD_s$ is Lipschitz on $\Xi_s$,
$\Leb d(\Xi_s\setminus \DDD_s)=0$ for every
$s\in (0,1)$.

For $t\in (0,1)$ and $\tau=t{-}s$ we also have $1{+}2\tau \xi_s\geq
(1{-}t)/(1{-}s)>0$ so that 
\[
x \mapsto \frac{\tau}{1{+}2\tau \xi_s(x)}\,\nabla \xi_s(x) = \nabla \phi_{s,t}(x) \text{
  with } \phi_{s,t}(x)=\GGG \frac 12 \EEE \log\big(1{+}2\tau \xi_s(x) \big) .
\]
is again Lipschitz on $\Xi_s$. Thus, Lemma \ref{le:nightmare} can be applied
and $\phi_{s,t}$ is differentiable in the sense of \eqref{eq:138} on
$\DDD_s$. Finally, we exploit the explicit representation of $\bfT_{s\to t}$ via
\eqref{eq:136}, namely for all $x\in \Xi_s$ we have  
\begin{equation}
  \label{eq:bfTst}
  \bfT_{s\to t} (x) = x + \bfarctan\Big( \frac{\tau\nabla \xi_s(x)}{1{+}2\tau
  \xi_s(x)}\Big) = x + \bfarctan\big( 
  \nabla \phi_{s,t}(x)\big).
\end{equation} 
Now the chain rule \eqref{eq:ChainRule} guarantees the differentiability of
$\bfT_{s\to t}$ on the set $\DDD_s$:

\begin{lemma}[Differentiability of $\bfT$]
\label{le:Diff.bfT}
For all $s,t\in (0,1)$ the mapping $\bfT_{s\to t}$ is differentiable on $\DDD_s$,
and we have
\begin{subequations}
 \label{eq:Diff.bfT}
  \begin{align}
   \label{eq:163}
   &\begin{aligned} 
     &\rmD\bfT_{s\to t}(x)= \bbT\big(t{-}s,\xi_s(x),\nabla
      \xi_s(x),\rmD^2\xi_s(x) \big) \quad \text{with}\\
     &\bbT(\tau,\xi,\GRAD,\sfA):=\bbI+\big(\rmD^2
     \Ell_1(z)\big)^{-1}_{\big|{z=\bfarctan(\tfrac{\tau\GRAD}{1{+}2\tau\xi})}}
     \Big( \frac{\tau\, \sfA }{1{+}2\tau\xi}  - \frac{2\tau^2\,
     \GRAD\oti  \GRAD }{(1{+}2\tau\xi)^2} \Big).
  \end{aligned}
\\[0.3em]
  \label{eq:139}
  &\bfT_{s\to t}(\DDD_s)=\DDD_t \text{ \ and \ }
  \rmD\bfT_{t\to s}(\bfT_{s\to t}(x))\rmD \bfT_{s\to  t}(x)=\bbI \text{ for
  }x\in \DDD_s.
\intertext{For every $t_0,t_1\in (0,1)$, $t_2\in [0,1]$ we also have}
  \label{eq:140}
  &\rmD\bfT_{t_1\to t_2}(\bfT_{t_0\to t_1}(x))\rmD \bfT_{t_0\to t_1}(x)=
  \rmD \bfT_{t_0\to  t_2}(x)\quad  \text{for }x\in \DDD_{t_0}.
  \end{align}
\end{subequations}  
\end{lemma}
\begin{proof} 
Recall $\tau= t{-}s$, then the explicit formula \eqref{eq:163} follows
from differentiating  $\nabla \Ell_1\big(x{-}\bfT_{s\to  t} (x) \big)= -
\frac{\tau}{1+2\tau \xi_s} \nabla \xi_s$. 
Since $\bfT_{s\to t}^{-1}=\bfT_{t\to s}$ there exists a constant $L$ such that
\begin{equation}
  \label{eq:144}
  L^{-1}|x{-}x'|\le |\bfT_{s\to t}(x)-\bfT_{s\to t}(x')|\le L|x{-}x'|
  \quad\text{for every }x,x'\in \Xi_s.
\end{equation}
If $x\in \DDD_s$ and $A=\rmD\bfT_{s\to t}(x)$, choosing $\eps>0$ we can find
$\varrho>0$ such that
\begin{equation}
  \label{eq:145}
  |\bfT_{s\to t}(x')-\bfT_{s\to t}(x)-A(x'{-}x)|\le \eps|x'{-}x|\quad
  \text{for every }x'\in \Xi_t\cap B_\varrho(x),
\end{equation}
so that choosing $\eps<\frac 1{2L}$ and $x'=x+v$ we get
\begin{align*}
  |Av|\ge |\bfT_{s\to t}(x{+}v)-\bfT_{s\to t}(x)|- \eps|v|
  \ge \frac 1{2L}|v|\quad
  \text{for every }v\in B_\varrho(0)\cap (\Xi_t-x)
\end{align*}
Using the fact that $0$ is a density point of $\Xi_t-x$ we conclude that $A$ is
invertible with $|A^{-1}|\le 2L$.  For every $y'\in \Xi_t$ with
$L|y'{-}y|<\varrho$ and $x'=\bfT_{t\to s}(y')$, we get $|x'{-}x|<\varrho$ and
\eqref{eq:145} yields
\begin{align*}
  |\bfT_{t\to s}(y')-\bfT_{t\to s}(y)-A^{-1}(y'{-}y)|
  &=\big|A^{-1} \big(A(x'{-}x)-\bfT_{s\to t}(x')+\bfT_{s\to t}(x)\big)\big|
  \\&
  \le {2L} \eps |x'{-}x|\le 2L^2\eps |y'{-}y|    
\end{align*}
showing that $y\in\DDD_t$ and $A^{-1}=\rmD\bfT_{t\to s}(y)$. Hence, \eqref{eq:139}
is established.  

Equation \eqref{eq:140} then follows by the concatenation property
\eqref{eq:26}.
\end{proof}

The explicit formula \eqref{eq:163} shows that $\rmD \bfT_{s\to t}$ is the product of the positive matrix $\rmD^2\Ell_1(z)^{-1}$ and a symmetric
matrix, hence it is always real diagonalizable. The following result shows that
the determinant and hence all eigenvalues stay positive for $s,t\in (0,1)$. 
In fact, we now derive differential equations with respect to $t\in (0,1)$ for
the \TRAGRO\ pairs $(\bfT_{s\to t}(x), q_{s\to t}(x)) \in \R^d\ti
(0,+\infty)$ as well as for $\rmD \bfT_{s\to t}(x) \in \R^{d\ti d}$ and $\det
\rmD \bfT_{s\to t}(x) $. Recall that $t\mapsto (\bfT_{s\to t}(x), q_{s\to
  t}(x))$ is analytic for $t\in (0,1)$ by Theorem~\ref{thm:HJ1}(\ref{th:HJ.label5})
and (\ref{th:HJ.label6}).  

The following relations will be crucial to derive the curvature estimate
needed for our main result on geodesic $\HK$-convexity.

\begin{theorem}[The characteristic system on the contact set $\Xi_s$]
\label{thm:characteristics}
We fix $s\in (0,1)$, $x\in \Xi_s$, and $y\in \DDD_s$ (cf.~\eqref{eq:def.Ds}) and define the maps
\[
 \bfT(t):=\bfT_{s\to t}(x), \quad q(t):=q_{s\to t}(x), \quad
 \sfB(t) := \rmD \bfT_{s\to t}(y), \ \ \text{ and }\ \delta(t):= \det \sfB(t).
\]
Then, we have the initial conditions $\bfT(s)=x$, $q(s)=1$, $\sfB(s)=\bbI$, and
$\delta(s)=1$, and for $t\in (0,1)$ the following differential equations are
satisfied: 
\begin{subequations}
 \label{eq:CharSystem}
 \begin{align}
  \label{eq:157}
  & \dot \bfT(t)=\nabla\xi_t(\bfT(t))  &&\text{and } \ 
     \ddot\bfT(t)= -4\xi_t(\bfT(t))\,\nabla\xi_t(\bfT(t)), 
  \\ 
  \label{eq:160}
  &\dot q(t)=2\xi_t(\bfT(t))\,q(t)  &&\text{and }\ 
    \ddot q(t)=  |\nabla\xi_t(\bfT(t))|^2\, q(t), 
  \\
  \label{eq:158}
   & \dot\sfB(t) =\rmD^2\xi_t(\bfT(t)) \, \sfB(t) 
                                       &&\text{and }\ 
     \ddot\sfB(t) = -4\Big(\nabla\xi_t\oti\nabla\xi_t+ 
  \xi_t\rmD^2\xi_t\Big) {\circ} \bfT(t)\cdot \sfB(t),
  \\
  \label{eq:161}
  &\dot \delta(t)= \Delta \xi_t(\bfT(t)) \, \delta(t),&
  & \frac{\ddot \delta(t)}{\delta(t)}=\Big((\Delta \xi_t)^2 {-} 
  |\rmD^2\xi_t|^2 {-} 4|\nabla\xi_t|^2 {-} 4\xi_t\Delta
    \xi_t\Big) {\circ} \bfT(t),
  \end{align}
\end{subequations}
where $\Delta \xi_t(z)= \mathrm{tr}\big(\rmD^2 \xi_t(z)\big)$ and
$|\rmD^2\xi_t(z)|^2 = \sum_{i,j} \big(\pl_{x_i}\pl_{x_j}\xi_t(z)\big)^2$.  
\end{theorem}
\begin{proof}
We use \eqref{eq:26} and the Taylor expansion
\[
  \bfarctan\Big(\frac{ h \GRAD }{ 1{+}2h \xi }   \Big)\ 
  =  h   \GRAD  -  2h^2  \xi  \GRAD  + O(h^3)\quad\text{as }h\to 0.
\]
Setting $y=\bfT(t)=\bfT_{s\to t}(x)$ and using the fact that $y\in \Xi_t$,
\eqref{eq:bfTst} yields 
\[
  \bfT_{t\to t+h}(y)=y+h \nabla\xi_t(y)-2h^2 \xi_t(y)\nabla\xi_t(y)+
  O(|h|^3)\ \text{ as } h \to 0 .
\]
With the composition rule \eqref{eq:26} we have $\bfT_{s\to t+h} (x)=\bfT_{t\to
  t+h}(y)$ and compute  
\begin{align*}
  \dot \bfT(t)
  &= \lim_{h\to0}\frac{\bfT_{s\to {t+h}}(x)-\bfT_{s\to t}(x)}h
  = \lim_{h\to0}\frac{\bfT_{t\to {t+h}}(y)-y}h = \nabla\xi_t(y).
\end{align*}
This identity yields the first equation in \eqref{eq:157}. For the second relation in
\eqref{eq:157} we use
\begin{align*}
  \ddot\bfT(t)
  &=  \lim_{h\to0}\frac{\bfT_{s\to {t+h}}(x)-2\bfT_{s\to  t}(x)+\bfT_{s\to t-h}(x)}{h^2}
  \\& =
  \lim_{h\to0}\frac{\bfT_{t\to {t+h}}(y)-2y+\bfT_{t\to
      {t-h}}(y)}{h^2}  =-4\xi_t(y)\nabla\xi_t(y).
\end{align*}
  
The relations \eqref{eq:160} for $q(t)=q_{s\to t}$ follow similarly, using
the scalar product rule for $q_{s\to t}$ in \eqref{eq:26} and by
taking the square root of \eqref{eq:137}, namely 
\[
  q_{t\to t+h}(y)=1+2h \xi_t(y)+ \frac{h^2} 2
  |\nabla\xi_t(y)|^2  + o(h^2)\quad\text{as }h\to0.
\]
 
To show that $\sfB(t)$ satisfies \eqref{eq:158}, we exploit the matrix
product rule  \eqref{eq:140} and expand $\rmD \bfT_{t\to t+h}(y)$ in
\eqref{eq:163} to obtain 
\begin{equation}
  \rmD \bfT_{t\to t+h}(y) = \bbI + h \rmD^2 \xi_t - 2h^2 \Big(
  \nabla\xi_t{\otimes}\nabla \xi_t + \xi_t \rmD^2 \xi_t\Big) +
  o(h^2)   \quad\text{as }h\to0.
  \label{eq:159}
\end{equation}
For this note that $y-\bfT_{t\to t+h}(y)=O(|h|)$ so that
$\rmD^2\Ell_1\big(y{-}\bfT_{t\to t+h}(y)\big)= \bbI + O(|h|^2)$ as $\Ell_1$ is
even. Thus, \eqref{eq:158} follows as in the previous two cases. 

For the determinant $\delta(t)$ we again have a scalar product rule, and it
suffices to expand  $\det(\rmD \bfT_{t\to t+h}(x)) $
at $h=0$. For this we can use the classical expansion 
$  \det (\bbI{+}h\bbA)= 1 +
  h\trace \bbA +\frac12 h^2  \big(  (\trace \bbA)^2 - \trace(\bbA^2)\big) +
  O(h^3)$, and obtain
\begin{equation}
  \label{eq:162}
  \det \rmD \bfT_{t\to t+h}=
  1+h \Delta \xi_t+
  \frac12 h^2  \Big(  
  (\Delta\xi_t)^2 - |\rmD^2\xi_t|^2 -4|\nabla\xi_t|^2-4\xi \Delta\xi_t\Big) +o(h^2).
\end{equation}
As before this shows \eqref{eq:161}, and the theorem is proved. 
\end{proof}

In this section, we have studied the forward solutions $t\mapsto
\xi_t$ for $t\in (0,1)$ and its contact sets $\Xi_t$ with a corresponding backward
solution $\bar \xi_t$. We
obtained differentiability properties
in these sets or in the slightly smaller sets $D_t$
and derived transport relations for
important quantities such as $q_{s\to t}$ and $\delta_s(t)=\det \rmD \bfT_{s\to
  t}(x)$. In the following section, we still have to show that the contact sets
$\Xi_t$ are sufficiently big, if we define $\xi_t=\HopfLax t \xi_0$ and
$\bar\xi_t = \RopfLax{1-t} \xi_1$ for an optimal pair $(\xi_0,\xi_1)$. This
will be done in Theorem~\ref{thm:HL}.

\section{Geodesic curves}
\label{se:GeodCurves}

In this section, we improve the characterization of
Hellinger--Kantorovich geodesic curves as discussed already in
\cite[Sec.\,8.6]{LiMiSa18OETP}. More precisely, we consider constant-speed geodesics 
$\mu:[0,1]\to\MMM(\OOmega)$ that satisfy
\[
\forall\, s,t \in [0,1]: \quad
\HK(\mu(s),\mu(t)) = |s{-}t| \, \HK(\mu_0,\mu_1).
\]
We first show the optimality of potentials $\xi_t$ and $\bar\xi_t$
obtained from the forward or backward Hamilton--Jacobi equation in Theorem~\ref{thm:HL}. With this, we are
able to show in Theorem~\ref{thm:HL2} that for subparts $(s,t) \subset [0,1]$
with $\tau=t{-}s<1$ the 
corresponding LET problem has a unique solution in Monge form, which implies
that $(\MMM(\R^d),\HK)$ has the strong non-branching property. Finally, in
Theorem \ref{th:restriction} and Corollary \ref{cor:start} we provide
restrictions and splittings of geodesic curves needed for the main theorem in
Section \ref{se:Geod.Conv}.


\subsection{Geodesics and Hamilton--Jacobi equation}
 
The next result clarifies the connection with the forward and backward Hopf--Lax
flows $\xi_t$ and $\bar\xi_t$ studied in Theorem~\ref{thm:HJ1} and
the importance of the contact set $\Xi_t$ defined in \eqref{eq:116} (see also
\cite[Thm.\,8.20]{LiMiSa18OETP} and \cite[Chap.\,7]{Vill09OTON} for a similar
result in the framework of Optimal Transport and displacement
interpolation). We emphasize that despite the non-uniqueness of the
geodesics $(\mu_t)_{t\in [0,1]}$ (see \cite[Sec.\,5.2]{LiMiSa16OTCR}) in the
following result, $\xi_t$ and $\bar\xi_t$ only depend on
$\mu_0$ and $\mu_1$ and the optimal potentials $\varphi_0$ and
$\varphi_1$. 

The  result brings together the results of Sections \ref{se:2nd.Optim} and
\ref{sec:HJ} by starting with an optimal pair $(\varphi_0,\varphi_1)$ from
Section \ref{se:2nd.Optim} and considering the corresponding solutions $\xi_t$
and $\bar\xi_t$ of the forward and backward Hamilton--Jacobi equation starting
with $\xi_0 = \chGtrafo_1(\varphi_0)$ and $\bar\xi_1 = \Gtrafo_1(\varphi_1)$,
respectively.  
First, we observe that ``intermediate'' pairs $(\xi_s,\xi_t)$ or
$(\bar\xi_s,\bar \xi_t)$ are optimal for connecting the intermediate points
$\mu_s$ and $\mu_t$ on an arbitrary geodesic connecting $\mu_0$ and
$\mu_1$. Second, we observe that certain results obtained in Section \ref{sec:HJ}
for $s,t\in (0,1)$ also hold in the limit points $s,t\in \{0,1\}$. Finally, we
show that the contact set $\Xi_t$ is large enough in the sense that it contains $\mafo{supp}(\mu_t)$ (see Example
\ref{ex:ContactDirac} for some instructive case with $\varrho=\pit$).

\begin{theorem}
\label{thm:HL}
For $\mu_0,\mu_1\in \MMM(\OOmega)$ consider a tight optimal pair
$(\varphi_0,\varphi_1)$ of (lower,upper) semi-continuous potentials as in
Theorem~\ref{thm:regularity}. With
$\xi_0:=\chGtrafo_1(\varphi_0)= \frac12(\ee^{2\varphi_0}-1)$ and
$\bar\xi_1:= \Gtrafo_1(\varphi_1)= \frac 12(1-\ee^{-2\varphi_1})$ we
define $\xi_t=\HopfLax t \xi_0$ and $\bar\xi_t = \RopfLax{1-t} \bar\xi_1$ as in
\eqref{eq:37} and the contact sets $\Xi_t=\{\xi_t=\bar\xi_t\}$ as in
\eqref{eq:116}. Finally, consider an arbitrary geodesic
$(\mu_t)_{t\in [0,1]}$ connecting $\mu_0$ to $\mu_1$. Then, the following
holds:
\begin{enumerate}[{\upshape(1)}]
\item For all $s,t \in [0,1]$ with $s<t$ both pairs
  $(\psxi_s,\psxi_t)$ and $(\bar\psxi_s,\bar\psxi_t)$ are optimal for
  \eqref{eq:8} and \eqref{eq:10} for connecting $\mu_s$ to $\mu_t$, viz.\ 
    \begin{equation}
      \label{eq:93}
      \frac1{2(t{-}s)}\HK^2(\mu_s,\mu_t)=
      \int \xi_t\dd\mu_t-\int\xi_s\dd\mu_s=
      \int \bar\xi_t\dd\mu_t-\int\bar\xi_s\dd\mu_s
    \end{equation}
       \item $S_t=\supp(\mu_t)\subset \Xi_t$ for every $t\in [0,1]$.
  \end{enumerate}
\end{theorem}
\begin{proof}
  \noindent
\underline{Assertion (1).} It is sufficient to consider the forward flow $\HopfLax t{}$. Fixing

$t\in (0,1)$ we have
\begin{equation}
  \label{eq:94}
  \frac1{2t}\HK^2(\mu_0,\mu_t) \geq \int\xi_t\dd\mu_t-\int \xi_0\dd\mu_0 
  \text{ and }
  \frac1{2(1{-}t)}\HK^2(\mu_t,\mu_1)\geq  \int\xi_1\dd\mu_1-\int \xi_t\dd\mu_t.
\end{equation}
On the other hand, the geodesic property and the optimality of $(\xi_0,\xi_1)$
yield
\begin{align*}
    \int\xi_1\dd\mu_1-\int \xi_0\dd\mu_0=\frac 12\HK^2(\mu_0,\mu_1)=
    \frac1{2t}\HK^2(\mu_0,\mu_t)+ \frac1{2(1{-}t)}\HK^2(\mu_t,\mu_1)
\end{align*}
showing that the inequalities in \eqref{eq:94} are in fact equalities, 
in particular \eqref{eq:93} with $s=0$. For $s>0$ we still get \eqref{eq:93}  since
$\frac 1{2(t-s)}\HK^2(\mu_s,\mu_t) =
\frac1{2t}\HK^2(\mu_0,\mu_t)-\frac1{2s}\HK^2(\mu_0,\mu_s)$ if $0<s<t\le 1$. 
\medskip

\noindent
\underline{Assertion (2).} Equation \eqref{eq:93} for $s=0$ yields
$ \int (\xi_t-\bar\xi_t)\dd\mu_t=0$ for all $ t\in (0,1)$, so that 
$\xi_t \leq \bar\xi_t$ and the continuity of $\xi_t,\bar\xi_t$ yield
$\xi_t=\bar \xi_t$ on $S_t = \supp \mu_t$.
The cases $t=0$ and $1$ follow by the relations between $\xi_i$ and $\varphi_i$ 
and the fact that $\varphi_0=\varphi_1^{\shortleftarrow L_1}$, 
$\varphi_1=\varphi_0^{\shortrightarrow L_1}$. 
\end{proof}

Note that the inclusion $S_t=\supp(\mu_t)\subset \Xi_t$ is in general a
strict inclusion. This can be seen for the case $|z_1{-} z_0| = \pit$ in Example
\ref{ex:ContactDirac}, where $\Xi_t =[z_0,z_1]$, however, there exists a pure Hellinger
geodesic with $\mafo{supp}(\mu_t)= \{z_0,z_1\}$ for $t\in (0,1)$. 
 
We can now exploit all the regularity features of the maps $\bfT_{s\to
  t}$ and $q_{s\to t}$ on the contact set $\Xi_t$ (cf.\ Theorem~\ref{thm:HJ1}). 
A first important consequence is that, given an $\HK$ geodesic $(\mu_t)_{t\in
  [0,1]}$ and $s\in (0,1)$, the $\HK$ problem 
between $\mu_s$ and $\mu_t$ for any $t \in [0,1]$ has only one
solution, which can be expressed in Monge form
(see \cite[Lem.\,7.2.1]{AmGiSa08GFMS} for the corresponding
properties for the $\rmL^2$-Wasserstein distance in $\R^d$).

\begin{theorem}[Regularizing effect along geodesics]
  \label{thm:HL2}
  Under the assumptions of Theorem~\ref{thm:HL}, if $s\in (0,1)$ and
  $t\in [0,1]$, then the \TRAGRO\ pair
  $(\bfT_{s\to t}, q_{s\to t})$ of Theorem~\ref{thm:HJ1} is the unique solution
  of the Monge formulation \eqref{eq:24} of the Entropy-Transport problem between $\mu_s$ and
  $\mu_t$.  In particular, the optimal Entropy-Transport problem between
  $\mu_s$ and $\mu_0$ or between $\mu_s$ and $\mu_1$ has a
  unique solution, \GGG and this solution is in Monge form. \EEE
\end{theorem}
\begin{proof}
  Let us consider the case $0<s<t\le 1$, $\tau=t-s<t$.  By Theorem~\ref{thm:HL},
  the pair $(\bar \xi_s,\bar \xi_t)$ is optimal for $(\mu_s,\mu_t)$ and
  $\supp(\mu_s)\subset \Xi_s$.  Using the transformations
  \begin{equation}
    \label{eq:106}
    \varphi_0:=\frac 1{2\tau}\log(1{+}2\tau\bar\xi_s)\quad \text{and} \quad 
    \varphi_\tau:=-\frac 1{2\tau}\log(1{-}2\tau\bar\xi_t),
  \end{equation}
  we see that $(\varphi_{0},\varphi_{\tau})$
  is a pair of potentials satisfying the
  assumptions of Theorem~\ref{thm:regularity}(2).
  Since $1-2\tau\xi_t\ge
  1-\tau/t>0$ we deduce that $\varphi_\tau$ is bounded from above,
  so that $\mu_t''=0$ thanks to \eqref{eq:105} (where the measures $\mu_t'$ and 
  $\mu_t''$ are defined as in \eqref{eq:43}).
  
  Moreover, we know that $\mu_s'$ is concentrated on $\{\varphi_0>-\infty\}$;
  since it is also concentrated on $\Xi_s$ we deduce that $\mu_s'$ is
  concentrated on $D_0' = \mathrm{dom}(\nabla\varphi_0)$, 
  so that we can apply Corollary \ref{cor:Monge},
  recalling the expression of $\bfT,q$ given by \eqref{eq:transport.map}.
\end{proof}

The above theorem allows us to deduce the fact that
$(\MMM(\R^d),\HK)$ has a strong non-branching property. It is shown in
\cite[Sec.\,5.2]{LiMiSa16OTCR} that the set of geodesics connecting two Dirac
measures $\delta_{y_0}$ and $\delta_{y_1}$ is very large if
$|y_1{-}y_0|=\pi/2$: it is convex but does not lie in a finite-dimensional
space. The following result shows that all these geodesics are mutually
disjoint except for the two endpoints $\mu_0$ and $\mu_1$. 

\begin{corollary}[Strong non-branching]
  \label{cor:nonbranching}
  If for some $s\in (0,1)$ we have $\HK(\mu_0,\mu_s)=s\HK(\mu_0,\mu_1)$ and
  $\HK(\mu_s,\mu_1)=(1{-}s)\HK(\mu_0,\mu_1)$, then
  there exists a unique geodesic curve $t\mapsto \mu(t)$
  such that $\mu(0)=\mu_0$, $\mu(s)=\mu_s$, and $\mu(1)=\mu_1$.
\end{corollary}

The next result shows that from a given geodesic we may construct new
geodesics by multiplying the measures $\mu_t$ by a suitably transported
function. This will be useful in the proof of the main Theorem~\ref{th:Geod.Cvx}.

\begin{theorem}[Restriction of geodesics]
\label{th:restriction}
Let $(\mu_t)_{t\in [0,1]}$ be an $\HK$ geodesic.  For a given $s\in (0,1)$ let
$\nu_s\in \MMM(\R^d)$ with $\supp(\nu_s)\subset \supp(\mu_s)$. Then the curve
$[0,1]\ni t\mapsto \nu_t:=(\bfT_{s\to t},q_{s\to t})_\star \nu_s$ is
also an $\HK$ geodesic.  If in addition $\nu_s=\varrho_s\mu_s$ for
some Borel function $\varrho_s:\supp(\mu_s)\to [0,+\infty]$, then
$\nu_t' = \varrho_t\mu_t$ with
$\varrho_t(y)= \varrho_s(\bfT_{t\to s} (y)) $ for every $t\in (0,1)$.
\end{theorem}
\begin{proof}  
We keep the same notation of Theorem~\ref{thm:HL},  let
$0<t_1<s<t_2<1$, and set $\tau_1:=s-t_1$, $\tau_2:=t_2-s$, and
$\tau=\tau_1{+}\tau_2$.  We clearly have
\begin{align*}
  \frac1{2\tau}\HK^2(\nu_{t_2},\nu_{t_1})
  &\ge \int\xi_{t_2}\dd\nu_{t_2}-
  \int\xi_{t_1}\dd\nu_{t_1}
  \\&=
  \Big(\int\xi_{t_2}\dd\nu_{t_2}-\int \xi_s\dd\nu_s\Big)+
  \Big(\int \xi_s\dd\nu_s-\int\xi_{t_1}\dd\nu_{t_1}\Big)
\end{align*}
The conclusion then follows, if we show 
$\int\xi_{t_2}\dd\nu_{t_2}-\int \xi_s\dd\nu_s\ge \frac
1{2\tau_2}\HK^2(\nu_{t_2},\nu_s)$ and
$\int\xi_{s}\dd\nu_{s}-\int \xi_{t_1}\dd\nu_{t_1}\ge \frac
1{2\tau_2}\HK^2(\nu_s,\nu_{t_1})$.  We check the first inequality, the
second follows similarly.

Define $q_2:=q_{s\to t_2}$ and $\bfT_2:=\bfT_{s\to t_2}$. Using the fact that
(i) $\nu_{t_2}=(\bfT_{2},q_{2})_\star \nu_s$ and (ii) identity \eqref{eq:137} we
obtain
\begin{align*}
      \int (1 {-} 2\tau_2 \xi_{t_2})\dd\nu_{t_2}
      \overset{\text{(i)}}=
      \int \Big(1{-}2\tau_2\xi_{t_2}\big(\bfT_{s\to t_2}(x)\big)\Big)q^2_{s\to
                     t_2} (x) \dd \nu_{s}(x) 
      \overset{\text{(ii)}}=
      \int (1{+}2\tau_2\xi_s) \dd \nu_s .
  \end{align*}
Combining \eqref{eq:136} and \eqref{eq:137}, we arrive at 
\begin{equation}
  \label{eq:166}
  \int\big(1-2\tau\xi_{t_2}\big) \dd\nu_{t_2}=
  \int \big(1+2\tau_2\xi_s\big) \dd\nu_s=
  \int q_2\cos(|x{-}\bfT_2(x)|)\dd\nu_s.
\end{equation}
With this we find
\begin{align*}
  \HK^2(\nu_{t_2},\nu_s)
  &\overset{\text{\eqref{eq:95}}}\leq
  \int \big( q_2^2+1-2q_2\cos(|x{-}\bfT_2(x)|\big) \dd\nu_s
  \\& \overset{\text{\eqref{eq:166}}}= \nu_{t_2}(\R^d) + 
  \nu_{s} (\R^d) - \int(1{-}2\tau\xi_{t_2})\dd\nu_{t_2}-    
  \int (1{+}2\tau_2\xi_s)\dd\nu_s
  \\&=2\tau\Big(\int\xi_{t_2}\dd\nu_{t_2}-\int \xi_s\dd\nu_s\Big).
\end{align*}
Hence, we have shown $\frac1{2\tau}\HK^2(\nu_{t_2},\nu_{t_1}) =
\int\xi_{t_2}\dd\nu_{t_2}-\int \xi_s\dd\nu_s $, which implies that
$(\nu_t)_{t\in (0,1)}$ is a geodesic as well.
   
We can then pass to the limits $t_1\downarrow0$ and $t_2\uparrow 1$ as
follows.  Notice that the curve $t\mapsto \nu_{t}$, $t\in (0,1)$, is
converging in $(\MMM(\R^d),\HK)$ to a limit $\nu_0$ and $\nu_1$ for
$t\downarrow0$ and $t\uparrow 1$, since $(\nu_t)$ is a geodesic. Moreover, for
every $\zeta\in \rmC_b(\R^d)$ we can pass to the limit $t\uparrow 1$ in
\begin{equation}
  \label{eq:165}
  \int \zeta\dd\nu_t=\int \zeta(\bfT_{s\to t}(x))q^2_{s\to t}(x)\dd\nu_s(x),
\end{equation}
since $\lim_{t\uparrow1}\bfT_{s\to t}(x)= \bfT_{s\to 1}(x)$
and $\lim_{t\uparrow1}q_{s\to t}(x)= q_{s\to 1}(x)$ and $q$ is uniformly
bounded.  A similar argument holds for the case $t\downarrow0$.
  
In order to check the identity concerning the density $\varrho_t'$ of $\nu_t$,
we use \eqref{eq:165} and find
\begin{align*}
  \int\zeta\dd\nu_t
  &=
  \int \zeta(\bfT_{s\to t}(x))q^2_{s\to t}(x)\dd\nu_s
  =\int \zeta(\bfT_{s\to t}(x))q^2_{s\to
    t}(x)\varrho_s(x)\dd\mu_s
  \\&=
  \int \zeta(\bfT_{s\to t}(x))q^2_{s\to
    t}(x)\varrho_t(\bfT_{s\to t}(x))\dd\mu_s(x)=
  \int \zeta(y)\varrho_t(y)\dd\mu_t(y).
\end{align*}
The case $t \in [0,s]$ is analogous. 
\end{proof}

The next result provides the fundamental formula for the representation of
densities along geodesics. Generalizing the celebrated formulas for the
Kantorovich--Wasserstein geodesics, the densities are again obtained by
transport along geodesics, but now with non-constant speed and an additional
\AAA growth \EEE factor $a_s(t,x) =q^2_{s\to t}(x)$ to account for the
annihilation and creation of mass. Recall that
$\DDD_s= \mafo{dom} (\rmD^2\xi_s) \subset \Xi_s$ has full Lebesgue measure in
$\Xi_s$, i.e.\ $\Leb d(\Xi_s\setminus \DDD_s)=0$.

\begin{corollary}[Representation of densities along geodesics]
\label{cor:start}
For $\mu_0,\mu_1\in \MMM(\R^d)$ consider a geodesic $(\mu_t)_{t\in [0,1]}$ 
connecting $\mu_0$ to $\mu_1$. Assume that at least one of the
following properties holds:
\vspace{-0.3em}
\begin{enumerate}[\upshape(a)]\itemsep-0.3em
\item there exists $s\in (0,1)$ such that $\mu_s=c_s\Leb d\ll\Leb d$;
\item $\mu_0=c_0\Leb d\ll\Leb d$ and $\mu_{1}''\ll \Leb d$.
\end{enumerate}
Then, we have
\begin{enumerate}[\upshape(1)]
\item $\mu_t\ll\mathcal L^d$ for every $t\in (0,1)$, viz.\ $\mu_t =
  c(t,\cdot) \Leb d$. 

\item For every $s\in (0,1)$ the density $c(t,\cdot)$ can be
  expressed via the formula 
  \begin{subequations}
  \label{eq:150all}
   \begin{align}
    \label{eq:150}
    c(t,y)\big|_{y= \bfT_{s\to t}(x)} =    c(s,x)\frac{\alpha_s(t,x)}{\delta_s(t,x)}
    \quad \text{for every }x\in \DDD_s,\quad t\in (0,1),
  \end{align}
  with $ \DDD_s = \mathrm{dom}_s(\rmD\nabla\xi_s))=
  \mathrm{dom}_s(\rmD^2\xi_s) $ (cf.\ \eqref{eq:def.Ds})  and 
  \begin{equation}
    \label{eq:151}
    \alpha_s(t,x):= \big(1+2(t{-}s)\psxi_s(x)\big)^2+(t{-}s)^2|\nabla\xi_s(x)|^2,
    \quad    \delta_s(t,x):=\det \rmD\bfT_{s\to t}(x).
  \end{equation}
  \end{subequations}
  Moreover, we have $\rmD^2\xi_s(x)=0$ and $\delta_s(t,x)=1$ for
  $\Leb d$-a.e.\,$x\in \GGG \Xi^0_s\supset \Xi^\pm$;
  in particular 
  \begin{equation}
    \label{eq:169}
    c(t,x)=\frac{t^2}{s^2}c(s,x)\quad\text{for  }x\in \Xi^+\quad 
    \text{and} \quad 
    c(t,x)=\frac{(1{-}t)^2}{(1{-}s)^2}c(s,x)\quad
    \text{for }x\in \Xi^-.
  \end{equation}

\item If $\mu_0\ll\Leb d$ (resp.\,$\mu_1\ll\Leb d$) \eqref{eq:150all} and
  \eqref{eq:169} hold up to $t=0$ (resp.\,up to $t=1$).

\item If $\mu_1''=0$ the representations in \eqref{eq:150all} also hold
      for $s=0$ by restricting $x$ in $D_0''=\mathrm{dom}(\rmD^2\varphi_0)$, and we have the formula
  \begin{equation} 
   \rmD \bfT_{0\to t}(x) =
    \bbT\big(t,\psxi_0(x),\nabla\psxi_0(x),\rmD^2\psxi_0(x)\big) \quad
    \text{for every }x\in D_0'',  
  \end{equation}
where $\bbT $ is defined in \eqref{eq:163}. 
\end{enumerate}
\end{corollary}
\begin{proof}
  \underline{Assertion (1).}  In the case (a) holds for $s\in (0,1)$, there
  exists a bi-Lipschitz map $\bfT_{s,t}:\Xi_s\to\Xi_t$ and bounded \AAA growth
  factors \EEE $q_{s,t}:\Xi_s\to [a,b]$ with $0<a<b<\infty$ such that
  $\mu_t=(\bfT_{s\to t},q_{s\to t})_\star \mu_s$. In particular, for every
  Borel set $A$ we have
\begin{equation}
  \label{eq:BorelSetA}
  \mu_t(A)\le b^2\mu_s(\bfT_{s\to t}^{-1}(A))  =b^2\mu_s(\bfT_{t\to s}(A)).
\end{equation}
If $\Leb d(A)=0$ then $\Leb d(\bfT_{t\to s}(A))=0$ because $\bfT_{t\to s}$ is
Lipschitz. Hence, using $\mu_s \ll \Leb d$ we find 
$\mu_s(\bfT_{t\to s}(A)=0$, such that \eqref{eq:BorelSetA} gives
$\mu_t(A)=0$. With this we conclude $\mu_s\ll\Leb d$.  

In the case of assumption (b), we argue as before but with $\mu_0=c_0\Leb
d$ for $s=0$. Using  the fact
that $q_{t\to 0}$ is locally bounded from below and that $\bfT_{t\to 0}$ is locally
Lipschitz on $A_t:=\Xi_t\setminus \Xi^+$,  we deduce that
$\mu_t\res A_t\ll \Leb d.$ On the other hand we have $\mu''_1\ll \Leb d$
and the restriction of $\bfT_{t\to 1}$
to $\Xi^+$ coincides with the identity and $q_{t\to1}$ is bounded from below
thanks to \eqref{eq:137}. Thus, we obtain $\mu_t\ll \Leb d$.

\medskip\noindent \underline{Assertion (2).}  The representation \eqref{eq:150}
follows by Theorem~\ref{thm:HL2} and Corollary \ref{cor:Monge}.

Relation \eqref{eq:169} can be deduced directly by Theorem~\ref{thm:HL2}.  In
order to prove that $\rmD^2\xi_s=0$ $\mu_s$-a.e.\,in $\Xi^\pm$ it is sufficient
to consider density points of $\Xi^\pm$, since $\mu_s\ll\Leb d$, and to compute
the differential of $\nabla\xi_s$ on $\Xi^\pm$, where it is constant.

\medskip\noindent \underline{Assertions (3) and (4).} Both assertions follow from Corollary
\ref{cor:Monge}.
\end{proof}

As a last application, we will also discuss the propagation of
the singular part with respect to $\Leb d$, which will be needed in the proof
of the main result in Theorem~\ref{th:Geod.Cvx}.

\begin{corollary}[Propagation of the singular part]
\label{cor:singular}
Let $\mu_0,\mu_1\in \MMM(\R^d)$ and let $(\mu_t)_{t\in [0,1]}$ be a geodesic
connecting $\mu_0$ to $\mu_1$ and let $\mu_s=\mu_s^\rma+\mu_s^\perp$ be the
decomposition of $\mu_s$ with respect to the Lebesgue measure $\Leb d$ at some point
$s\in (0,1)$.  For every $t\in [0,1]$ we set 
\begin{equation}
  \label{eq:170}
  \wt\mu_t:=(\bfT_{s\to t},q_{s\to t})_\star \mu_s^a \quad \text{and} \quad 
  \wh\mu_t:=(\bfT_{s\to t},q_{s\to t})_\star \mu_s^\perp.
\end{equation}
Then, the curves $(\wt\mu_t)_{t\in (0,1)}$ and $(\wh\mu_t)_{t\in (0,1)}$ are $\HK$
geodesics, we have $\wh\mu_t \perp \Leb d$ for $t\in [0,1]$ and
$\mu_t= \wt\mu_t + \wh\mu_t$ provides the Lebesgue decomposition for $t\in (0,1)$,
viz.\ $\mu^\rma_t=\wt\mu_t$ and $\mu_t^\perp = \wh\mu_t$.
\end{corollary}
\begin{proof}
Let us decompose $\Xi_s$ in the disjoint union of two Borel sets $A,B$ such
that $\mu_s^\rma=\mu_s\res A$ and $\mu_s^\perp=\mu_s\res B$ with $\Leb d(B)=0$.
By Theorem \ref{th:restriction} we clearly have
$\mu_t=\wt\mu_t + \wh\mu_t$. On the one hand, $\wt\mu_t \ll \Leb d$
by Corollary \ref{cor:start} for all $t\in (0,1)$. On the other hand, 
for all $t\in [0,1]$ the measure  $\wh\mu_t$ is concentrated on the set
$\bfT_{s\to t}(B)$ which is $\Leb d$-negligible, since $\bfT_{s\to t}$ is
Lipschitz. If follows that $\wh\mu_t \perp \Leb d$, so that $\wh\mu_t =
\mu_t^\rma$  and $\wh\mu_t=\mu_t^\perp$ for all $t\in (0,1)$. 

 The fact that $(\mu_t^a)$ and $(\mu_t^\perp)$ are geodesics
follows by Theorem \ref{th:restriction} as well.
\end{proof}

\subsection{Convexity of the Lebesgue density along \texorpdfstring{$\HK$}{HK}-geodesics}
\label{su:Representation}

In this subsection, we consider geodesics $(\mu_t)_{t\in
  [0,1]}$ such that 
$\mu_s\ll \Leb d$ for some, and thus for all,  $s\in (0,1)$.  We fix $s$ and introduce the functions $\alpha_s,\delta_s$ as in \eqref{eq:151} and the
functions
\begin{equation}
  \label{eq:153}
  \left.\begin{aligned}
    \gamma_s(t,x):={}&\alpha_s^{1/2}(t,x)=q_{s\to t}(x),\\
    \rho_s(t,x):={}&\alpha_s^{1/2}(t,x)\delta_s^{1/d}(t,x)=
    q_{s\to t}(x)\big( \! \det \rmD\bfT_{s\to t}(x)\big)^{1/d}    
  \end{aligned}
\right\} 
\quad \text{ for } x\in D_s.
\end{equation}

We now exploit the explicit differential relations for
$\gamma_s(t,x)=q_{s\to t}(x) $ and $\delta_s(t,x)= \det \rmD\bfT_{s\to t}(x)$
provided in Theorem~\ref{thm:characteristics} and derive lower estimates for
$\ddot \gamma_s$ and $\ddot \rho_s$. It remains unclear whether the given
choice for $\gamma_s$ and $\rho_s$ is the only possible, however it turns out
that for these variables the following curvature estimates are relatively
simple and hence the final convexity calculus goes through. 
For comparison, we mention that in the Kantorovich--Wasserstein case we have
$\gamma_\mathsf{KW}(t)\equiv 1$ and
$\rho_\mathsf{KW}(t)= \big(\delta_\mathsf{KW}(t)\big)^{1/d} $ with
$\delta_\mathsf{KW}(t)=\det((1{-}t)\bbI{+}t\rmD\bfT_\mathsf{KW}(x))$, such that
$\ddot\rho_\mathsf{KW}(t) \leq 0$ since $\rmD\bfT_\mathsf{KW}(x)$ is
diagonalizable with nonnegative real eigenvalues, see
\cite[Eqn.\,(9.3.12)]{AmGiSa08GFMS}.

\begin{proposition}[Curvature estimates for $(\rho,\gamma)$]
\label{pr:Est.alpha.delta}
Let $(\rho_s,\gamma_s):(0,1)\times \DDD_s\to [0,\infty[^2$ be defined as above
along a geodesic. Then, we have for all $t\in (0,1)$ the relations 
\begin{equation}
  \label{eq:Est.a.d}
  \frac{\ddot\gamma_s(t,x)}{\gamma_s(t,x)} \geq 0 
  \quad \text{and} \quad 
  \left\{ \ba{cl} \ds \frac{\ddot\rho_s(t)}{\rho_s(t)}\leq
        \Big(1{-}\frac4d\Big)\frac{\ddot\gamma_s(t)}{\gamma_s(t)} 
    & \text{for }d\geq 2, \\[0.8em] 
    \ds \frac{\ddot\rho_s(t)}{\rho_s(t)} =
        \Big(1{-}\frac4d\Big)\frac{\ddot\gamma_s(t)}{\gamma_s(t)} 
      & \text{for }d=1.  
  \ea\right.  
%
\end{equation}
\end{proposition}
\begin{proof} As $s \in (0,1)$ and $x\in \DDD_s$ are fixed, we will simply
write $\rho(t)$ instead of $\rho_s(t,x)$ and similarly for the other
variables. Using the specific definition of
$\rho$ we obtain 
\[
\frac{\ddot \rho}\rho = \frac{\ddot \gamma}\gamma + \frac2d
\,\frac{\dot\gamma}\gamma \, \frac{\dot \delta}\delta + \frac1d\,
\frac{\ddot\delta}\delta + \frac1d\big(\frac1d-1\big)
\big(\frac{\dot\delta}\delta\big)^2 .
\]
We can now use the formulas provided in \eqref{eq:157}--\eqref{eq:161} giving
$\dot\gamma = 2\xi_t \gamma$ and $\ddot \gamma = |\nabla \xi_t|^2 \gamma$,
where $\xi_t$ and its derivatives are evaluated at $y=\bfT_{s\to t}
(x)$. Inserting this and \eqref{eq:161} for $\dot\delta$ and $\ddot\delta$ into
the above relation for $\ddot\delta/\delta$ we observe significant
cancellations and obtain
\begin{equation}
  \label{eq:ddot.g.ddot.r}
  \frac{\ddot\gamma}{\gamma}= |\nabla\psxi_t|^2
\quad \text{and}\quad 
\frac{\ddot\rho}{\rho} = \frac1{d^2}\big( (\Delta
\psxi_t)^2 {-}d|\rmD^2\psxi_t|^2\big)+\big(1-\frac4d\big)|\nabla\psxi_t|^2 .
\end{equation}
For $d=1$ we have $\rmD^2\psxi = \Delta \psxi$, while for
$d\geq 2$ all matrices $A\in \R^{d\ti d}$ satisfy
$d|A|^2 = d \sum_{i,j=1}^d A_{ij}^2 \geq (\trace A)^2 = \big(\sum_{1}^d
A_{ii}\big)^2$. Thus, the curvature estimates \eqref{eq:Est.a.d} follow.
\end{proof}

The above curvature estimates will be crucial in Section
\ref{se:Geod.Conv} for deriving our main result on geodesic convexity. We
remark that for $d\geq 2$ they are even slightly better that the ``sufficient
curvature estimates'' given in \eqref{eq:Suff.Curv} because of
$1-4/d\leq 1-4/d^2$ (with equality only for $d=1$).

We finally derive a useful result concerning the convexity of the density
$t \mapsto c(t,x)$ along geodesics. This provides a direct proof of the fact,
which was used in \cite{DimChi20TGMH} that the $\rmL^\infty$-norm along
geodesics is bounded by the $\rmL^\infty$-norm of the two endpoints. Indeed, we
show more, namely that the function $t \mapsto c(t,\bfT_t(x))$ is either
trivially constant or it is strictly convex.

\begin{theorem}[Convexity of densities along geodesics]
\label{thm:StrictCvx}\mbox{ } 

\emph{(1)} Under the assumption of Corollary \ref{cor:start}, 
for every $s\in (0,1)$ and $x\in \DDD_s\cup \Xi^\pm$ 
the function $c_s(t)=c(t,\bfT_{s\to t}(x))$ given by
\eqref{eq:150} or \eqref{eq:169}, respectively, is 
\GGG convex and positive in $(0,1)$; 
moreover, with a possible $\Leb d$-negligible exception, 
it is either constant or strictly
convex. \EEE 

\emph{(2)} If moreover $\mu_0\ll\Leb d$
(resp.\,$\mu_1\ll\Leb d$) then for $\mu_s$-a.e.\,$x$ their limit as
$t\downarrow 0$ (resp.\:as $t\uparrow1$) coincides with
$c_0\circ \bfT_{s\to 0}$ (resp.\:$c_1\circ \bfT_{s\to 1}$).
\end{theorem}
\begin{proof}
  \noindent
\underline{Assertion (1).} Since $x\in \OOmega$ and $s\in (0,1)$ play no role, we drop them  for
notational simplicity. We simply calculate the second derivative of the
function $t\mapsto c(t)=\gamma(t)^{d+2}c_s/\rho(t)^{ d }$. 
If $c_s=c(s,x)=0$  then $c(t,\bfT_{s\to t}(x))=0$ and the result is
obviously true. Hence, we may assume $c_s>0$ and obtain after an
explicit calculation
\begin{equation}
  \label{eq:ddot.wtc}
    \ddot{c}=  c\,\Big((d{+}2) \frac{\ddot \gamma}\gamma - d
          \frac{\ddot\rho}\rho + (d{+}1)(d{+}2)
          \big(\frac{\dot\gamma}\gamma\big)^2 - 2d(d{+}2)
          \frac{\dot\gamma}\gamma \,\frac{\dot\rho}\rho + d(d{+}1)
      \big(\frac{\dot\rho}\rho\big)^2  \Big).
\end{equation}
The quadratic form involving the first derivatives is positive
definite, and for the terms involving the second derivatives we can
use the curvature estimates in \eqref{eq:Est.a.d} to obtain
\[
  \ddot{c} \geq c \,\Big( \big( (d{+}2) \frac{\ddot \gamma}\gamma - d\big(1
  -\frac4d\big) \frac{\ddot \gamma}\gamma + 0 \Big)= 6 \, c \, 
  \frac{\ddot \gamma}\gamma.
\]
\GGG 
Notice that $t\mapsto \gamma(t)$ is the square root of 
the non-negative (and strictly positive in $(0,1)$)
quadratic polynomial $\alpha(\cdot,x)$ given by 
\eqref{eq:151}, so that $\gamma''\ge0$ 
and we conclude that
$\ddot{c}(t)\ge 0$ as well due to $c(t)>0$.

Moreover,  if $x\not\in \Xi^0_s$
then $|\nabla\xi_s(x)|>0$,
have $\ddot \gamma(t)>0$, and we deduce that
$\ddot{c}(t)>0$ obtaining the strict convexity of $c$.

If $x\in \Xi^0_s$ where $\nabla\xi_s(x)=0$,
we can use the representation \eqref{eq:169} for $c$
up to a $\Leb d$-negligible set.
%
%

\noindent
\underline{Assertion (2).} If $\mu_0=c_0\Leb d\ll \Leb d$, then $\delta_s(0,x)>0$ for
$\mu_s$-a.a.\ $x\in \DDD_s$  thanks to
the last statement of Corollary \ref{cor:Monge}
(which is a direct consequence of 
Theorem~\ref{thm:regularity}(\ref{th:reg.label4}))
and both $\delta_s(0,x)$ and
$\alpha_s(0,x)$ coincides with their limit as $t\downarrow0$. A further application of 
Corollary \ref{cor:start}(3) yields the result.
The case $t=1$ is completely analogous.
\end{proof}

The above result easily provides the following statement on 
\GGG convexity of $\rmL^\infty$ norms along $\HK$-geodesics. \AAA This
generalizes to a corresponding result for the Kantorovich--Wasserstein geodesics 
(which might been known, but the authors were not able to identify a
reference, see the Remark \ref{rem:KWL} below). 
\EEE
\begin{corollary}[Convexity of the $L^\infty$ norm along geodesics]
  \label{cor:Linfty}
  Let $\mu_0,\mu_1\in \MMM(\R^d)$ be absolutely continuous with respect to $\Leb d$
  with densities $c_i\in L^\infty(\R^d)$ and let $(\mu_t)_{t\in [0,1]}$ be a
  $\HK$ geodesic connecting $\mu_0$ to $\mu_1$.  Then $\mu_t=c_t\Leb d$ and
  $\|c_t\|_{L^\infty}\le (1{-}t)\|c_0\|_{L^\infty}+ t\|c_1\|_{L^\infty}$.
\end{corollary}
\begin{proof} The result for $(\MMM(\R^d),\HK)$ follows directly from Theorem
\ref{thm:StrictCvx}. 
\end{proof}
\begin{remark}
\label{rem:KWL}
    \AAA 
    Let
  $(\mu^\rmW_t)_{t\in [0,1]}$ be the Kantorovich--Wasserstein geodesic connection
  between 
  two probability measures $\mu_0,\mu_1 \in \calP_2(\R^d)$ with
  $\mu_i=c_i \Leb d$ and $c_0,c_1 \in L^\infty(\R^d)$. Similar to the previous result, $\mu^\rmW_t=c^\rmW_t\Leb d$ is absolutely continuous w.r.t.~$\Leb d$ and
  $\|c^\rmW_t\|_{L^\infty}\le (1{-}t)\|c^\rmW_0\|_{L^\infty}+
  t\|c^\rmW_1\|_{L^\infty}$.

  In fact, for $(\calP_2(\R^d), \sfW_2)$ we replace \eqref{eq:150all} by the simpler
formula for the Kantorovich--Wasserstein transport
\[
c^\rmW(t,\bfT_{s\to t}(x) ) = \frac{c^\rmW_s(x)}{\delta_s(t,x)} \quad \text{with }
\delta_s(x) = \det \bfT^\rmW_{ s\to t}(x),
\] 
see \cite[Prop.\,9.3.9]{AmGiSa05GFMS}.
Using $\mu_0=c_0 \Leb d$ we can choose $s=0$ and have $\bfT^\rmW_{0\to t}(x) =
x + t(\nabla \varphi(x){-}x)$ for a convex Kantorovich potential. Since for
every symmetric positive semidefinite matrix $D$ the function $t \mapsto
1/\det\big( (1{-}t)I + tD\big)$ is convex, the desired result follows with the
same arguments as for Theorem \ref{thm:StrictCvx}. 

\end{remark}
\EEE

\section{\GGG Preliminary \EEE discussion of the convexity conditions}
\label{se:Discussion}

In this section, \AAA we discuss the equivalence of two formulations of the
convexity conditions and give a few examples. The proof of sufficiency and
necessity of these conditions is then given in the following Section
\ref{se:Geod.Conv}.

For most parts of this section, \EEE we assume that
$E:{[0,\infty[}\to \R\cup\{\infty\}$ is lower semi-continuous and convex,
satisfies $E(0)=0$, and is twice continuously differentiable on the interior of
its domain $D(E):=\set{c\geq 0}{E(c)<\infty}$.  The following result gives a
characterization of the conditions \eqref{eq:HK.cond.N} on
$N_E:(\rho,\gamma)\mapsto (\rho/\gamma)^d E(\gamma^{d+2}/\rho^d)$ in terms of
the derivatives of $E$, namely $\eps_j(c)=c^j E^{(j)}(c)$ for $j=0$, $1$, and
$2$, which appear in
\begin{equation}
  \label{eq:bbB.def}
  \bbB(c) := 
\bma{cc} \eps_2(c)-\frac{d-1}d\big(\eps_1(c){-}\eps_0(c)\big) & \eps_2(c)
-\frac12\big(\eps_1(c){-}\eps_0(c) \big) \\[0.2em] 
\eps_2(c)-\frac12\big(\eps_1(c){-}\eps_0(c) \big) & \eps_2(c) 
+ \frac12 \eps_1(c) 
\ema.
\end{equation}
This characterization 
will then be used to derive a nontrivial monotonicity result 
in Proposition \ref{pr:N.E.monot}, which is a crucial building block
of the main geodesic convexity result. 

Note that the variables $\rho$ and $\gamma$  are related 
to the variable $c$ via $c= c_0\gamma^{d+2}/\rho^d$.

\begin{proposition}[Equivalent conditions on $E$]
\label{pr:N.E.conds} 
Let $N_E$ and $\bbB$ be defined in terms of $E$ as in \eqref{eq:N.def}
and \eqref{eq:bbB.def}, respectively. Then the
following conditions are equivalent: 
\begin{enumerate}[{\upshape(A)}]
\item $N_E$ satisfies \eqref{eq:HK.cond.N};
\item in the interior of the domain $D(E)$ we have $\bbB(c)\geq 0$ 
  and $(d{-}1)\big(\eps_1(c){-}\eps_0(c)\big)\geq 0$.
\end{enumerate}
\end{proposition}
\begin{proof} We first observe that the desired monotonicity of
  $\rho \mapsto N_E(\rho,\gamma)$ for $d\geq 2$ is indeed equivalent
  to the condition $\eps_1(c)\geq \eps_0(c)$. This follows easily from the
  relation
\[
\pl_\rho N_E(\rho,\gamma) = \frac{d\rho^{d-1}}{\gamma^d}
\;\!E\big(\frac{\gamma^{d+2}}{\rho^d} \big)
+\frac{\rho^d}{\gamma^d}\;\!E'\big(\frac{\gamma^{d+2}}{\rho^d}
\big) \big({-}d\frac{\gamma^{d+2}}{\rho^{d+1}}\big)=
 \frac{d\rho^{d-1}}{\gamma^d}\;\!\big( \eps_0(c)- \eps_1(c)\big).
\]

It remains to establish the equivalence between the convexity of $N_E$
and the positive semi-definiteness of $\bbB$. For this we note that $N_E$
is given as a linear function of $E$, hence the Hessian $\rmD^2N_E$ will
be a given as a linear combination of $E$, $E'$, and $E''$. Indeed, an
explicit calculation yields 
\[
\rmD^2 N_E(\rho,\gamma) = \frac{\rho^d}{\gamma^d}\bma{cc}d/\rho&-d/\gamma\\ 0 &-2/\gamma \ema^\top
\bbB\big(\frac{\gamma^{d+2}}{\rho^d}\big)  \bma{cc}d/\rho&-d/\gamma\\ 0
&-2/\gamma \ema. 
\]
With this, we see that $\rmD^2 N_E$ is positive semidefinite if and only if
$\bbB$ is. Hence, the assertion is proved.   
\end{proof}

From the semi-definiteness of the matrix $\bbB(c)$, we obtain as
necessary conditions the non-negativity of the two diagonal elements
which provide the McCann condition $\bbB_{11}=\eps_2-
\frac{d{-}1}{d}(\eps_1-\eps_0) \geq 0$ and the convexity conditions with
respect to the Hellinger--Kakutani distance
$\bbB_{22}=\eps_2+\frac12\eps_1\geq 0$.  Moreover, testing $\bbB$ 
with $(1,-1)^\top$ reveals the additional condition
\begin{equation}
\label{eq:HK.cond3/2}
\binom{1}{-1}\cdot\bbB(c)\binom{1}{-1}\geq 0 \quad
\Longleftrightarrow\quad  (d{+}2)\eps_1(c)-2\eps_0(c)\geq 0.
\end{equation}

\begin{proposition}[New necessary monotonicity]
\label{pr:N.E.monot}
Let $E$ be such that the conditions in Proposition \ref{pr:N.E.conds}
hold and let $N_E$ be defined via \eqref{eq:N.def}. 
Then,  the following three equivalent conditions hold: 
\begin{enumerate}[\upshape (A)]
\itemsep-0.2em
\item The function ${]0,\infty[}\ni c\mapsto c^{-2/(d+2)} E(c)$
  is non-decreasing.
\item For all $\rho,\gamma>0$ we have the inequality $
    \big( 1- \frac4{d^2}\big) \rho \pl_\rho N_E(\rho,\gamma) + \gamma
    \pl_\gamma N_E(\rho,\gamma)\geq 0$. 
\item  For all $\rho,\gamma>0$ the mapping 
${]0,\infty[}\ni s \mapsto N_E(s^{1-4/d^2} \rho,s \gamma)$ is non-decreasing. 
\end{enumerate}

\end{proposition}
\begin{proof} 
Expressing $\pl_\rho N_E$ and $\pl_\gamma N_E$ via
$\eps_0$ and $\eps_1$ and using $\delta=(\rho/\gamma)^d$ we obtain 
\[
\rho\pl_\rho N_E(\rho,\gamma)= -d\delta(\eps_1{-}\eps_0) \quad
\text{and} \quad \gamma \pl_\gamma N_E(\rho,\gamma)=
\delta\big((d{+}2)\eps_1-d \eps_0\big).
\]
Thus, we conclude $(1{-} \frac4{d^2}) \rho \pl_\rho N_E(\rho,\gamma) + \gamma
\pl_\gamma N_E(\rho,\gamma)\:=\: \frac{2\delta}d\,\big((d{+}2)\eps_1 - 2\eps_0\big)$,
which is positive because of \eqref{eq:HK.cond3/2}. Thus,
(B) is established and the monotonicity of $s
\mapsto N_E(s^{1-4/d^2} \rho,s \gamma)$ in (C) follows simply by
differentiation.

Statement (A) follows by applying (C) for $\rho=\gamma=1$ and choosing $s =
c^{2(d{+}2)/d}$. 
\end{proof}

The crucial monotonicity stated at the end of the above proposition means 
\begin{equation}
  \label{eq:E.d.monotone}
  0\leq c_1 < c_2\quad \Longrightarrow \quad E(c_1) \,\leq
  \big(\frac{c_1}{c_2}\big)^{2/(d+2)} E(c_2).
\end{equation}
It implies that if $E$ attains a negative
value it cannot be differentiable at $c=0$: If $E(c_1)<0$ then $E(c) \leq
(c/c_1)^{2/(d+2)} E(c_1)<0$, which leads to $E'(c)\searrow -\infty$ for
$c \searrow 0$.
 
In the following examples we investigate which functions $E$ satisfy the above
conditions.  The following two results will be used in Corollary \ref{co:m.q}
to obtain geodesic convexity for functionals of the form
$E(c)=\int_\Omega a c^r\dd x$. The third example shows that in case of the
Boltzmann entropy with $E(c)=c\log c$ the conditions do not hold and hence
geodesic convexity fails.

\begin{example}[Density function $E(c)=c^m$]
\label{ex:m.ge.1} 
We have $\eps_0(c) = c^m$, $\eps_1(c) = m c^m$, and $\eps_2(c)=m(m{-1}) c^m$,
which gives the matrix
\[
\bbB(c) = c^m \bma{cc} (m{-}1)\big(m-\frac{d-1}{d}\big) &
(m{-}1)\big(m-\frac12 \big)\\[0.2em] (m{-}1)\big(m-\frac12 \big) & m \big(m-\frac12 \big)  \ema .
\]
The Hellinger condition $\bbB_{22}(c)\geq 0$ holds for $m \not\in {]0,\frac12[}$,
while the McCann condition $\bbB_{11}(c)\geq 0$ holds for $m \not\in
{]\frac{d-1}d,1[}$.  Moreover, for $d\geq 2$ the monotonicity condition 
$\eps_1\geq \eps_0$ implies $m\geq 1$. 

Thus, the remaining cases are either $m\geq 1$ or $d=1$ and $m\leq 0$, and it
remains to check $\det \bbB(c) \geq 0$. An explicit calculation gives
\[
\det \bbB(c)= (m{-}1)\big(m-\frac12 \big)\, \frac{(d{+}2)m -d}{2d}  .
\]
Clearly, for $m\geq 1$ we have $\det \bbB(c)\geq 0$ for all space dimensions
$d\in \N$. Moreover, $\det \bbB(c)< 0$ for $m\leq 0$. 

In summary, we obtain geodesic convexity if and only if $m\geq 1$. 
\end{example}

\begin{example}[Density function $E(c)=- c^q$]\label{ex:q.le.1}
As in the previous example we have 
\[
\bbB(c) = c^q \bma{cc} (1{-}q)\big(q-\frac{d-1}{d}\big) &
(1{-}q)\big(q-\frac12 \big)\\[0.2em] (1{-}q)\big(q-\frac12\big) & q
\big(\frac12-q \big)  \ema . 
\]
The Hellinger condition $\bbB_{22}(c)\geq 0$ holds for $q \in [0,\frac12]$,
while the McCann condition $\bbB_{11}(c)\geq 0$ holds for $q \in
[\frac{d-1}d,1]$, which also implies the monotonicity 
$\eps_1\geq \eps_0$. With
\[
\det \bbB(c)= (1{-}q) \,\big(\frac12-q \big)\, \frac{(d{+}2)q -d}{2d} \: \AAA c^{2q} .
\]
we obtain the additional condition $q\geq d/(d{+}2)$ and summarize that
$E(c)=-c^q$ leads to a geodesically convex functional if and only if $q \in
\big[\max\{\frac {d-1}d,\frac d{d+2}\}, \frac12\big]$,  which has solutions only for $d=1$ and
$d=2$. 
\end{example}

\begin{example}[Boltzmann entropy]\label{ex:BoltzNotGLCvx} As a negative example where the geodesic
convexity fails, we consider the Boltzmann function $E(c)= c\log c$. We
compute $ \bbB_{22}(c)= \eps_2(c)+\frac12\eps_1(c)= \frac32 c +\frac12c\log c$,
which shows that the necessary Hellinger condition fails.  Moreover,
considering the measures $\mu_0 = 0$ and $\mu_1=c\Leb d$ for a non-negative
density $c\in \rmL^1(\Omega)$ we find that along the geodesic curve, given by
$\mu(s)= s^2\mu_1$, we have
\[
\calE(\mu(s)) = \int_\Omega E(s^2 c)\dd x = s^2 \calE(\mu_1) +
2s^2\log(s)\int_\Omega c\dd x, 
\]
which is clearly not convex if $\int_\Omega c \dd x = \mu_1(\Omega)>0$.
\end{example}

Finally, we discuss a few examples where the density function $E$ is
not smooth. Note that the conditions in \eqref{eq:HK.cond.N} form a
closed cone. Moreover, as for convex functions, the supremum $ E : c
\mapsto \sup\set{\wt E_\alpha(c)}{\alpha \in A}$ satisfies
\eqref{eq:HK.cond.N} if all $\wt E_\alpha$ do so.

\begin{example}[Nonsmooth $E$]\label{ex:E.nonsmooth}
In applications one is also interested in cases where $E$ is
nonsmooth.
For example the case $E_\kappa(c)=\kappa c$ for $c\in
[0,c^*]$ and $E(c)=\infty$ for $c>c_*$ is considered in \cite{DimChi20TGMH}. Clearly, $E_0$
satisfies our assumptions \eqref{eq:HK.cond.N} since $N_E $ only takes
the values $0$ and $\infty$ and the value $0$ is taken on the convex
set $\gamma^{d+2}\leq \MAT c^*\EEE\rho^d$. Thus, $E_\kappa $ generates a functional
$\calE_\kappa = \calE_0 + \kappa \mathscr M $ that is geodesically $2\kappa$-convex. 

A second example is given by $E(c)=\max\{ 0, c^2-c\}$. We first
observe that $ \wt E_1(c)=c$ and $\wt E_2(c)=c^2$ satisfy
\eqref{eq:HK.cond.N}. Hence, $c\mapsto \max\{\wt E_1(c),\wt E_2(c)\}=
E(c)+c$ satisfies \eqref{eq:HK.cond.N} as well. Thus, we know that $E$
generates a functional $\calE$ that is at least geodesically
$({-}2)$-convex. However, we may inspect the function $c\mapsto c^2-c$
in the region $c\geq 1$ directly and find that $E$ itself satisfies
\eqref{eq:HK.cond.N}.    
\end{example}
 
In practical applications, in particular for evolutionary variational
inequalities as treated in \cite{LasMie22?EVIH}, it is desirable to find the
optimal $\lambda$ for the geodesic $\lambda$-convexity. So far, we have treated
the case of geodesic $0$-convexity and now return to the general case, which
leads to the conditions 
\[
\bbB(c) \geq \bma{cc}0&0\\ 0 &\lambda c/2 \ema \quad \text{ and } \quad 
(d{-}1)\big(\eps_1(c) -\eps_0(c)\big) \geq 0. 
\]  
The monotonicity condition is clearly independent of $\lambda$. The first
equation still relies on the necessary McCann condition $\bbB_{11}(c)\geq
0$. If this holds with strict inequality we see that the optimal $\lambda$ is
characterized by 
\begin{equation}
  \label{eq:lambda.opt}
  \lambda_\text{opt} = \inf\Bigset { \frac{2 \det \bbB(c)}{c \bbB_{11}(c)}}{ c>0}.
\end{equation}

\begin{example}[$d=1$ and $E(c)= c^2-c^{2/5}$] From the previous examples, we
  know that $E^2(c)=c^2$ and $E^{2/5}(c)=-c^{2/5}$ are both geodesically
  $0$-convex, and we want to show that the sum is geodesically $\lambda$-convex
  for $\lambda>0$. As $\bbB$ is linear in $E$ we have
  $\bbB(c)= c^2 \bbB^{(2)} + c^{2/5} \bbB^{(2/5)}$ with constant matrices
  $\bbB^{(2)}$ and $\bbB^{(2/5)} $ that are both strictly positive
  definite. Thus,
\[
\ell(c):=\frac{2 \det \bbB(c)}{c\, \bbB_{11}(c)} >0 \quad \text{ for all } c>0.
\]
Moreover, we find $\ell(c) \sim 2 c^{-3/5} \det \bbB^{(2/5)}/ \bbB^{(2/5)}_{11}$ for
$c \approx 0$ and $\ell(c) \sim 2 c \det \bbB^{(2)}/ \bbB^{(2)}_{11}$ for $
c \gg 1$. Thus, by compactness $\lambda_\mafo{opt} = \inf \bigset{\ell(c)}
{c>0}$ is strictly positive. 

Numerically, we find $\lambda_\mafo{opt} \approx 0.638 $ which is attained at $c_*\approx 0.0319$.
\end{example}

\begin{remark}[Geodesic convexity via the Otto calculus]\label{rm:GLC}
\slshape Following the key ideas 
in \cite{OttWes05ECCW,DanSav08ECDC} a formal calculus
for reaction-diffusion systems was developed in \cite{LieMie13GSGC}.
It uses the dynamical formulation in Subsection
\ref{subsec:HKdynamic} and the associated Onsager  operator
$\bbK(c)\xi = -\alpha\div(c\nabla\xi) + \beta c \xi$ to characterize
the geodesic $\lambda$-convexity of the functional $\calE$ by
calculating the quadratic form $M(c,\cdot)$ (contravariant Hessian of
$\calE$):
\begin{gather*}
M(c,\xi)=\langle \xi,\rmD\bfV(c) \bbK(c)\xi\rangle-
\frac12 \rmD_c\langle \xi , \bbK(c)\xi\rangle 
[\bfV(c)] \quad  \text{with } \bfV(c)=\bbK(c) \rmD\calE(c). 
\end{gather*}
Then, one needs to show the estimate $M(c,\xi) \geq \lambda \langle
\xi,\bbK(c)\xi\rangle $. 

Following the methods in \cite[Sect.\,4]{LieMie13GSGC}, for
 $c\in \rmC^0_\rmc(\Omega)$ and smooth $\xi$  we obtain 
\begin{align*}
M(c,\xi)&=\int_\Omega \bigg[ \alpha^2 \Big( \big( A(c){-}H(c)\big)
(\Delta \xi)^2 + H(c)\big| \rmD^2 \xi\big|^2 \Big) \\    
&\hspace{5em}+ \alpha\beta \Big(B_1(c) |\nabla\xi|^2 + B_2(c)\xi\Delta \xi \Big)
+\beta^2 B_3(c) \xi^2 \bigg] \dd x,\\
\text{where }&A(c)= \eps_2(c), \quad H(u)= \eps_1(c) - \eps_0(c), \quad
B_1(c)=\frac{3}{2} \eps_1(c)- \eps_0(c),\\
&B_2(c)=-2\eps_2(c) + \eps_1(c)- \eps_0(c),\quad B_3(c)= \eps_2(c)+
\frac{1}{2}\eps_1(c).
\end{align*}
  Analyzing the condition $M(c,\xi)\geq \lambda \langle \xi,
\bbK (c)\xi\rangle$ we find the conditions 
\begin{equation}
  \label{eq:E-GLCvx}
\forall c\geq 0:\quad (d{-}1)H(c)\geq 0, \ \ B_1(c)\geq
\frac{\lambda}{\beta} c, \ \ \bma{cc}\!\! A(c)-\frac{d-1}d H(c)& \frac12
B_2(c)\\ \frac12 B_2(c)& B_3(c)-\frac{\lambda}{\beta} c \!\!\ema \geq 0,
\end{equation}
which for $\lambda=0$ give the same conditions as $\bbB(c)\geq 0$, see
Proposition \ref{pr:N.E.conds}.  Note that the middle estimate in
\eqref{eq:E-GLCvx} follows from the first and the third estimates because of
\[
\textstyle
B_1(c) = \frac32\eps_1 - \eps_0 = \frac{d-1}d (\eps_1{-}\eps_0) + \frac1{2d}
\big((d{+}2)\eps_2{-}2\eps_0\big). 
\]
\end{remark}

\section{\texorpdfstring{\AAA Proof of \EEE geodesic  convexity of $\calE$}
            {Geodesic  convexity of E}} 
\label{se:Geod.Conv}
 
In this section, we finally prove the necessity and sufficiency of the
conditions for geodesic convexity of functionals $\calE$ on $\MMM(\Omega)$ in
\eqref{eq:HK.cond.N}, where we now allow for a general closed and convex domain
$\Omega\subset \R^d$.  In order to keep the arguments clear, we first 
restrict ourselves to absolutely continuous measures $\mu_0$ and $\mu_1$. Thus, by
Corollary~\ref{cor:start} the connecting geodesic curves are also absolutely
continuous, and we can rewrite $\calE$ along the latter in the form
\[
\calE(\mu_t)= \int_\Omega E(c(t,y))\dd y = \int_\Omega e(t,x) \dd x,
\quad \text{where }e(t,x)= \delta(t,x)\;\!E\Big( c_*(x)
\frac{\alpha(t,x)}{\delta(t,x)}\Big).
\]
The general case will then be treated by using an approximation argument.

Under the assumption that $E$ is twice differentiable in the interior of its
domain, we show that for $\mu_0$-a.a.\ $x\in\Omega$ the function
$t\mapsto e(t,x)$ is convex. Since $\alpha(\cdot,x)$ and $\delta(\cdot)$ are
analytic functions on $[0,1]$, we can show convexity in this case by
establishing $\ddot e(t,x)\geq 0$. For this, we can fix $x\in \Omega$, drop
the dependence on $x$ for notational convenience, and set
\begin{equation}
  \label{eq:e.rho.gamma}
  e(t)= \delta (t)E\big(c_*\frac{\alpha(t)}{\delta(t)}\big)=
N_E \big( \rho(t),\gamma(t) \big) \quad \text{with
}\rho:=(c_*\alpha)^{1/2}\delta^{1/d}, \ \  \gamma:=(c_*\alpha)^{1/2}, 
\end{equation}
and $N_E$ from \eqref{eq:N.def}. Now, the classical chain rule implies the 
relation
\begin{equation}
\label{eq:ddot-e}
\ddot e= \Big\langle \binom{\dot\rho}{\dot\gamma}, \rmD^2N_E(\rho,\gamma)
\binom{\dot\rho}{\dot\gamma} \Big\rangle +\pl_\rho
N_E(\rho,\gamma)\ddot\rho +  \pl_\gamma N_E(\rho,\gamma)\ddot\gamma. 
\end{equation}
The aim is to show $\ddot e(t)\geq 0$ for all $t\in [0,1]$. By the
convexity of $N_E$ it suffices to treat the last two terms. 

For this we exploit the curvature estimates \eqref{eq:Est.a.d} on $\ddot
\gamma$ and $\ddot \rho$ as well as the monotonicities in \eqref{eq:HK.cond.N.c} 
and Proposition \ref{pr:N.E.monot}.

\subsection{Usage of the curvature estimates}
\label{su:SuffCurv}

We first show that it is sufficient to use the curvature estimates
\begin{equation}
  \label{eq:Suff.Curv}
  \frac{\ddot\gamma}{\gamma}\geq 0 \quad \text{and} \quad 
  \left\{  \ba{cl}  \ds
    \frac{\ddot\rho}\rho \leq \big(1{-}\frac4{d^2}\big) 
    \!\;\frac{\ddot \gamma}\gamma &  \text{ for } d\geq 2,
    \\[0.9em]
    \ds \frac{\ddot\rho}\rho = \big(1{-}\frac4{d^2}\big)
    \!\;\frac{\ddot\gamma}\gamma &\text{ for } d=1.
    \ea\right.  
\end{equation}
In particular, the equality condition for $d=1$ is different from the
inequality conditions for $d\geq 2$.  This will be used to compensate for the
missing monotonicity of $N_E$ in \eqref{eq:HK.cond.N.c} in the case $d= 1$.

Below we will see that the curvature estimates \eqref{eq:Suff.Curv} are
\emph{necessary} to complete our proof. Note that they are implied by the
curvature estimates derived in Proposition \ref{pr:Est.alpha.delta}. In fact,
both coincide for $d=1$, while for $d\geq 2$ the former are strictly weaker as
the latter because of $1-4/d < 1-4/d^2$.

\begin{proposition}[$\ddot e\geq 0$ via curvature estimates]
\label{pr:SuffCurvE}
Assume that $N_E$ satisfies \eqref{eq:HK.cond.N} and that $t\mapsto
(\rho(t),\gamma(t))$ satisfies \eqref{eq:Suff.Curv}, then $\ddot e\geq 0$
in \eqref{eq:ddot-e}.
\end{proposition}
\begin{proof}
As the first term (involving $\rmD^2 N_E$) on the right-hand side of
\eqref{eq:ddot-e} is non-negative, we only have to show that the last
two terms have a non-negative sum. For this we rearrange terms as follows:
\begin{align*}
&\pl_\rho N_E(\rho,\gamma) \ddot \rho {+} 
 \pl_\gamma N_E(\rho,\gamma) \ddot \gamma 
=  \big(\!{-}\rho\pl_\rho N_E\big) 
   \Big(\! \big(1{-}\frac4{d^2}\big)
      \frac{\ddot \gamma}\gamma -\frac{\ddot \rho}\rho\Big) 
 +  \Big(\!\gamma \pl_\gamma N_E {+} \big(1{-}\frac4{d^2}\big)
    \rho  \pl_\rho N_E  \Big)\frac{\ddot\gamma}\gamma\,.
\end{align*} 
The right-hand side is the sum of two products, both of which are
non-negative. Indeed, the first product equals $0$ in the case $d=1$
independently of the sign of $\pl_\rho N_E$, because
the second factor is $0$. In the case $d\geq2$ both factors are
non-negative (using $\pl_\rho N_E\leq 0$ and the second curvature estimate in
\eqref{eq:Suff.Curv}), so the first product is non-negative again.

In the second product both terms are non-negative by Proposition
\ref{pr:N.E.monot}(B) and the first curvature estimate in
\eqref{eq:Suff.Curv}. Thus, $\ddot e\geq 0$ in \eqref{eq:ddot-e} is proved.
\end{proof}

\subsection{\texorpdfstring{The main results on geodesic $\lambda$-convexity}
                                     {The main results on geodesic lambda-convexity}}
\label{su:MainResults}

We are now ready to establish our main result on the geodesic convexity of
functionals $\calE$ given in terms of a density $E$. We now make our general
assumptions of $E$ precise.
\begin{subequations}
\label{eq:E.conds}
\begin{equation}
  \label{eq:E.conds.A}
\begin{aligned}
 & E:{[0,\infty[}\to \R\cup\{\infty\} \text{ is lower semi-continuous,
   convex,}\\
&E(0)=0, \text{ and there exists } c_\circ>0 \text{ such that } 
   E(c_\circ) < \infty.
\end{aligned} 
\end{equation}
We also want to include the case that $E$ is not necessarily superlinear, so we
introduce the recession constant
\[
E'_\infty:= \lim_{c\to \infty} \frac1c E(c) \in \R\cup\{\infty\}.
\]
The case $E'_\infty=\infty$ is the superlinear case where the functional
$\calE(\mu)$ is always $+\infty$, if $\mu$ has a singular part, i.e.\
$\mu^\perp\neq 0$ in the decomposition $\mu= c\Leb d + \mu^\perp$ with
$\mu^\perp \perp \Leb d$.

We introduce a closed (convex) domain
$\Omega\subset \R^d$, and 
we consider the set of measures $\mu$ with support contained in $\Omega$,
which we identify with $\MMM(\Omega)$.
In the case that the right derivative 
$\displaystyle E'_0:=\lim_{c\downarrow0}\frac 1c
E(c)$ of $E$ at $0$ is not finite, we 
further have to impose that $\Omega$ has finite Lebesgue measure. Therefore, we will assume that 
\begin{equation}
  \label{eq:E.conds.B}
  \begin{aligned}
 &
 \Omega \text{ is a closed convex set with nonempty
  interior and}\\
&  \text{$\Omega \text{ is also bounded, i.e.~}\Leb d(\Omega
 )<\infty$, if $E_0'=-\infty$.}
\end{aligned} 
\end{equation}
\end{subequations}
Thus, the functionals $\calE$ are defined as follows
\begin{equation} 
  \label{eq:calE.full}
 \hspace*{-2em}  \calE(\mu)=\int_{\Omega}  E(c(x)) \dd x +
E'_\infty \mu^\perp(\Omega)\quad\text{for}\quad \mu=c\Leb d{+}\mu^\perp \text{ with } 
  \mu^\perp \perp \Leb d. 
\end{equation}
It is well known that \eqref{eq:E.conds} guarantees that $\calE$ is a weakly
lower semi-continuous functional on $\MMM(\Omega)$.  In particular, condition
\eqref{eq:E.conds.B} is necessary to guarantee that the negative part
$x\mapsto \min\{ E(c(x)),0\}$ is integrable, because for $ c\in \rmL^1(\Omega)$
the functions $x\mapsto -\sqrt{c(x)}$ may not lie in $\rmL^1(\Omega)$. We refer
to Example \ref{ex:E.withPot} for a case where \eqref{eq:E.conds.B} can be
avoided by using a confining potential.

We are now in the position to formulate our main result on the geodesic
$\lambda$-convexity of integral functionals $\calE$ on the
Hellinger--Kantorovich space $(\MMM(\Omega),\HK)$.  The proof consists of three
steps. First, we assume that $E$ is twice continuously differentiable in its
domain. Restricting to geodesic curves connecting absolutely continuous
measures, we can use the above differentiable theory giving $\ddot e\geq 0$. In
Step~2, we generalize to possibly non-differentiable density functions $E$, but
keep absolutely continuous measures. For smoothing a given $E$, we use that
whenever $E$ solves the conditions \eqref{eq:HK.cond.N} and \eqref{eq:E.conds}
then $c\mapsto E(rc)$ does so for each $r\in [0,1]$. With a multiplicative
convolution we construct a smooth $E_\delta$ to which Step~1 applies. Finally,
Step 3 handles the case where $\mu_0^\perp$ or $\mu_1^\perp$ are non-zero by a
standard approximation argument of general measures using absolutely continuous
measures.
 
\begin{theorem}[Geodesic convexity of $\calE$] 
\label{th:Geod.Cvx} 
Assume that $E:{[0,\infty[}\to \R\cup\{\infty\}$ and $\Omega\subset \R^d$
satisfy \eqref{eq:E.conds.A} and \eqref{eq:E.conds.B}, respectively.  
If for a $\lambda_*\in \R$ the function
\[
N_{\lambda_*,E}(\rho,\gamma):= \big(\frac\rho\gamma\big)^d \,E\Big(
\frac{\gamma^{d+2}}{\rho^d} \Big) - \frac{\lambda_*}2 \gamma^2, \quad
\text{ for }\rho,\gamma>0, 
\]
satisfies the conditions \eqref{eq:HK.cond.N.b} and \eqref{eq:HK.cond.N.c},
then the functional $\calE$ defined in \eqref{eq:calE.full} is geodesically
$\lambda_*$-convex on $( \MMM(\Omega), \HK)$.
\end{theorem} 
\begin{proof}
  Without loss of generality, we set $\lambda_*=0$ throughout the proof and
  shortly write $N_E=N_{\lambda_*,E}$.
  
\emph{Step 1: The smooth and absolutely-continuous case.}  We first assume that $E$ is
twice continuously differentiable in the interior ${]0,c_E[}$ of its domain and
that the measures $\mu_0$ and $\mu_1$ are absolutely continuous with respect to
$\Leb d$, i.e.\ $\mu_j= c_j\Leb d$ for $c_j \in \rmL^1(\Omega)$.
 
We fix $s\in (0,1)$ adopting the notation of Corollary \ref{cor:start}.  Then,
the geodesic curve $t \mapsto \mu_t=c(t,\cdot)\Leb d$ satisfies
\[
  \calE(\mu_t)=\int_\Omega E(c(t,y))\dd y= \int_\Omega
  E\big(c_s(x)\frac{\alpha_s(t,x)}{\delta_s(t,x)} \big) \delta_s(t,x)\dd x =
  \int_\Omega e(t,x) \dd x
\]
with $e(t,x)=N_E(\rho(t,x),\gamma(t,x))$ as above.  We want to show that for
a.a.\ $x\in \Omega$ the function $t\mapsto e(t,x)$ is convex.

As shown in Theorem \ref{thm:StrictCvx} the functions $t\mapsto 
\wt c(t,x)=c_0(x)\frac{\alpha(t,x)}{\delta(t,x)}= \gamma(t,x)^{d+2}/\rho(t,x)^d$
are either constant or strictly convex.  If the function $\wt
c(\cdot,x)$ is constant then either $c_0(x)=0$ or
$(\rho(\cdot,x),\gamma(\cdot,x))$  is constant. In both cases,
$e(\cdot,x)$ is constant as well, and hence convex. 

In the strictly convex case, the values of $\wt c(t,x)$ for $t\in {]0,1[}$ lie
in the interior of the domain of $E$, where $E$ is twice differentiable. Hence,
combining Propositions \ref{pr:SuffCurvE} and \ref{pr:Est.alpha.delta} shows
that $t\mapsto e(t,x)$ is convex for a.a.\ $x\in \Omega$. Since integration
over $\Omega$ maintains convexity we conclude that $t \mapsto \calE(\mu_t)$ is
convex, too.

\emph{Step 2: The nonsmooth but absolutely-continuous case.}  We still assume
$\mu_j=c_j\Leb d$, but now consider an $E$ that is not necessarily twice
differentiable, but still satisfies \eqref{eq:E.conds}. We choose a function
$\chi\in \rmC^\infty_\rmc(\R)$ satisfying $\chi(r)\geq 0$,
$\int_{-2}^{-1} \chi(r)\dd r =1$ , and $\chi(r)=0$ for $r\not\in [-2,-1]$. Now
for $\delta \in {]0,1/2[}$ we define the smoothings
\[
E_\delta(c) = \int_0^1 \chi_\delta(r) E(rc) \dd r, \quad \text{where }
\chi_\delta(r)=\frac1\delta \chi\big(\frac1\delta(r{-}1)\big).
\]
Hence, $\chi_\delta$ has support in $[1{-}2\delta, 1{-}\delta]$. If the closure
of the domain of $E$ is $[0,c_E]$, then $E_\delta$ is well-defined and
$\rmC^\infty$ on ${]0,c_E/(1{-}\delta)[}$.  Moreover, for all $c\in [0,c_E]$ we
have $E_\delta(c)\to E(c)$ for $\delta \searrow 0$. We easily check, that
$E_\delta$ still satisfies the assumption \eqref{eq:HK.cond.N} and
\eqref{eq:E.conds}. Moreover, $E_\delta(c)$ can be estimated by $E(c)$ via
\begin{equation}
  \label{eq:Edelta.E}
  \exists\, K>0\ \forall\, \delta\in {]0,1/4[} \ \forall\,  c\geq 0: 
\qquad |E_\delta(c) | \leq K\big(c + |E(c)| \big). 
\end{equation}
To see this, we first consider the largest interval ${[0,c_1[}$ on which $E$ is
non-increasing. Then $0=E(0)\geq E_\delta(c) \geq E(c)$ which implies
\eqref{eq:Edelta.E} with $K=1$. If $c_1=\infty$ then we are done. If
$c_1<\infty$, then $E$ starts to increase and there exists
$c_2\in {[c_1,\infty[}$ with $E(c)\geq 0$ for $c\geq c_2$.  Using the
construction of $E_\delta$, we obtain for all $c\geq2c_2\geq c_2/(1{-}2\delta)$
the lower bound $E_\delta(c)\geq 0$. Using \eqref{eq:E.d.monotone} we easily
get $E_\delta(c) \leq E(c)$.

It remains to cover the case $c\in [c_1,3c_2]$. If $c_1=0$ then $E(c)\geq 0$
for all $c$, which means $c_2=0$ as well, then \eqref{eq:Edelta.E} follows
immediately from the above arguments. If $c_1>0$, a uniform continuity argument
gives the estimate $|E_\delta(c)-E(c)| \leq M$ for $c\in [c_1,3c_2]$.  Then,
choosing $K=M/c_1$ provides \eqref{eq:Edelta.E}.

With this preparation, Lebesgue's dominated convergence theorem implies
\[
 \mu=c\Leb d \text{ with }c\in \rmL^1(\Omega) \ \text{and} \ \calE(\mu)<\infty 
\quad \Longrightarrow \quad \int_\Omega E_\delta(c(x))\dd x \to
\calE(\mu) \text{ as }\delta \downarrow 0.
\]

Taking any constant-speed geodesic $[0,1] \ni t\mapsto \mu_t =
c(t,\cdot)\Leb d$, we know by Step~1 that the curves 
\[
\ol e_\delta: t \mapsto \int_\Omega E_\delta(c(t,x)) \dd x 
\]
are convex. As $\ol e_\delta(t)\to \calE(\mu_t)$ we conclude that
$t \mapsto \calE(\mu_t)$ is convex on $[0,1]$. 

\emph{Step 3: Pure growth.} The curve $t \mapsto t^2 \mu_1$ is the unique
geodesic connecting $\mu_0=0$ and $\mu_1$. Using the Lebesgue decomposition
$\mu_1= c_1\Leb d + \mu^\perp$ we see that
\[
t\mapsto \calE(\mu_t) = \calE(t^2c_1\Leb d) + t^2 E_\infty' \mu^\perp,
\]
is convex on [0,1] by Step 2 for the first term and by $E_\infty'\geq 0$. The
nonnegativity of $E_\infty'= \lim_{c\to \infty} E(c)/c$  follows from
\eqref{eq:E.conds.A} and Proposition 
\ref{pr:N.E.monot}(A), namely for $c\geq c_\circ$ we have 
\[
\frac1c E(c) = \frac1{c^{d/(d+2)}} \,c^{-2/(d+2)} E(c) 
\overset{\text{(i)}}\geq  \frac1{c^{d/(d+2)}} \,c_\circ^{-2/(d+2)} E(c_\circ) 
\to 0 \ \text{ for } c \to \infty. 
\]    

\emph{Step 4: The general case allowing for singular measures.} Singular
measures can only occur for $E$ with sublinear growth. Hence, we assume
$ E'_\infty \in \R$ from now on. In particular $\mathscr E$ is finite
everywhere, and using $E(c)\le E_\infty' c$ we have
$ \mathscr E(\mu) \leq E_\infty'\mu(\Omega)$.

As in Corollary \ref{cor:singular}, we consider an arbitrary geodesic
$(\mu_t)_{t\in [0,1]}$ connecting $\mu_0 $ and $\mu_1$. For a fixed $s \in (0,1)$, we 
decompose  $\mu_s$ as $\mu_s^\rma+\mu_s^\perp$. Then, $\mu_t= \wt\mu_t +
\wh\mu_t$ splits into two geodesics with disjoint supports and $\wt\mu_s=\mu_s^\rma$
and $\wh\mu_s= \mu_s^\perp$, see Corollary \ref{cor:singular}. 
Moreover, we have $\wh\mu_t \perp \Leb d$ and
$\wt\mu_t \ll \Leb d$ for all $t\in (0,1)$. This implies the relation
\[
  \mathscr E(\mu_t)=\mathscr E(\wt\mu_t)+\mathscr E(\wh\mu_t)=
  \mathscr E(\wt\mu_t)+E_\infty'\wh\mu_t (\Omega) .
\]
Since $(\wh\mu_t)_{t\in [0,1]}$ is a geodesic and the total mass functional
$\mathscr M(\mu)= \mu(\Omega)$ is convex (see \eqref{eq:calMcvx}) and
$E_\infty'\geq 0$, the last term $t\mapsto E_\infty'\wh\mu_t(\Omega)$ is
convex. Hence, it is sufficient to check the convexity of
$t\mapsto \calE(\wt\mu_t)$.

Since $\wt\mu_t\ll\Leb d$ for all $t\in (0,1)$, the function
$ t \mapsto \calE(\wt\mu_t) $ is convex in the open interval $(0,1)$ by
Step 2. Hence, to show convexity on $[0,1]$ it is sufficient to check that
\[
\limsup_{t\downarrow0}\mathscr E(\wt\mu_t) \leq \mathscr E(\wt\mu_0) \qquad
\text{and} \qquad 
\limsup_{t\uparrow1}\mathscr E(\wt\mu_t)\le \mathscr E(\wt\mu_1),
\] 
because $\HK$ convergence implies weak convergence and $\calE$ is weakly
l.s.c.

Let us focus on the limit $t\downarrow 0$ as the limit $t\uparrow 1$
is completely analogous. The problem is that $\wt\mu_t \ll \Leb d$ for
$t\in(0,1)$ only, but $\wt\mu_0$ may have a singular part. Hence, we forget the
decomposition $\mu_t=\wt\mu_t+\wh\mu_t$ and use a different one. Before that, we
restrict to the case $\mu_0(\Xi^+)=0$ because on $\Xi^+$ we have pure
growth and this case is covered by Step~3. 

Now, we exploit the Lebesgue decomposition of
$\mu_0=\mu_0^\rma+\mu_0^\perp$ at $t=0$ and 
consider two disjoint Borel sets $A,B \subset \Omega\setminus \Xi^+$ such
that $\mu_0^\rma=\mu_0\res A$ and $\mu_0^\perp=\mu_0\res B$.
We define the
corresponding disjoints sets $A_t:=\bfT_{t\to0}^{-1}(A)$ and
$B_t:=\bfT_{t\to 0}^{-1} (B)$ as well as the measures
$\nu_t^A:=\mu_t\res A_t$ and $\nu_t^B:=\mu_t\res B_t$.  By Theorem
\ref{th:restriction}, we obtain two 
geodesics $\nu_t^A$, $\nu_t^B$ concentrated on disjoint sets giving
$\mathscr E(\mu_t)= \mathscr E(\nu_t^A)+\mathscr E(\nu_t^B)$.  Since
$\nu_t^A\ll\Leb d$ for every $t\in [0,1)$ we deduce that
$t\mapsto \mathscr E(\nu_t^A)$ is convex up to $0$ by Step 2.  
Concerning $\mathscr E(\nu_t^B)$, we use $\calE(\mu)\leq E'_\infty\mu(\Omega)$
and find
\[
  \limsup_{t\downarrow0}\mathscr E(\nu_t^B)
  \leq E_\infty'\limsup_{t\downarrow0}\nu_t^B(\Omega) =
  E_\infty'\nu_0^B(\Omega) =\mathscr E(\nu_0^B),
\]
where we exploited $\nu_0^B\perp\Leb d$ in the last identity. 

This finishes the proof of the main theorem
\end{proof}

The following result is a direct consequence of the main result by 
using the results of Examples \ref{ex:m.ge.1} and \ref{ex:q.le.1},
respectively. In particular, this establishes the result announced in
\cite[Thm.\,2.14]{DimChi20TGMH}. 

\begin{corollary}[Power-law functionals] \label{co:m.q}
Assume that $\Omega\subset \R^d$ and $E:{[0,\infty[}\to \R$
satisfy \eqref{eq:E.conds} and let $\calE$ be defined via
\eqref{eq:calE.full}.   
\begin{enumerate}[{\upshape(1)}]
\item If $E(c)=c^m$ with $m\geq 1$, then $\calE$ is geodesically convex on
  $( \MMM(\Omega), \HK)$.
\item If $\calL^{d}(\Omega)<\infty$, $d\in \{1,2\}$, and $E(c)=-c^q$ with
  $d/(d{+}2) \leq q \leq 1/2$, then $\calE$ is geodesically convex on
  $( \MMM(\Omega), \HK)$.
\end{enumerate} 
\end{corollary}

\begin{example}\label{ex:E.withPot}\slshape We have seen above that the
  density $E(c)=-\sqrt{c}$ produces a geodesically convex functional
  in dimensions $d=1$ and $2$, if $\Leb d(\Omega)<\infty$. 
  The restriction of finite volume for $\Omega$ can be dropped
  by using a confining potential $V$ as follows: Let
\[
\calE_{1/2,V}(\mu) = \int_{\R^d} \big({-}\sqrt{c(x)}\big) \dd x + \int_{\R^d}
V \dd \mu \ \text{ for } \mu=c\Leb d + \mu^\perp,
\]
where $V\in \rmC(\R^d)$ satisfies for $m >d$ and $A\in \R$ the
lower bound $V(x)\geq a_0|x|^m -A$ on $\R^d$. Then it is easy to see that
$\calE_{1/2,V}$ is well-defined and weakly lower semi-continuous. 

Moreover, in \cite[Prop.\,20]{LiMiSa16OTCR} it was shown for a
continuous $V:\R^d \to \R$ with $\inf  V> - \infty$ that  the
linear mapping $ \mu \mapsto \int_{\R^d} V \dd \mu$ is geodesically
$\lambda_V$-convex on 
$(\MMM(\Omega),\HK)$ if and only if the mapping $\wt V:[x,r]\mapsto r^2 V(x)$
is geodesically $\lambda_V$-convex on the metric cone space $(\mfC,\msd_\mfC)$. 
For smooth $V$, this amounts to the estimate 
\[
\bma{cc}
\nabla ^2 V(x)+2V(x)\bbI_d&\nabla V(x)\\[0.2em]
\nabla V(x)^\top& 2V(x) \ema \geq \lambda_V\bbI_{d+1}.
\] 
Thus, for $V$ satisfying both of the above assumptions, the functional
$\calE_{1/2,V}$ is geodesically $\lambda_V$-convex on $(\MMM(\R^2), \HK)$ for
$d\in \{1,2\}$. For $d=1$ we may choose $V(x)= \alpha + \beta |x|^2$ with
$\beta>0$ and obtain $\lambda_V=2\alpha$. 
\end{example}

\subsection{Necessity of the conditions on \texorpdfstring{$E$}{E}}
\label{su:Necessity}

Theorem~\ref{th:Geod.Cvx} states that the conditions \eqref{eq:HK.cond.N} and
\eqref{eq:E.conds} on the density $E:[0,\infty) \to (-\infty,+\infty]$ 
are sufficient for the geodesic convexity of the integral functional $\calE$. 
We finally show that the conditions are also sufficient. To simplify the
analysis we restrict ourselves to the smooth case where $E:\mafo{dom}(E)\to \R$
lies in $\rmC^2$. Thus, we can obtain conditions by differentiation along
suitably chosen geodesic curves. For this, the characteristic equations
\eqref{eq:CharSystem} derived in Theorem~\ref{thm:characteristics} will be the main tool. 

\begin{theorem}[Necessity of conditions on $E$]
\label{th:Necessity} 
Consider  a closed,
convex domain $\Omega\subset \R^d$ with nonempty interior and a density
function $E:[0,\infty) \to (-\infty,+\infty]$ such 
that \eqref{eq:E.conds} holds and that 
$E$ is $\rmC^2$ on the interior of its domain. If 
the induced functional $\calE:\calM(\Omega) \to  (-\infty,+\infty]$ defined in
\eqref{eq:calE.full} is geodesically convex on $(\calM(\Omega),\HK)$, then $E$
satisfies the conditions \eqref{eq:HK.cond.N}. 
\end{theorem}
\begin{proof} We first observe that it is sufficient to show that for
$(\rho_*,\gamma_*)$ with $c_*=\gamma_*^{2+d}/\rho_*^d\in \Int(\mafo{dom}(E))$ we
have the inequalities
\begin{equation}
  \label{eq:N.Inequalities}
  \rmD^2 N_E(\rho_*,\gamma_*)\geq 0 \quad \text{and} \quad (d{-}1) \pl_\rho
  N_E(\rho_*,\gamma_*)\leq 0. 
\end{equation}
By the scaling properties of $N_E(\rho,\gamma)=(\rho/\gamma)^d
E(\gamma^{2+d}/\rho^d)$ it is sufficient to look at the case
$(\rho_*,\gamma_*)=(c_*^{1/2},c_*^{1/2})$. 

The main idea is to construct suitable geodesic curves $\mu_t$ such that the
convexity of $t\mapsto \calE(\mu_t)$ gives the desired inequality. For this we
choose a point $x_*\in \Int(\Omega)$ and $r_*>0$ such that 
$B_{3r_*}(x_*)\subset \Omega$. Without loss of generality we assume $x_*=0$ and
write $B_r$ in place of $B_r(x_*)$ for $r\in (0,3r_*]$. 

We further choose an $s \in (0,1)$ and a smooth function $\xi_s \in
\rmC^3(B_{3r_*})$. Then, there exists an $\eps>0$ such that there is a unique smooth
solution $\xi:(s{-}\eps, s{+}\eps)\ti B_{2r_*}\to \R$ of the Hamilton--Jacobi
equation \eqref{eq:genHJeqn}. With this $\xi_t= \xi(t,\cdot)$ and $r\in
(0,r_*)$ we can construct a geodesic curve 
\[
(s{-}\eps, s{+}\eps) \ni t\mapsto \mu_t^{(r)} \quad \text{with } \mu_s^{(r)}=c_* \Leb
d\res B_r \ \text{ and } \ \mafo{sppt}(\mu_t^{(r)}) \subset B_{2r_*}
\]
of absolutely continuous measures $\mu_t =c^{(r)}(t,\cdot) \Leb d$, see
\eqref{eq:150} in Corollary \ref{cor:start}. If necessary $\eps$ needs to be
reduced to avoid mass flowing outside $B_{2r_*}$. For this geodesic we have
\[
\calE(\mu_t^{(r)})=\int_\Omega E(c^{(r)}(t,y))\dd y = \int_{B_r} e(t,x) \dd x \quad
\text{with } e(t,x)= N_E(\rho(t,x),\gamma(t,x)),
\]
where $\gamma(t,x)=c_*^{1/2}q_{s\to t}(x)$ and $\rho(t,x)= c_*^{1/2}q_{s\to
  t}(x) \big(\delta_{s\to t}(x)\big){}^{1/d}$ with $q$ and $\delta$ from
\eqref{eq:CharSystem}. Note that $q$, $\delta$, and $e$ do not depend on $r$,
cf.\ Theorem \ref{th:restriction}. 

By the smoothness of $\xi$, and hence of $\rho$ and $\gamma$, we may pass to
the limit $r\downarrow 0$ in the convex functions 
$t \mapsto \frac1{\Leb d(B_r)} \calE(\mu_t^{(r)})= \frac1{\Leb d(B_r)}
\int_{B_r} e(t,x) \dd x$. Thus, the limit 
\[
t \mapsto e(t,0) = N_E(\rho(t,0),\gamma(t,0)) \quad \text{is convex on } (s{-}\eps,
s{+}\eps).
\] 
In particular, the second derivative is non-negative which means 
\[
0 \leq \ddot e(s,0) = \rmD^2
N_E(c_*^{1/2},c_*^{1/2})\Big[\text{\smaller$\ds\binom{\dot\rho}{\dot\gamma}$},
\text{\smaller$\ds\binom{\dot\rho}{\dot\gamma}$} \Big] 
 + \pl_\rho N_E(c_*^{1/2},c_*^{1/2} ) \ddot\rho +
\pl_\gamma N_E(c_*^{1/2},c_*^{1/2} ) \ddot\gamma,
\]
where now $\dot\rho=\pl_t\rho(s,0)$, $\dot\gamma$, $\ddot\rho$, and $\ddot
\gamma$ are
given by \eqref{eq:CharSystem} and \eqref{eq:ddot.g.ddot.r} in terms of $\xi$ only:
\[
\frac{\dot\rho}{c_*^{1/2}} =2\xi_s{+}\Delta \xi_2, \quad 
\frac{\dot\gamma}{c_*^{1/2}}= 2\xi_s, \quad 
\frac{\ddot\rho}{c_*^{1/2}} =\frac{(\Delta\xi_s)^2 {-}
d|\rmD^2\xi_s|^2}{d^2} + \frac{d{-}4}d |\nabla \xi_s|^2,\quad 
\frac{\ddot\gamma}{c_*^{1/2}}=|\nabla \xi_s|^2,
\]
where $\xi_s$ and its derivatives are evaluated at $x=x_*=0$. 

To obtain the convexity of $N_E$ we can now choose the functions $\xi_s$ such
that $\ddot\rho=\ddot\gamma=0$, which is the case for $\xi_s(x) = \alpha +
\beta |x|^2$, which implies $\nabla \xi_s(0)=0 $, $\Delta \xi_s(0)= 2d\beta$,
and $|\rmD^2\xi_s|^2=d\beta^2$. Moreover, $\dot\rho$ and $\dot\gamma$ can be
chosen arbitrarily by adjusting $\alpha,\,\beta\in \R$. Thus, $\rmD^2N\geq 0$
is established.

To prove the second estimate in \eqref{eq:N.Inequalities} we may assume $d\geq
2$, as there is nothing to show for $d=1$.  Choosing the function
$\xi_s=\alpha (x_1^2-x_2^2)$ we obtain $\dot\rho=\dot\gamma=\ddot\gamma=0$ and 
$\ddot\rho= -8c_*^{1/2} \alpha^2/d$. This implies $\pl_\rho N_E\leq 0$ and the
theorem is established.   
\end{proof}

\subsection{A more direct sufficiency proof for \texorpdfstring{$2\leq d\leq 4$}{2 <= d <= 4}}  
\label{su:MoreDirect}

The above proof of Theorem~\ref{th:Geod.Cvx} strongly relies on differentiating
 $e(t,x)=N_E(\rho(t,x),\gamma(t,x))$ with respect to $t$. In the case
$2\leq d\leq 4$, this can be avoided since we have the curvature estimates
  \begin{equation}
    \label{eq:Curvature.d4}
  \text{(a) }  t\mapsto \gamma(t,x) \text{ is convex \quad and \quad
    (b) } t\mapsto
    \rho(t,x) \text{ is concave}, 
  \end{equation}
where we used $d\leq 4$ in \eqref{eq:Est.a.d} for (b).
With (a) and (b),  we can further exploit
\begin{enumerate}[{\upshape (i)}]
  \item  the convexity of $N_E$,
\item the
monotonicity of $\rho \mapsto N_E(\rho,\gamma)$ (non-increasing, cf.\
\eqref{eq:HK.cond.N.c} for $d\geq 2$), and 
\item the monotonicity of $s \mapsto N_E(s^{1-4/d^2}\rho,s\gamma)$
(non-decreasing, cf.\ Proposition \ref{pr:N.E.monot}).
\end{enumerate}
Choosing $t_0,t_1,\theta\in [0,1]$ with $t_0<t_1$, we set
$t_\theta:=(1{-}\theta)t_0+\theta t_1$ and have to show 
\begin{equation}
  \label{eq:DirectProof1}
N_E(\rho(t_\theta),\gamma(t_\theta)) \leq (1{-}\theta)
N_E(\rho(t_0),\gamma(t_0)) + \theta N_E(\rho(t_1),\gamma(t_1)). 
\end{equation}
We start with the right-hand side and use convexity (i) first:
\begin{align*}
&(1{-}\theta) N_E(\rho(t_0),\gamma(t_0)) + \theta N_E(\rho(t_1),\gamma(t_1))
\overset{\text{(i)}}\geq  
 N_E\big((1{-}\theta)\rho(t_0){+}\theta \rho(t_1), 
       (1{-}\theta)\gamma(t_0){+}\theta \gamma(t_1) \big).
\end{align*}
With the convexity (a) of $\gamma$ we have 
  $s:=\gamma(t_\theta)/\big[(1{-}\theta)\gamma(t_0){+}\theta
  \gamma(t_1)\big] \in [0,1]$ and continue
\begin{multline*}
N_E\big((1{-}\theta)\rho(t_0){+}\theta \rho(t_1), 
       (1{-}\theta)\gamma(t_0){+}\theta \gamma(t_1) \big)\\
\overset{\text{(iii)}}\geq N_E\big( s^{1-4/d^2}
  \big[(1{-}\theta)\rho(t_0){+}\theta \rho(t_1)\big] ,
  s\big[(1{-}\theta)\gamma(t_0){+}\theta
  \gamma(t_1)\big] \big)\\
 = N_E\big( s^{1-4/d^2}
  \big[(1{-}\theta)\rho(t_0){+}\theta \rho(t_1)\big] ,\gamma(t_\theta)\big).
\end{multline*} 
Using the monotonicity (ii) (for $d\geq 2$) we can increase the first argument
using $s^{1-4/d^2} \leq 1$ (because of $s\in [0,1]$ and $d\leq 4$) 
and then exploit the concavity in (b) of $\rho$ (i.e.\
$\rho(t_\theta)\geq (1{-}\theta)\rho(t_0){+}\theta \rho(t_1)$) giving
\begin{align*}  
N_E\big( s^{1-4/d^2}
  \big[(1{-}\theta)\rho(t_0){+}\theta \rho(t_1)\big] ,
  \gamma(t_\theta) \big)
 \overset{\text{(ii), (a+b)}}\geq N_E(\rho(t_\theta),\gamma(t_\theta)).
\end{align*}
Thus, we have proved the desired convexity \eqref{eq:DirectProof1} for the case $d\in
\{2,3,4\}$.

\subsection*{Acknowledgments}
M.L.\ was partially supported by the German Research Foundation (DFG) under the
German National Excellence Strategy – The Berlin Mathematics Research Center
MATH+ (EXC-2046) in project AA2-10.

A.M.\ was partially supported by 
DFG via subproject C05 within CRC\,1114 \emph{Scaling Cascades in Complex
Systems} (Project no.\;235221301) and by ERC AdG\,267802 \emph{AnaMultiScale}.
A.M. is grateful to the Institute of Advanced Study, Technische
Universit\"at M\"unchen, for hosting him in the final stage of this work. 

G.S. was partially supported by the Institute for Advanced Study, Technische
Universit\"at M\"unchen, funded by the German Excellence Initiative, and by the
MIUR-PRIN 2017 project \emph{Gradient flows, Optimal Transport and Metric
  Measure Structures.}

\footnotesize

\newcommand{\etalchar}[1]{$^{#1}$}
\def\cprime{$'$}
\providecommand{\bysame}{\leavevmode\hbox to3em{\hrulefill}\thinspace}

\end{document}